\documentclass[twocolumn,natbib]{svjour3}          
\smartqed  
\usepackage{graphicx}
\usepackage{amsmath,amssymb}
\usepackage{subeqnarray}
\usepackage{cases}
\usepackage{multirow}
\usepackage{amssymb}
\usepackage{subfigure}
\usepackage{graphicx}
\usepackage{epsfig,times}
\usepackage{algorithm,algorithmic}

\newcommand{\rank}{\mbox{rank}}
\newcommand{\comment}[1]{{}}

\begin{document}

\title{Solving Principal Component Pursuit in Linear Time via $l_1$ Filtering
}


\author{Risheng Liu         \and
        Zhouchen Lin        \and
        Siming Wei          \and
        Zhixun Su 
}


\institute{R. Liu \at
              School of Mathematical Sciences, Dalian University of Technology.\\
              \email{rsliu0705@gmail.com}           
           \and
           Z. Lin (corresponding author)\at
              Key Lab. of Machine Perception (MOE), Peking University.\\
              This work was done in Microsoft Research Asia.\\
              \email{zlin@pku.edu.cn}           
           \and
           S. Wei \at
           College of Computer Science and Technology, Zhejiang University.\\
           \email{tobiawsm@gmail.com}           
           \and
           Z. Su \at
           School of Mathematical Sciences, Dalian University of Technology.\\
           \email{zxsu@dlut.edu.cn}           
}

\date{Received: date / Accepted: date}

\maketitle

\begin{abstract}
In the past decades, exactly recovering the intrinsic data
structure from corrupted observations, which is known as robust
principal component analysis (RPCA), has attracted tremendous
interests and found many applications in computer vision.
Recently, this problem has been formulated as recovering a
low-rank component and a sparse component from the observed data
matrix. It is proved that under some suitable conditions, this
problem can be exactly solved by principal component pursuit
(PCP), i.e., minimizing a combination of nuclear norm and $l_1$
norm. Most of the existing methods for solving PCP require
singular value decompositions (SVD) of the data matrix, resulting
in a high computational complexity, hence preventing the
applications of RPCA to very large scale computer vision problems.
In this paper, we propose a novel algorithm, called $l_1$
filtering, for \emph{exactly} solving PCP with an $O(r^2(m+n))$
complexity, where $m\times n$ is the size of data matrix and $r$
is the rank of the matrix to recover, which is supposed to be much
smaller than $m$ and $n$. Moreover, $l_1$ filtering is
\emph{highly parallelizable}. It is the first algorithm that can
\emph{exactly} solve a nuclear norm minimization problem in
\emph{linear time} (with respect to the data size). Experiments on
both synthetic data and real applications testify to the great
advantage of $l_1$ filtering in speed over state-of-the-art
algorithms.
\end{abstract}


\section{Introduction}

Robustly recovering the intrinsic low-dimensional structure of
high-dimensional visual data, which is known as robust principal
component analysis (RPCA), plays a fundamental role in various
computer vision tasks, such as face image alignment and
processing, video denoising, structure from motion, background
modeling, photometric stereo and texture representation (see e.g.,
\cite{Wright-NIPS2009}, \cite{Hui2010Video},
\cite{fernando2003rsl}, \cite{peng2010rasl}, \cite{wu-robust}, and
\cite{zhang2011texture}, to name just a few). Through the years, a
large number of approaches have been proposed for solving this
problem. The representative works include \cite{fernando2003rsl},
\cite{Nie2011L1}, \cite{Aanes-2002-fac},
\cite{Baccini96anl1-norm}, \cite{Ke2005l1},
\cite{Skocaj-2007-weight}, and \cite{Markus2009FRPCA}. The main
limitation of above mentioned methods is that there is no
theoretical guarantee for their performance. Recently, the
advances in compressive sensing have led to increasingly interests
in considering RPCA as a problem of exactly recovering a low-rank
matrix $\mathbf{L}_0$ from corrupted observations
$\mathbf{M}=\mathbf{L}_0 + \mathbf{S}_0$, where $\mathbf{S}_0$ is
known to be sparse (\cite{Wright-NIPS2009},
\cite{candes2009robust}). Its mathematical model is as follows:
\begin{equation}
\min \rank(\mathbf{L}) + \lambda \|\mathbf{S}\|_{l_0},\quad s.t.
\quad \mathbf{M}=\mathbf{L}+\mathbf{S},\label{eq:RPCA}
\end{equation}
where $\|\cdot\|_{l_0}$ is the $l_0$ norm of a matrix, i.e., the number of nonzero entries in the matrix.

Unfortunately, problem (\ref{eq:RPCA}) is known to be NP-hard. So
\cite{candes2009robust} proposed using principal component pursuit
(PCP) to solve (\ref{eq:RPCA}), which is to replace the rank
function and $l_0$ norm with the nuclear norm (which is the sum of
the singular values of a matrix, denoted as $\|\cdot\|_*$) and
$l_1$ norm (which is the sum of the absolute values of the
entries), respectively. More specifically, PCP is to solve the
following convex problem instead:
\begin{equation}
\min \|\mathbf{L}\|_* + \lambda \|\mathbf{S}\|_{l_1},\quad s.t. \quad \mathbf{M}=\mathbf{L}+\mathbf{S}.\label{eq:PCP}
\end{equation}
They also \emph{rigorously proved} that under fairly general
conditions and $\lambda=1/\sqrt{\max(m,n)}$, PCP can \emph{exactly
recover} the low-rank matrix $\mathbf{L}_0$ (namely the underlying
low-dimensional structure) with an overwhelming probability, i.e.,
the difference of the probability from 1 decays exponentially when
the matrix size increases. This theoretical analysis makes PCP
distinct from previous methods for RPCA.

All the existing algorithms for RPCA need to compute either SVD or
matrix-matrix multiplications on the whole data matrix. So their
computation complexities are all at least quadratic w.r.t. the
data size, preventing the applications of RPCA to large-scale
problems when the time is critical. In this paper, we address the
large-scale RPCA problem and propose a truly linear cost method to
solve the PCP model (\ref{eq:PCP}) when the data size is very
large while the target rank is relatively small. Such kind of data
is ubiquitous in computer vision.

\subsection{Main Idea}\label{sec:outline}
Our algorithm fully utilizes the properties of low-rankness. The
main idea is to apply PCP to a randomly selected submatrix of the
original noisy matrix and compute a low rank submatrix. Using this
low rank submatrix, the true low rank matrix can be estimated
efficiently, where the low rank submatrix is part of it.

\begin{figure}[t]
\centering
\begin{tabular}{c@{\extracolsep{0.2em}}c@{\extracolsep{0.2em}}c@{\extracolsep{0.2em}}c}
\includegraphics[width=0.11\textwidth,
keepaspectratio]{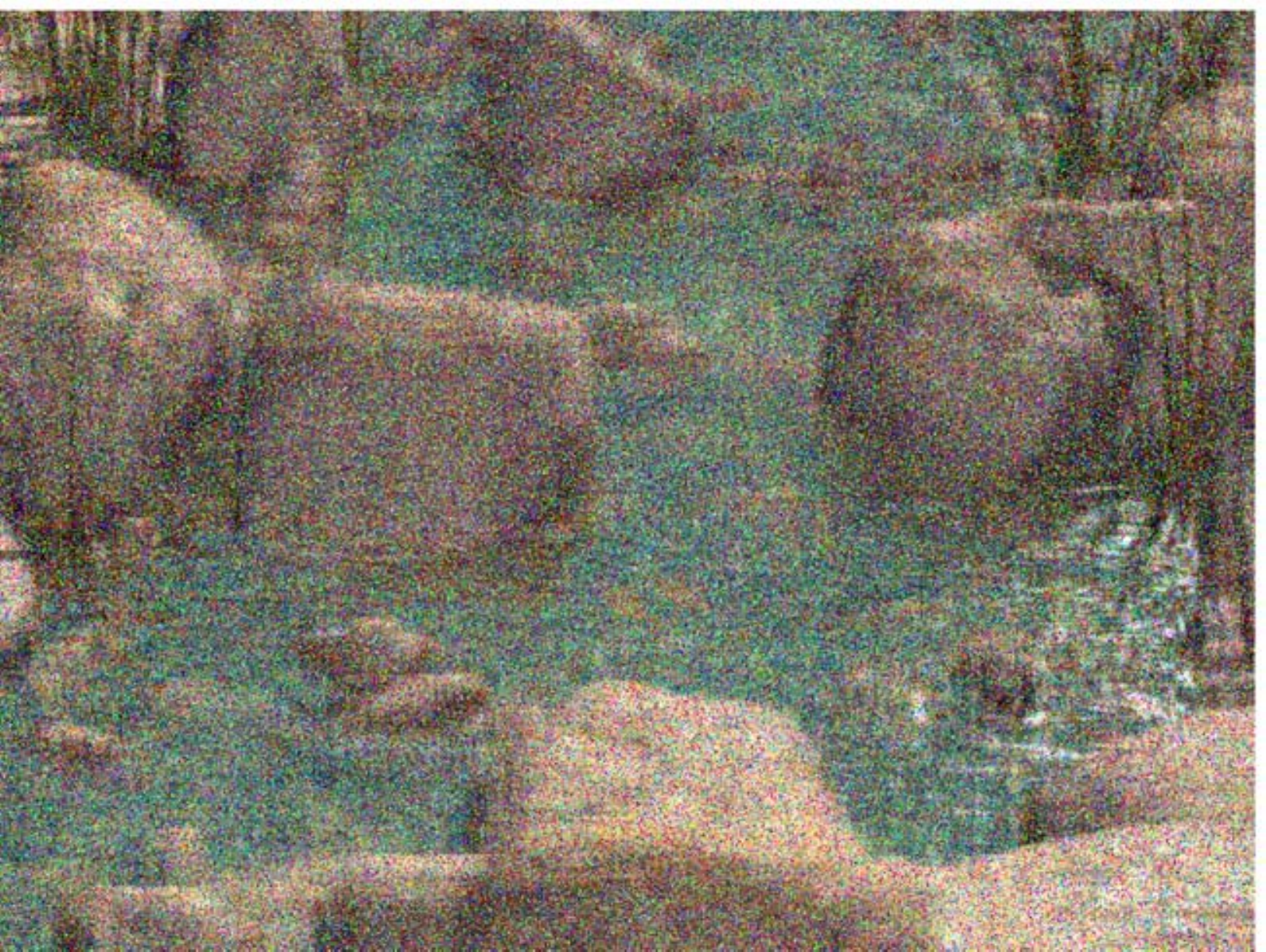}
&\includegraphics[width=0.11\textwidth,
keepaspectratio]{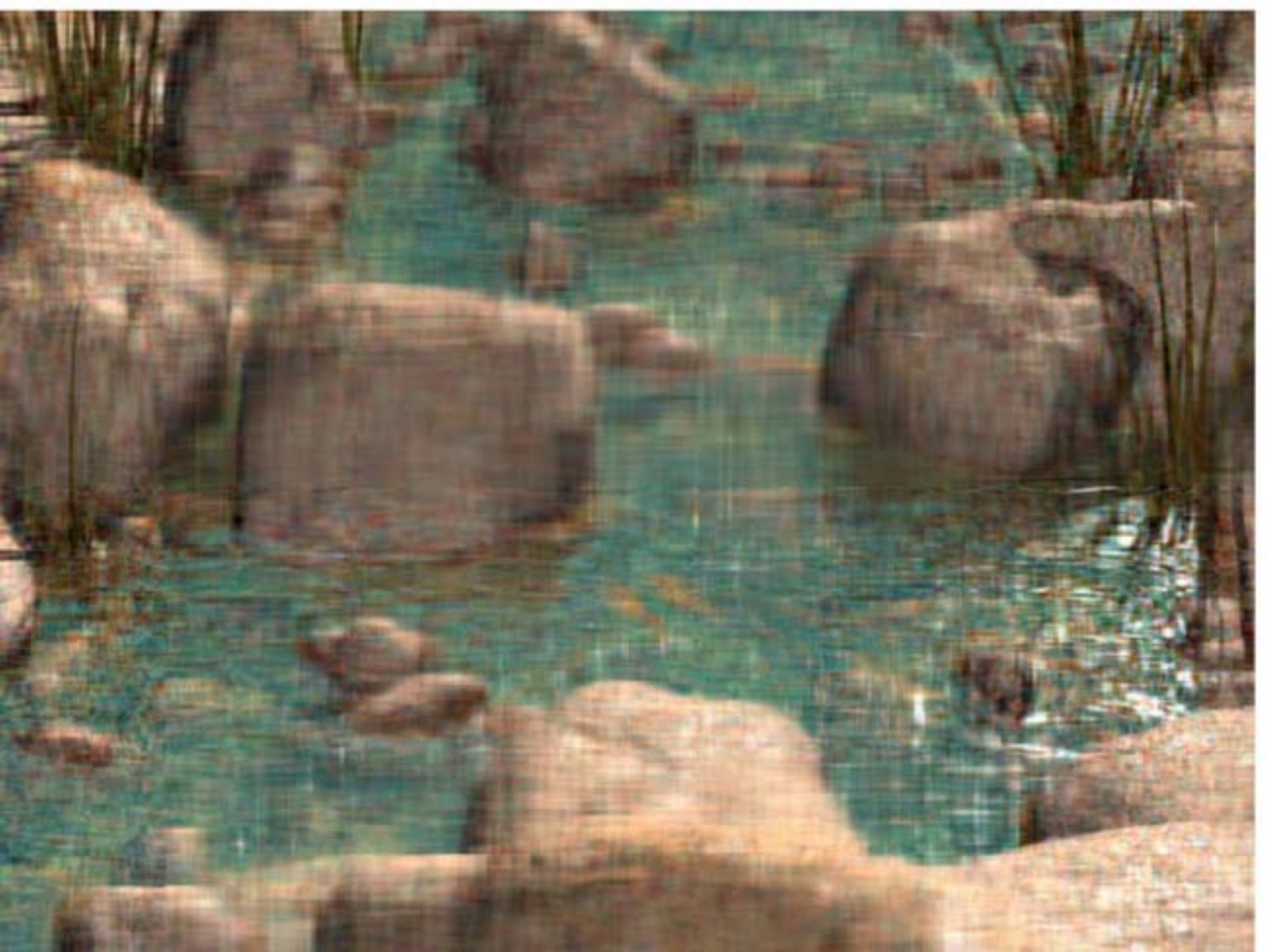}
&\includegraphics[width=0.11\textwidth,
keepaspectratio]{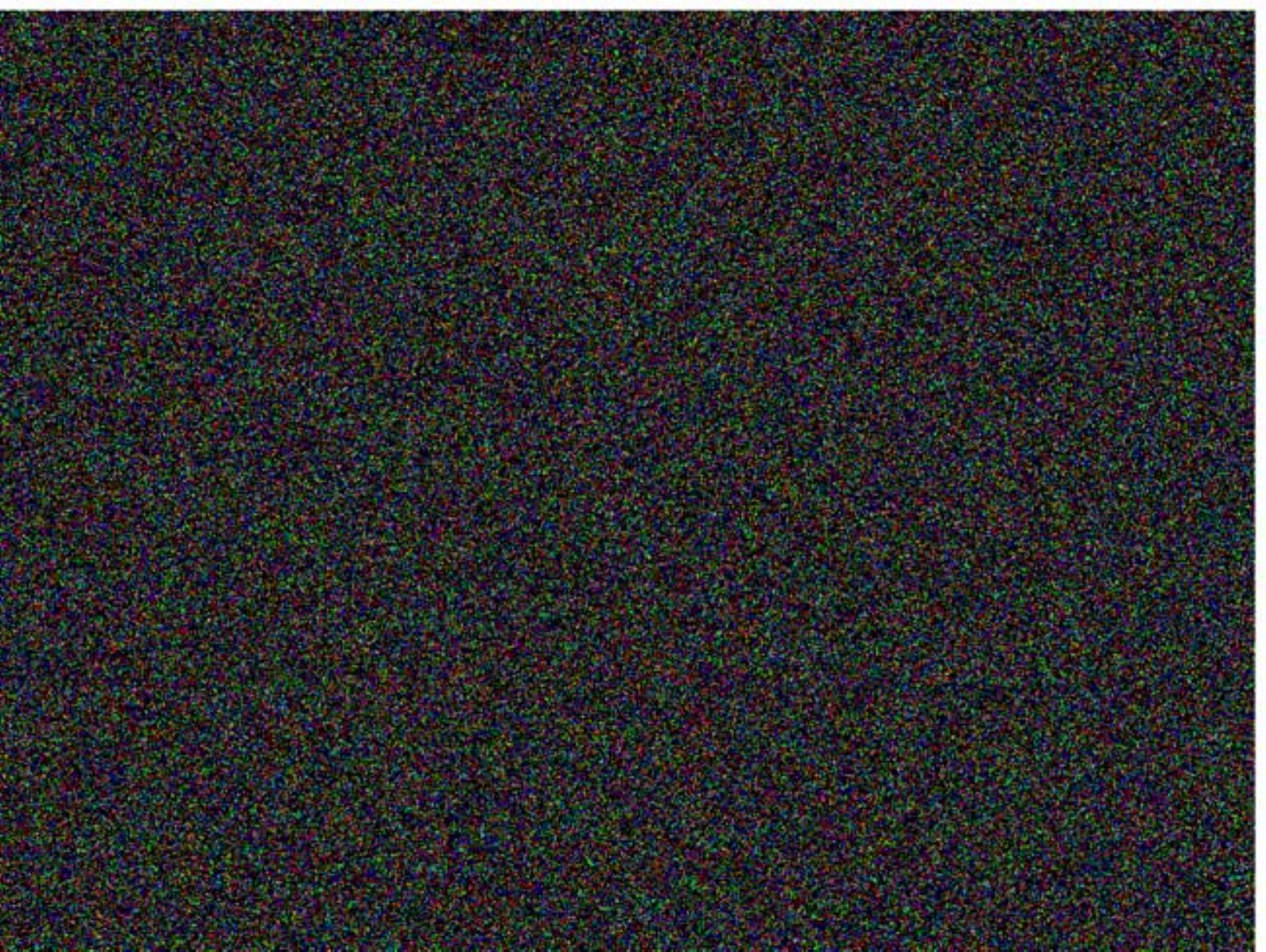}
&\includegraphics[width=0.11\textwidth,
keepaspectratio]{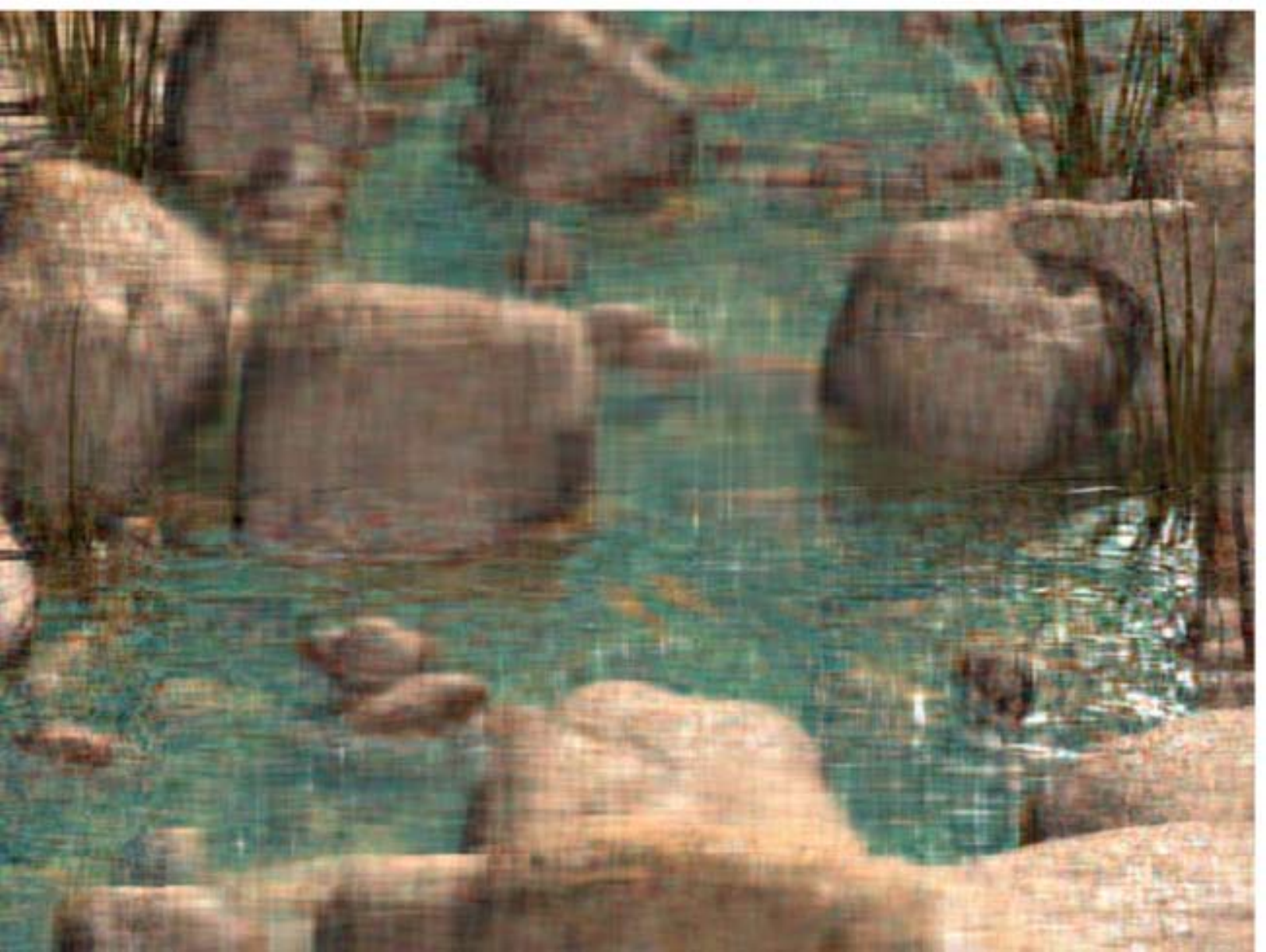}\\
(a) $\mathbf{M}$  & (b) $\mathbf{L}_0$  & (c) $\mathbf{S}_0$ & (d) $\mathbf{L}^*$\\
\includegraphics[width=0.11\textwidth,
keepaspectratio]{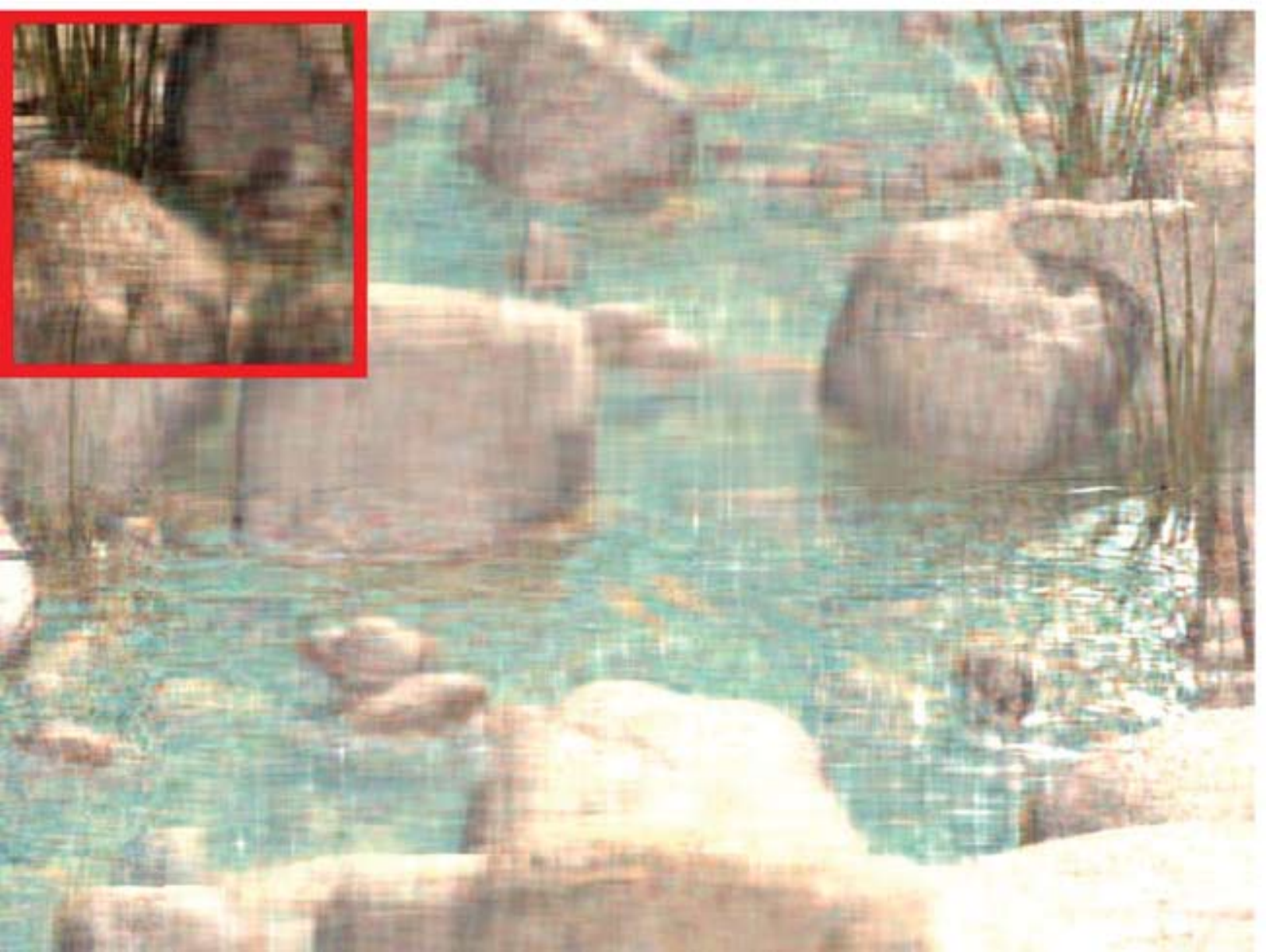}
&\includegraphics[width=0.11\textwidth,
keepaspectratio]{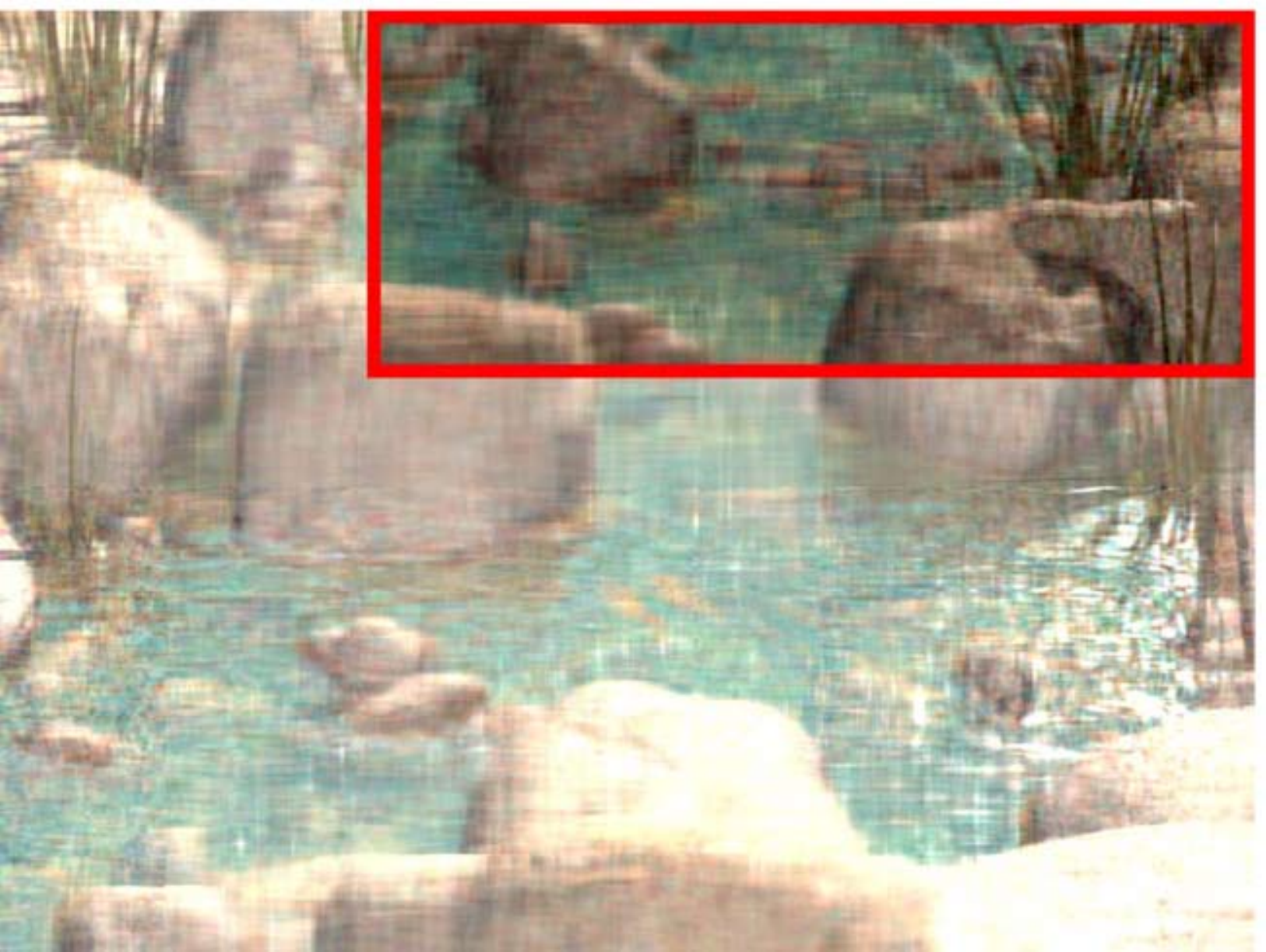}
&\includegraphics[width=0.11\textwidth,
keepaspectratio]{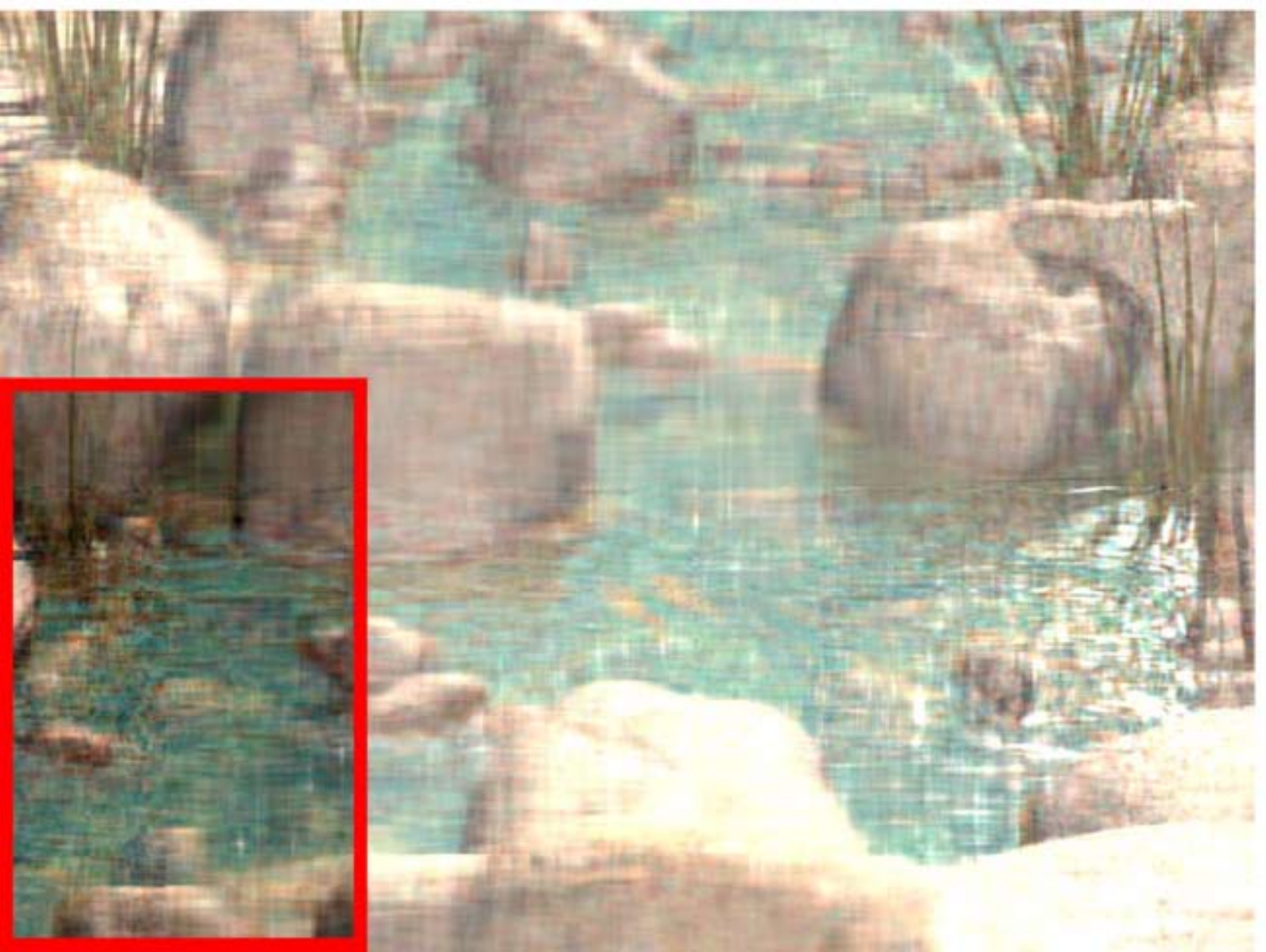}
&\includegraphics[width=0.11\textwidth,
keepaspectratio]{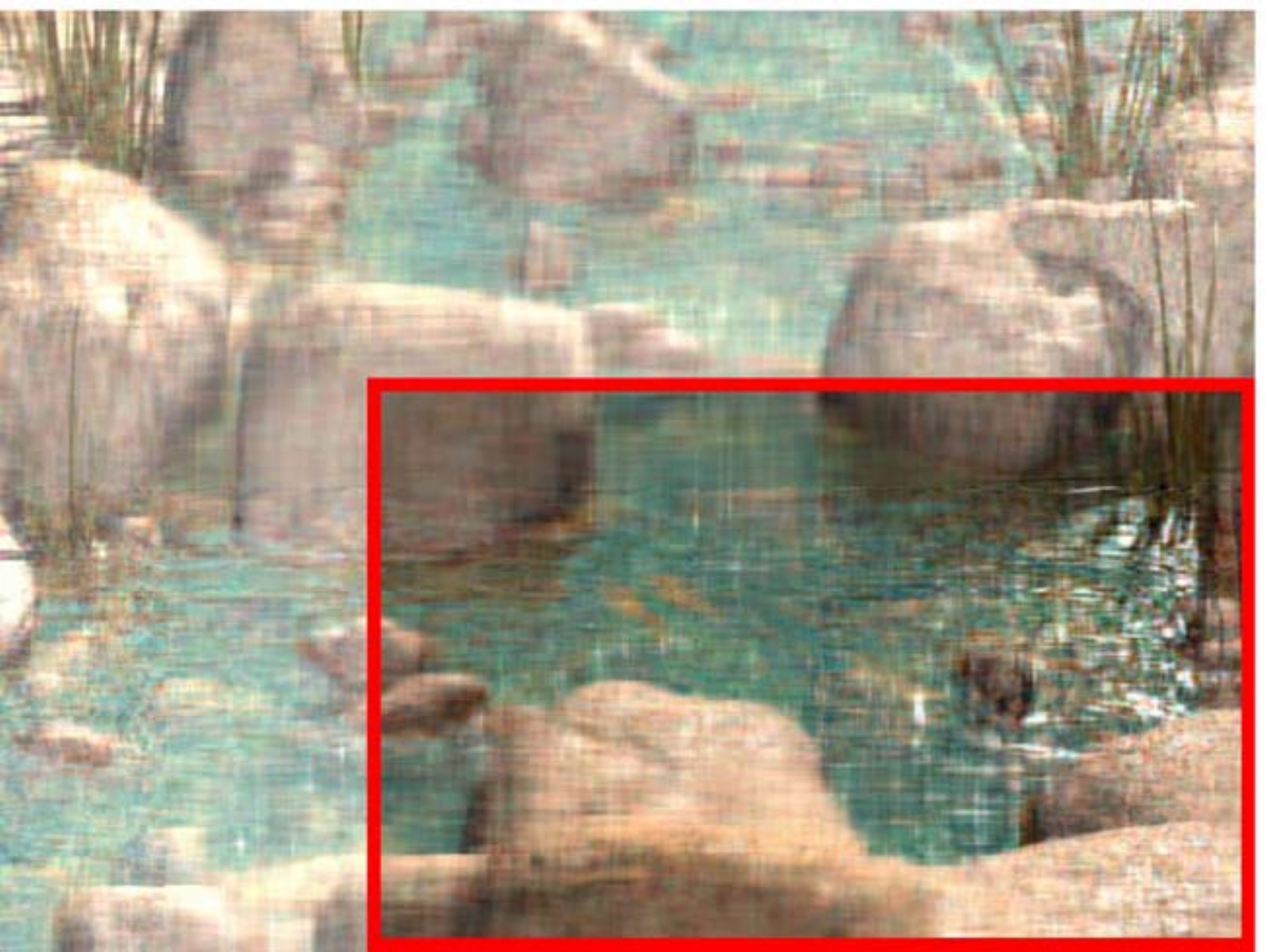}\\
\includegraphics[width=0.11\textwidth,
keepaspectratio]{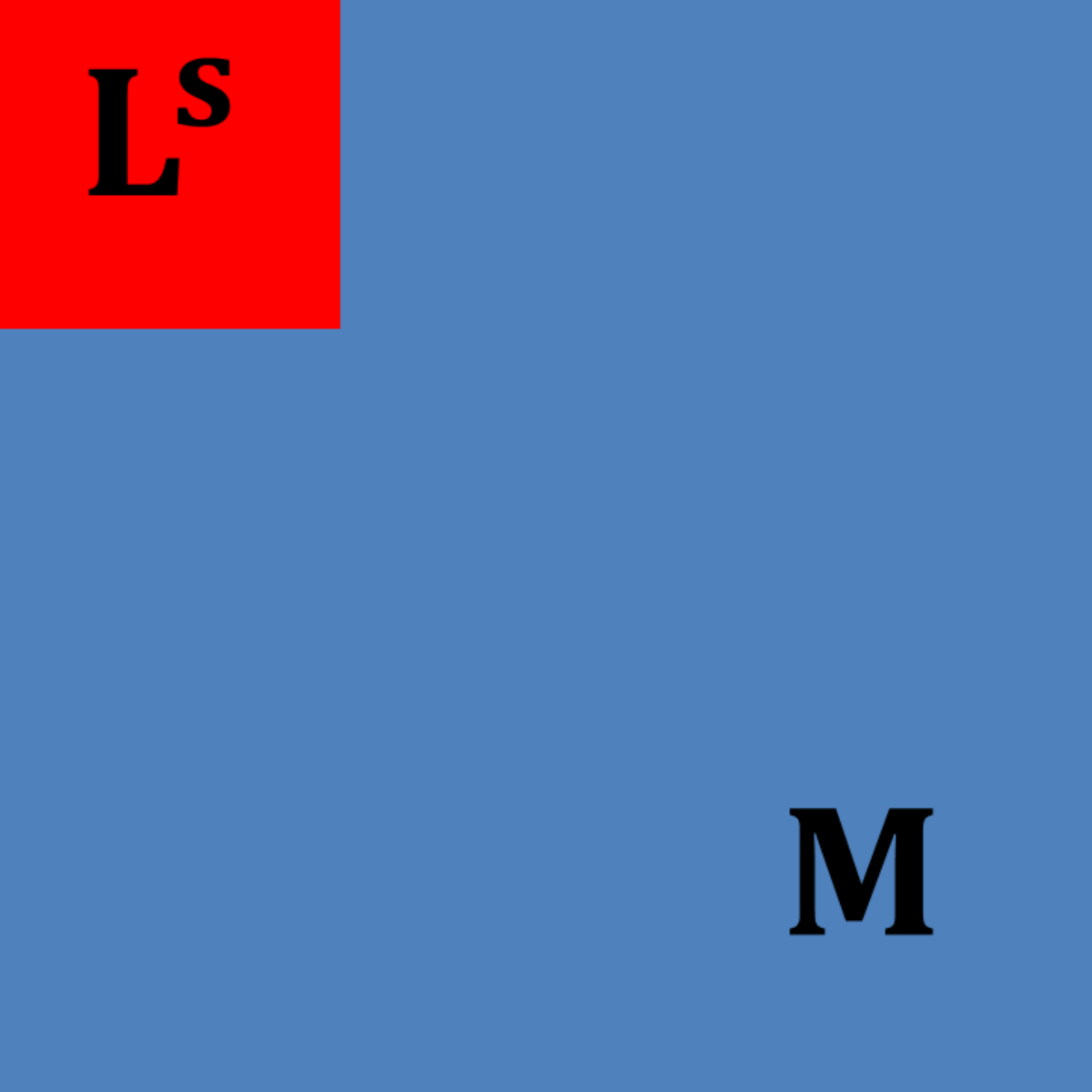}
&\includegraphics[width=0.11\textwidth,
keepaspectratio]{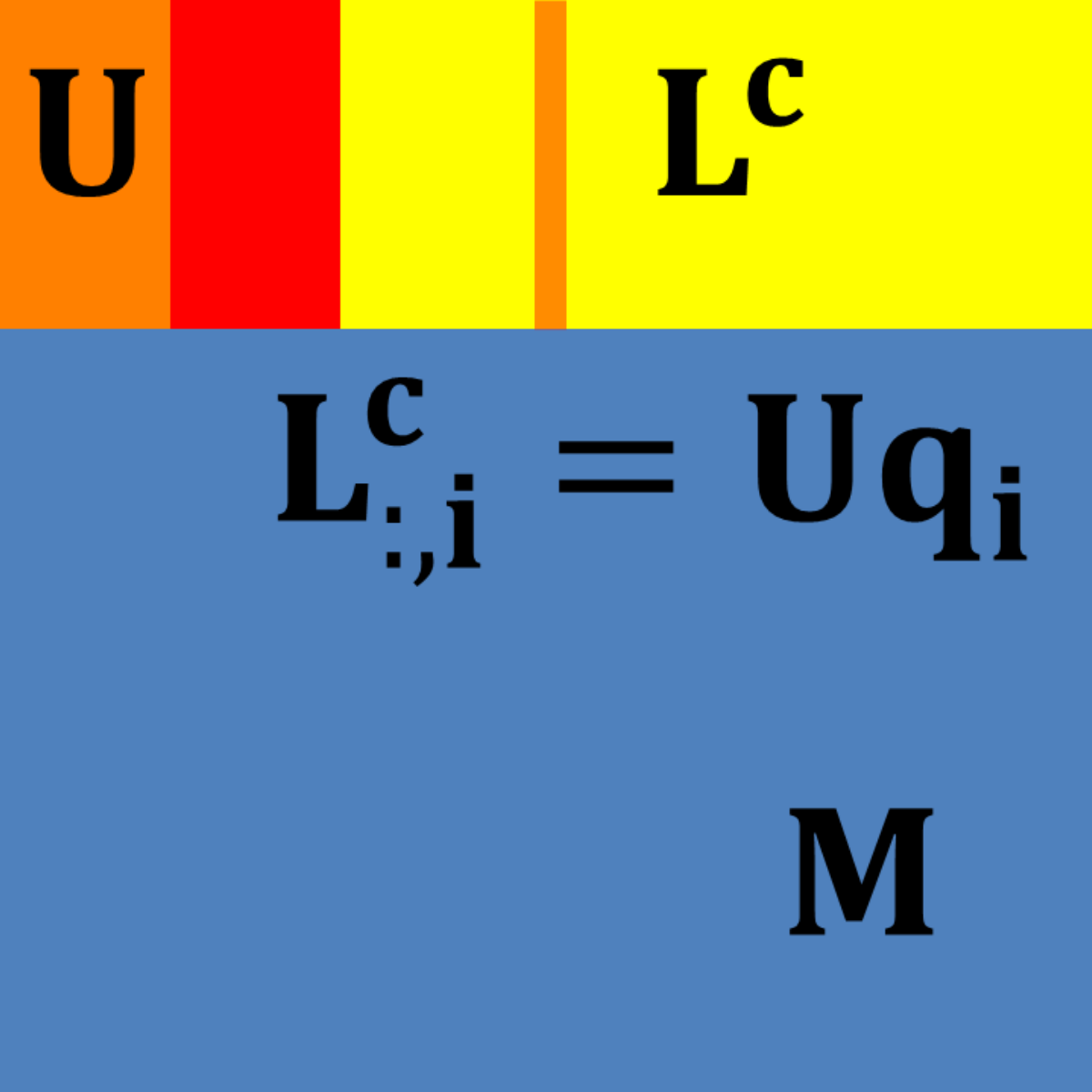}
&\includegraphics[width=0.11\textwidth,
keepaspectratio]{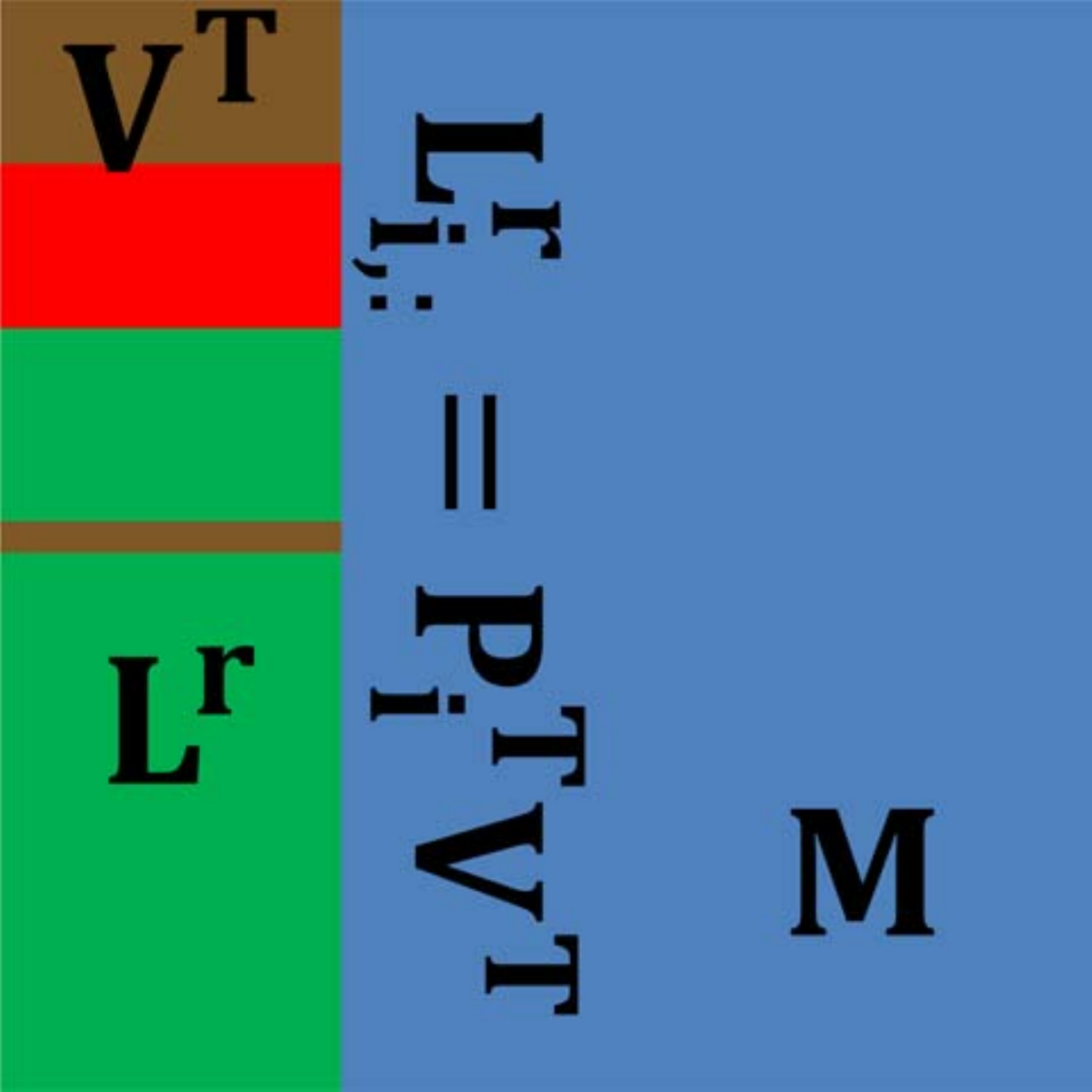}
&\includegraphics[width=0.11\textwidth,
keepaspectratio]{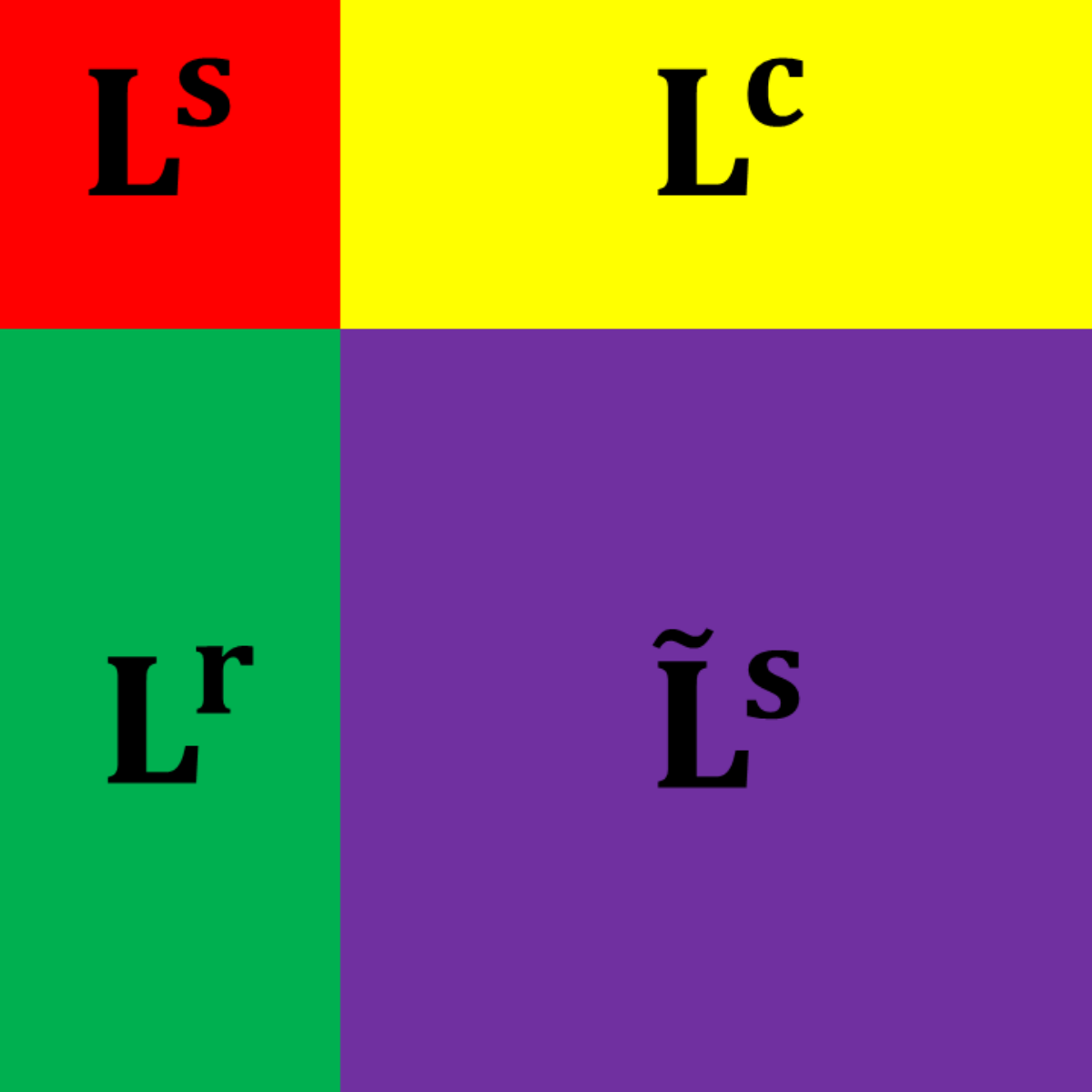}\\
(e) $\mathbf{L}^s$ & (f) $\mathbf{L}^c$ & (g) $\mathbf{L}^r$ & (h)
$\tilde{\mathbf{L}}^s$
\end{tabular}\\
\caption{Illustration of the proposed $l_1$ filtering method.
A large observed data matrix $\mathbf{M}$ (a) is the sum of a low-rank matrix
$\mathbf{L}_0$ (b) and a sparse matrix $\mathbf{S}_0$ (c).
The method first recovers a seed matrix (a submatrix of $\mathbf{L}_0$) $\mathbf{L}^s$ (e). Then the submatrices $\mathbf{L}^c$ (f) and $\mathbf{L}^r$ (g)
can be recovered by column and row filtering, respectively, where
$\mathbf{U}$ and $\mathbf{V}^T$ are the column space and row space
of $\mathbf{L}^s$, respectively. Then the complement matrix
$\tilde{\mathbf{L}}^s$ (h) can be represented by $\mathbf{L}^s$,
$\mathbf{L}^c$ and $\mathbf{L}^r$. Finally, we obtain the computed low-rank matrix $\mathbf{L}^*$ (d), which is identical to $\mathbf{L}_0$ with an overwhelming probability.}\label{fig:lena}
\end{figure}

Specifically, our method consists of two steps (illustrated in
Figure~\ref{fig:lena}). The first step is to recover a
submatrix\footnote{Note that the ``submatrix'' here does not
necessarily mean that we have to choose consecutive rows and
columns from $\mathbf{M}$.} $\mathbf{L}^s$ (Figure~\ref{fig:lena}
(e)) of $\mathbf{L}_0$. We call this submatrix the seed matrix
because all other entries of $\mathbf{L}_0$ can be further
calculated by this submatrix. The second step is to use the seed
matrix to recover two submatrices $\mathbf{L}^c$ and
$\mathbf{L}^r$ (Figures~\ref{fig:lena} (f)-(g)), which are on the
same rows and columns as $\mathbf{L}^s$ in $\mathbf{L}_0$,
respectively. They are recovered by minimizing the $l_1$ distance
from the subspaces spanned by the columns and rows of
$\mathbf{L}^s$, respectively. Hence we call this step $l_1$
filtering. The remaining part $\tilde{\mathbf{L}}_s$
(Figure~\ref{fig:lena} (h)) of $\mathbf{L}_0$ can be represented
by $\mathbf{L}^s$, $\mathbf{L}^c$ and $\mathbf{L}^r$, using the
generalized Nystr\"{o}m method (\cite{Wang2009Nystrom}). As
analyzed in Section~\ref{sec:comp}, our method is of linear cost
with respect to the data size. Besides the advantage of linear
time cost, the proposed algorithm is also highly parallel: the
columns of $\mathbf{L}^c$ and the rows of $\mathbf{L}^r$ can be
recovered fully independently. We also prove that under suitable
conditions, our method can \emph{exactly} recover the underling
low-rank matrix $\mathbf{L}_0$ with an overwhelming probability.
To our best knowledge, this is the first algorithm that can
\emph{exactly} solve a nuclear norm minimization problem in
\emph{linear time}.

\section{Previous Works}\label{sec:previous_work}
In this section, we review some previous algorithms for solving
PCP. The existing solvers can be roughly divided into three
categories: classic convex optimization, factorization and
compressed optimization.

For small sized problems, PCP can be reformulated as a
semidefinite program and then be solved by standard interior point
methods. However, this type of methods cannot handle even moderate
scale matrices due to their $O(n^6)$ complexity in each iteration.
So people turned to first-order algorithms, such as the dual
method (\cite{Ganesh-2009}), the accelerated proximal gradient
(APG) method (\cite{Ganesh-2009}) and the alternating direction
method (ADM) (\cite{lin2009augmented}), among which ADM is the
most efficient. All these methods require solving the following
kind of subproblem in each iteration
\begin{equation}
\min\limits_\mathbf{A} \eta\|\mathbf{A}\|_*
+\frac{1}{2}\|\mathbf{A}-\mathbf{W}\|_F^2,\label{eq:Proxy}
\end{equation}
where $\|\cdot\|_F$ is the Frobenious norm.
\cite{cai2008singular} proved that the above problem has a closed
form solution
\begin{equation}
\mathbf{A}=\mathbf{U}\mathcal{S}_{\eta}(\mathbf{\Sigma})\mathbf{V}^T,\label{eq:SVT}
\end{equation}
where $\mathbf{U}\mathbf{\Sigma}\mathbf{V}^T$ is the singular
value decomposition of $\mathbf{W}$ and
$\mathcal{S}_{\eta}(x)=\mbox{sgn}(x)\max(|x|-\eta,0)$ is the soft
shrinkage operator. Therefore, these methods all require computing
SVDs for some matrices, resulting in $O(mn\min(m,n))$ complexity,
where $m\times n$ is the matrix size.

As the most expensive computational task required by solving
(\ref{eq:PCP}) is to perform SVD, \cite{lin2009augmented} adopted
partial SVD (\cite{Larsen-1998-PROPACK}) to reduce the complexity
at each iteration to $O(rmn)$, where $r$ is the target rank.
However, such a complexity is still too high for very large data
sets. \cite{Drineas2006LTSVD} developed a fast Monte Carlo
algorithm, named linear time SVD (LTSVD), which can be used for
solving SVDs approximately (also see \cite{Halko-2011-random}).
The main drawback of LTSVD is that it is less accurate than the
standard SVD as it uses random sampling. So the whole algorithm
needs more iterations to achieve the same accuracy. As a
consequence, the speed performance of LTSVD quickly deteriorates
when the target rank increases (see Figure~\ref{fig:rank_ratio}).
Actually, even adopting LTSVD the whole algorithm is still
quadratic w.r.t. the data size because it still requires
matrix-matrix multiplication in each iteration.

To address the scalability issue of solving large-scale PCP
problems, \cite{shen-MatrixFac} proposed a factorization based
method, named low-rank matrix fitting (LMaFit). This approach
represents the low-rank matrix as a product of two matrices and
then minimizes over the two matrices alternately. Although they do
not require nuclear norm minimization (hence the SVDs), the
convergence of the proposed algorithm is not guaranteed as the
corresponding problem is non-convex. Moreover, both the
matrix-matrix multiplication and the QR decomposition based rank
estimation technique require $O(rmn)$ complexity. So this method
does not essentially reduce the complexity.

Inspired by compressed optimization, \cite{mu-2011-RandomProj}
proposed reducing the problem scale by random projection (RP).
However, this method is highly unstable -- different random
projections may lead to radically different results. Moreover, the
need to introduce additional constraint to the problem slows down
the convergence. And actually, the complexity of this method is
also $O(pmn)$, where $p\times m$ is the size of the random
projection matrix and $p>r$. So this method is still not of linear
complexity with respect to the matrix size.

\section{The $l_1$ Filtering Algorithm}\label{sec:l1f}
Given an observed data matrix $\mathbf{M} \in \mathbb{R}^{m\times n}$, which is the sum of a low-rank matrix $\mathbf{L}_0$ and a sparse matrix $\mathbf{S}_0$, PCP is to recover $\mathbf{L}_0$ from $\mathbf{M}$.
What our approach differs from traditional ones is that the underlying low-rank matrix $\mathbf{L}_0$ is reconstructed from a seed matrix.
As explained in Section~\ref{sec:outline}, our $l_1$ filtering
algorithm consists of two steps: first recovering a seed matrix, second performing $l_1$ filtering on
corresponding rows and columns of the data matrix. Below we
provide details of these two steps.

\subsection{Seed Matrix Recovery}\label{sec:submr}
Suppose that the target rank $r$ is very small compared with the
data size: $r\ll \min(m,n)$. We first randomly sample an $(s_r
r)\times (s_c r)$ submatrix $\mathbf{M}^s$ from $\mathbf{M}$,
where $s_r>1$ and $s_c >1$ are the row and column oversampling
rates, respectively. Then the submatrix $\mathbf{L}^s$ of the
underlying matrix $\mathbf{L}_0$ can be recovered by solving a
small sized PCP problem:
\begin{equation}
\min\limits_{\mathbf{L}^s,\mathbf{S}^s} \|\mathbf{L}^s\|_* +
\tilde{\lambda}\|\mathbf{S}^s\|_{l_1}, \quad s.t. \quad \mathbf{M}^s =
\mathbf{L}^s + \mathbf{S}^s,\label{eq:SPCP}
\end{equation}
e.g., using ADM (\cite{lin2009augmented}), where
$\tilde{\lambda}=1/\sqrt{\max(s_r r,s_c r)}$.

By Theorem~1.1 in (\cite{candes2009robust}), the seed matrix
$\mathbf{L}^s$ can be exactly recovered from $\mathbf{M}^s$ with
an overwhelming probability when $s_r$ and $s_c$ increases. In
fact, by that theorem $s_r$ and $s_c$ should be chosen at the
scale of $O(\ln^2 r)$. For the experiments conducted in this
paper, whose $r$'s are very small, we simply choose $s_c=s_r=10$.

\subsection{$l_1$ Filtering}
For ease of illustration, we assume that $\mathbf{M}^s$ is the
top left $(s_r r)\times (s_c r)$ submatrix of $\mathbf{M}$. Then
accordingly $\mathbf{M}$, $\mathbf{L}_0$ and $\mathbf{S}_0$ can be
partitioned into:
\begin{equation}
\mathbf{M} = \begin{bmatrix}\mathbf{M}^s &
\mathbf{M}^c\\
\mathbf{M}^r & \tilde{\mathbf{M}}^s\end{bmatrix},\ \mathbf{L}_0=
\begin{bmatrix}
\mathbf{L}^s &
\mathbf{L}^c\\
\mathbf{L}^r & \tilde{\mathbf{L}}^s
\end{bmatrix},\
\mathbf{S}_0 =
\begin{bmatrix}\mathbf{S}^s &
\mathbf{S}^c\\
\mathbf{S}^r &
\tilde{\mathbf{S}}^s\end{bmatrix}.\label{eq:partition}
\end{equation}
Since $\rank(\mathbf{L}_0)=\rank(\mathbf{L}^s)=r$, there must
exist matrices $\mathbf{Q}$ and $\mathbf{P}$, such that
\begin{equation}
\mathbf{L}^c = \mathbf{L}^s\mathbf{Q} \quad\mbox{and}\quad \mathbf{L}^r =
\mathbf{P}^T\mathbf{L}^s.\label{eq:row-column-recovery}
\end{equation}
As $\mathbf{S}_0$ is sparse, so are $\mathbf{S}^c$ and
$\mathbf{S}^r$. Therefore, $\mathbf{Q}$ and $\mathbf{P}$ can be
found by solving the following problems:
\begin{equation}
\min\limits_{\mathbf{S}^c,\mathbf{Q}}\|\mathbf{S}^c\|_{l_1}, \quad
s.t. \quad \mathbf{M}^c = \mathbf{L}^s\mathbf{Q} +
\mathbf{S}^c,\label{eq:l1c}
\end{equation}
and
\begin{equation}
\min\limits_{\mathbf{S}^r,\mathbf{P}}\|\mathbf{S}^r\|_{l_1}, \quad
s.t. \quad \mathbf{M}^r = \mathbf{P}^T\mathbf{L}^s +
\mathbf{S}^r,\label{eq:l1r}
\end{equation}
respectively. The above two problems can be easily solved by ADM.

With $\mathbf{Q}$ and $\mathbf{P}$ computed, $\mathbf{L}^c$ and
$\mathbf{L}^r$ are obtained as (\ref{eq:row-column-recovery}).
Again by $\rank(\mathbf{L}_0)=\rank(\mathbf{L}^s)=r$, the
generalized Nystr\"{o}m method (\cite{Wang2009Nystrom}) gives:
\begin{equation}
\tilde{\mathbf{L}}^s=\mathbf{L}^r (\mathbf{L}^s)^\dag
\mathbf{L}^c,\label{eq:l1s}
\end{equation}
where $(\mathbf{L}^s)^\dag$ is the Moore-Penrose pseudo inverse of
$\mathbf{L}^s$.

In real computation, as the SVD of $\mathbf{L}^s$ is readily
available when solving (\ref{eq:SPCP}), due to the singular value
thresholding operation (\ref{eq:SVT}), it is more convenient to
reformulate (\ref{eq:l1c}) and (\ref{eq:l1r}) as
\begin{equation}
\min\limits_{\mathbf{S}^c,\tilde{\mathbf{Q}}}\|\mathbf{S}^c\|_{l_1},
\quad s.t. \quad \mathbf{M}^c = \mathbf{U}^s\tilde{\mathbf{Q}} +
\mathbf{S}^c,\label{eq:l1c'}
\end{equation}
and
\begin{equation}
\min\limits_{\mathbf{S}^r,\tilde{\mathbf{P}}}\|\mathbf{S}^r\|_{l_1},
\quad s.t. \quad \mathbf{M}^r =
\tilde{\mathbf{P}}^T(\mathbf{V}^s)^T +
\mathbf{S}^r,\label{eq:l1r'}
\end{equation}
respectively, where $\mathbf{U}^s\mathbf{\Sigma}^s(\mathbf{V}^s)^T$ is the
skinny SVD of $\mathbf{L}^s$ obtained from (\ref{eq:SVT}) in the
iterations. Such a reformulation has multiple advantages. First,
as
$(\mathbf{U}^s)^T\mathbf{U}^s=(\mathbf{V}^s)^T\mathbf{V}^s=\mathbf{I}$,
it is unnecessary to compute the inverse of
$(\mathbf{U}^s)^T\mathbf{U}^s$ and $(\mathbf{V}^s)^T\mathbf{V}^s$
when updating $\tilde{\mathbf{Q}}$ and $\tilde{\mathbf{P}}$ in the
iterations of ADM. Second, computing (\ref{eq:l1s}) also becomes
easy if one wants to form $\tilde{\mathbf{L}}^s$ explicitly
because now
\begin{equation}
\tilde{\mathbf{L}}^s=\tilde{\mathbf{P}}^T(\mathbf{\Sigma}^s)^{-1}\tilde{\mathbf{Q}}.\label{eq:l1s'}
\end{equation}

To make the algorithm description complete, we sketch in
Algorithm~\ref{alg:adm} the ADM for solving (\ref{eq:l1c'}) and
(\ref{eq:l1r'}), which are both of the following form:
\begin{equation}
\min\limits_{\mathbf{E},\mathbf{Z}}\|\mathbf{E}\|_{l_1}, \quad
s.t. \quad \mathbf{X} = \mathbf{A}\mathbf{Z} +
\mathbf{E},\label{eq:l1g}
\end{equation}
where $\mathbf{X}$ and $\mathbf{A}$ are known matrices and
$\mathbf{A}$ has orthonormal columns, i.e.,
$\mathbf{A}^T\mathbf{A}=\mathbf{I}$. The ADM for (\ref{eq:l1g}) is
to minimize on the following augmented Lagrangian function
\begin{equation}
\|\mathbf{E}\|_{l_1}+\langle\mathbf{Y},\mathbf{X} -
\mathbf{A}\mathbf{Z} - \mathbf{E}\rangle +
\frac{\beta}{2}\|\mathbf{X} - \mathbf{A}\mathbf{Z} -
\mathbf{E}\|_F^2,
\end{equation}
with respect to $\mathbf{E}$ and $\mathbf{Z}$, respectively, by
fixing other variables, and then update the Lagrange multiplier
$\mathbf{Y}$ and the penalty parameter $\beta$.\footnote{The ADM
for solving PCP follows the same methodology. As a reader can
refer to (\cite{lin2009augmented}, \cite{yuan2009sparse}) for
details, we omit the pseudo code for using ADM to solve PCP.}
\begin{algorithm}
   \caption{Solving (\ref{eq:l1g}) by ADM}\label{alg:adm}
\begin{algorithmic}
   \STATE {\bfseries Input:} $\mathbf{X}$ and $\mathbf{A}$.
   \STATE {\bfseries Initialize:} Set $\mathbf{E}_0$, $\mathbf{Z}_0$ and $\mathbf{Y}_0$ to zero matrices. Set $\varepsilon >0$, $\rho > 1$ and $\bar{\beta} \gg \beta_0 >0$.
    \WHILE {$\|\mathbf{X} - \mathbf{A}\mathbf{Z}_{k} - \mathbf{E}_k\|_{l_\infty}/\|\mathbf{X}\|_{l_\infty}\geq \varepsilon$}
   \STATE {\bfseries Step 1:} Update $\mathbf{E}_{k+1}=\mathcal{S}_{\beta_k^{-1}}(\mathbf{X}-\mathbf{A}\mathbf{Z}_k+\mathbf{Y}_k/\beta_{k}$), where $\mathcal{S}$ is the soft-thresholding operator~(\cite{cai2008singular}).
   \STATE {\bfseries Step 2:} Update $\mathbf{Z}_{k+1}=\mathbf{A}^T(\mathbf{X}-\mathbf{E}_{k+1}+\mathbf{Y}_k/\beta_{k})$.
   \STATE {\bfseries Step 3:} Update $\mathbf{Y}_{k+1}=\mathbf{Y}_k + \beta_k(\mathbf{X} - \mathbf{A}\mathbf{Z}_{k+1} - \mathbf{E}_{k+1})$ and $\beta_{k+1}=\min(\rho\beta_k, \bar{\beta})$.
   \ENDWHILE
\end{algorithmic}
\end{algorithm}

Note that it is easy to see that (\ref{eq:l1c'}) and (\ref{eq:l1r'})
can also be solved in full parallelism as the columns and rows of
$\mathbf{L}^c$ and $\mathbf{L}^r$ can computed independently,
thanks to the decomposability of the problems. So the recovery of
$\mathbf{L}^c$ and $\mathbf{L}^r$ is very efficient if one has a
parallel computing platform, such as a general purpose graphics
processing unit (GPU).

\subsection{The Complete Algorithm}
Now we are able to summarize in Algorithm~\ref{alg:l1fmrs} our
$l_1$ filtering method for solving PCP, where steps 3 and 4 can be
done in parallel.
\begin{algorithm}
   \caption{The $l_1$ Filtering Method for Solving PCP (\ref{eq:PCP})}
\begin{algorithmic}\label{alg:l1fmrs}
   \STATE {\bfseries Input:} Observed data matrix $\mathbf{M}$.
   \STATE {\bfseries Step 1:} Randomly sample a submatrix
   $\mathbf{M}^s$.
   \STATE {\bfseries Step 2:} Solve the small sized PCP problem
   (\ref{eq:SPCP}), e.g., by ADM, to recover the seed matrix $\mathbf{L}^s$.
   \STATE {\bfseries Step 3:} Reconstruct $\mathbf{L}^c$ by solving (\ref{eq:l1c'}).
   \STATE {\bfseries Step 4:} Reconstruct $\mathbf{L}^r$ by solving (\ref{eq:l1r'}).
   \STATE {\bfseries Step 5:} Represent $\tilde{\mathbf{L}}^s$ by (\ref{eq:l1s'}).
   \STATE {\bfseries Output:} Low-rank matrix $\mathbf{L}$ and sparse matrix $\mathbf{S}=\mathbf{M}-\mathbf{L}$.
\end{algorithmic}
\end{algorithm}

\subsection{Complexity Analysis}\label{sec:comp}
Now we analyze the computational complexity of the proposed
Algorithm~\ref{alg:l1fmrs}. For the step of seed matrix recovery,
the complexity of solving (\ref{eq:SPCP}) is only $O(r^3)$. For
the $l_1$ filtering step, it can be seen that the complexity of
solving (\ref{eq:l1c'}) and (\ref{eq:l1r'}) is $O(r^2n)$ and
$O(r^2m)$, respectively. So the total complexity of this step is
$O(r^2(m+n))$. As the remaining part $\tilde{\mathbf{L}}_s$ of
$\mathbf{L}_0$ can be represented by $\mathbf{L}^s$,
$\mathbf{L}^c$ and $\mathbf{L}^r$, using the generalized
Nystr\"{o}m method (\cite{Wang2009Nystrom})\footnote{Of course, if
we explicitly form $\tilde{\mathbf{L}}^s$ then this step costs no
more than $rmn$ complexity. Compared with other methods, our rest
computations are all of $O(r^2(m+n))$ complexity at the most,
while those methods all require at least $O(rmn)$ complexity
\emph{in each iteration}, which results from matrix-matrix
multiplication.} and recall that $r \ll \min(m,n)$, we conclude
that the overall complexity of Algorithm~\ref{alg:l1fmrs} is
$O(r^2(m+n))$, which is only of linear cost with respect to the
data size.

\subsection{Exact Recoverability of $l_1$ Filtering}\label{sec:Recoverability}
The exact recoverability of $\mathbf{L}_0$ using our $l_1$
filtering method consists of two factors. First, exactly
recovering $\mathbf{L}^s$ from $\mathbf{M}^s$. Second, exactly
recovering $\mathbf{L}^c$ and $\mathbf{L}^r$. If all
$\mathbf{L}^s$, $\mathbf{L}^c$, and $\mathbf{L}^r$ can be exactly
recovered, $\mathbf{L}_0$ is exactly recovered.

The exact recoverability of $\mathbf{L}^s$ from $\mathbf{M}^s$ is
guaranteed by Theorem~1.1 of (\cite{candes2009robust}). When $s_r$
and $s_c$ are sufficiently large, the chance of success is
overwhelming.

To analyze the exact recoverability of $\mathbf{L}^c$ and
$\mathbf{L}^r$, we first observe that it is equivalent to the
exact recoverability of $\mathbf{S}^c$ and $\mathbf{S}^r$. By
multiplying annihilation matrices $\mathbf{U}^{s,\perp}$ and
$\mathbf{V}^{s,\perp}$ to both sides of (\ref{eq:l1c'}) and
(\ref{eq:l1r'}), respectively, we may recover $\mathbf{S}^c$ and
$\mathbf{S}^r$ by solving
\begin{equation}
\min\limits_{\mathbf{S}^c}\|\mathbf{S}^c\|_{l_1}, \quad s.t. \quad
\mathbf{U}^{s,\perp}\mathbf{M}^c =
\mathbf{U}^{s,\perp}\mathbf{S}^c,\label{eq:l1c''}
\end{equation}
and
\begin{equation}
\min\limits_{\mathbf{S}^r}\|\mathbf{S}^r\|_{l_1}, \quad s.t. \quad
\mathbf{M}^r(\mathbf{V}^{s,\perp})^T =
\mathbf{S}^r(\mathbf{V}^{s,\perp})^T,\label{eq:l1r''}
\end{equation}
respectively. If the oversampling rates $s_c$ and $s_r$ are large
enough, we are able to choose $\mathbf{U}^{s,\perp}$ and
$\mathbf{V}^{s,\perp}$ that are close to Gaussian random matrices.
Then we may apply the standard theory in compressed sensing
(\cite{candes2007intro}) to conclude that if the oversampling rates
$s_c$ and $s_r$ are large enough and $\mathbf{S}^c$ and
$\mathbf{S}^r$ are sparse enough\footnote{As the analysis in the compressed sensing theories is qualitative and the bounds are actually pessimistic, copying those inequalities here is not very useful. So we omit the mathematical descriptions for brevity.}, $\mathbf{S}^c$ and
$\mathbf{S}^r$ can be exactly recovered with an overwhelming
probability.

We also present an example in Figure~\ref{fig:lena} to illustrate
the exact recoverability of $l_1$ filtering. We first truncate the
SVD of an $1024 \times 768$ image ``Water'' \footnote{The image is
available at
http://www.petitcolas.net/fabien/\\watermarking/image\_database/.}
to get a matrix of rank 30 (Figure~\ref{fig:lena} (b)). The
observed image (Figure~\ref{fig:lena} (a)) is obtained from
Figure~\ref{fig:lena} (b) by adding large noise to 30$\%$ of the
pixels uniformly sampled at random (Figure~\ref{fig:lena} (c)).
Suppose we have the top-left $300 \times 300$ submatrix as the
seed (Figure~\ref{fig:lena} (e)), the low-rank image
(Figure~\ref{fig:lena} (d)) can be exactly recovered by $l_1$
filtering. Actually, the relative reconstruction errors in
$\mathbf{L}^*$ is only $7.03\times 10^{-9}$.

\comment{As both the seed matrix recovery and $l_1$ filtering can
be done successfully with an overwhelming probability, if the
oversampling rates $s_c$ and $s_r$ are large enough and
$\mathbf{S}^c$ and $\mathbf{S}^r$ are sparse enough, we can
conclude that $\mathbf{L}_0$ and $\mathbf{S}_0$ can be
\emph{exactly} recovered with an overwhelming probability. }

\subsection{Target Rank Estimation}
The above analysis and computation are all based on a known value
of the target rank $r$. For some applications, we could have an
estimate on $r$. For example, for the background modeling problem
(\cite{Wright-NIPS2009}), the rank of the background video should be
very close to one as the background hardly changes; and for the
photometric stereo problem (\cite{wu-robust}) the rank of the
surface normal map should be very close to three as the normals
are three dimensional vectors. However, the rank $r$ of the
underlying matrix might not always be known. So we have to provide
a strategy to estimate $r$.

As we assume that the size $m'\times n'$ of submatrix
$\mathbf{M}^s$ is $(s_r r)\times (s_c r)$, where $s_r$ and $s_c$
should be sufficiently large in order to ensure the exact recovery
of $\mathbf{L}^s$ from $\mathbf{M}^s$, after we have computed
$\mathbf{L}^s$ by solving (\ref{eq:SPCP}), we may check whether
\begin{equation}
m'/r'\geq s_r \quad\mbox{and}\quad n'/r'\geq s_c \label{eq:rank_test}
\end{equation}
are satisfied, where $r'$ is the rank of $\mathbf{L}^s$. If yes,
$\mathbf{L}^s$ is accepted as a seed matrix. Otherwise, it implies
that $m'\times n'$ may be too small with respect to the target
rank $r$. Then we may increase the size of the submatrix to $(s_r
r')\times (s_c r')$ and repeat the above procedure until
(\ref{eq:rank_test}) is satisfied or
\begin{equation}
\max(m'/m, n'/n) > 0.5. \label{eq:rank_test1}
\end{equation}
We require (\ref{eq:rank_test1}) because the speed advantage of
our $l_1$ filtering algorithm will quickly lost beyond this size
limit (see Figure~\ref{fig:rank_ratio}). If we have to use a
submstrix whose size should be greater than $(0.5m)\times(0.5n)$,
then the target rank should be comparable to the size of data,
hence breaking our low-rank assumption. In this case, we may
resort to the usual method to solve PCP.

Of course, we may sample one more submatrix to cross validate the
estimated target rank $r$. When $r$ is indeed very small, such a
cross validation is not a big overhead.


\section{Experimental Results}\label{sec:exp}
In this section, we present experiments on both synthetic data and
real vision problems (structure from motion and background
modeling) to test the performance of $l_1$ filtering. All the
experiments are conducted and timed on the same PC with an AMD
Athlon\textsuperscript{\textregistered} II X4 2.80GHz CPU that has
4 cores and 6GB memory, running Windows 7 and Matlab (Version
7.10).

\subsection{Comparison Results for Solving PCP}

We first test the performance of $l_1$ filtering on solving PCP problem (\ref{eq:PCP}).
The experiments are categorized into the following three classes:
\begin{enumerate}
\item Compare with classic numerical solvers (e.g., ADM
\cite{lin2009augmented} and its variation, denoted as LTSVD ADM,
which uses linear-time SVD \cite{Drineas2006LTSVD} to solve the
partial SVD in each iteration) on randomly generated low-rank and
sparse matrices. \item Compare with factorization based solver
(e.g., LMaFit \cite{shen-MatrixFac}) on recovering either randomly
generated or deterministic
 low-rank matrix from its sum with a random sparse matrix.
\item Compare with random projection based solver (e.g., random
projection \cite{mu-2011-RandomProj}) on recovering randomly
generated low-rank and sparse matrices.
\end{enumerate}

In the experiments synthetic data, we generate random test data in
the following way: an $m \times m$ observed data matrix
$\mathbf{M}$ is synthesized as the sum of a low-rank matrix
$\mathbf{L}_0$ and a sparse matrix $\mathbf{S}_0$. The rank $r$
matrix $\mathbf{L}_0$ is generated as a product of two $m \times
r$ matrices whose entries are i.i.d. Gaussian random variables
with zero mean and unit variance. The matrix $\mathbf{S}_0$ is
generated as a sparse matrix whose support is chosen uniformly at
random, and whose $p$ non-zero entries are i.i.d. uniformly in
$[-500,500]$. The rank ratio and sparsity ratio are denoted as
$\rho_r = r/m$ and $\rho_s = p/m^2$, respectively.

\subsubsection{$l_1$ Filtering vs. Classic Convex Optimization}\label{sec:adm}

Firstly, we compare our approach with ADM on the whole
matrix\footnote{The Matlab code of ADM is provided by the authors
of (\cite{lin2009augmented}) and all the parameters in this code
are set to their default values.}, which we call the standard ADM,
and its variation, which uses linear-time SVD (LTSVD)\footnote{The
Matlab code of linear-time SVD is available in the FPCA package at
http://www.columbia.edu/$\sim$sm2756/FPCA.htm.} for solving the
partial SVD, hence we call the LTSVD ADM. We choose these two
approaches because the standard ADM is known to be the most
efficient classic convex optimization algorithm to solve PCP
exactly and LTSVD ADM has a linear time cost in solving
SVD\footnote{However, LTSVD ADM is still of $O(rmn)$ complexity as
it involves matrix-matrix multiplication in each iteration. See
also Section~\ref{sec:previous_work}.}. For LTSVD ADM, in each
time to compute the partial SVD we uniformly oversample $5r$
columns of the data matrix without replacement. Such an
oversampling rate is important for ensuring the numerical accuracy
of LTSVD ADM at high probability. For all methods in comparison,
the stopping criterion is
$\|\mathbf{M}-\mathbf{L}^*-\mathbf{S}^*\|_{F}/\|\mathbf{M}\|_F
\leq 10^{-7}$.

Table~\ref{tab:toy} shows the detailed comparison among the three
methods, where
$\mbox{RelErr}=\|\mathbf{L}^*-\mathbf{L}_0\|_F/\|\mathbf{L}_0\|_F$
is the relative error to the true low-rank matrix $\mathbf{L}_0$.
It is easy to see that our $l_1$ filtering approach has the
highest numerical accuracy and is also much faster than the
standard ADM and LTSVD ADM. Although LTSVD ADM is faster than the
standard ADM, its numerical accuracy is the lowest among the three
methods because it is probabilistic.
\begin{table*}[th]
\begin{center}
\caption{Comparison among the standard ADM (S-ADM for short),
LTSVD ADM (L-ADM for short) and $l_1$ filtering method ($l_1$ for
short) on the synthetic data.
We present CPU time (in seconds) and the numerical
accuracy of tested algorithms. $\mathbf{L}_0$ and $\mathbf{S}_0$
are the ground truth and $\mathbf{L}^*$ and $\mathbf{S}^*$ are the
solution computed by different methods. For the $l_1$ filtering
method, we report its computation time as $t = t_1 + t_2$, where
$t$, $t_1$ and $t_2$ are the time for total computation, seed
matrix recovery and $l_1$ filtering, respectively.}\label{tab:toy}
\begin{tabular}{|c||c|c|c|c|c|c|c|}\hline
Size ($m$) & Method & RelErr & $\rank(\mathbf{L}^*)$ &
$\|\mathbf{L}^*\|_*$ & $\|\mathbf{S}^*\|_{l_0}$ &
$\|\mathbf{S}^*\|_{l_1}$ & Time\\\hline\hline
 \multirow{4}{*}{2000} &
\multicolumn{7}{|c|}{$\rank(\mathbf{L}_0)=20$, \ \
$\|\mathbf{L}_0\|_*=39546 $, \ \ $\|\mathbf{S}_0\|_{l_0}=40000$, \
\ $\|\mathbf{S}_0\|_{l_1}= 998105$}\\\cline{2-8}
& S-ADM          & 1.46 $\times 10^{-8}$ & 20 & 39546  & 39998  & 998105 & 84.73 \\
& L-ADM       & 4.72 $\times 10^{-7}$  & 20 & 39546 & 40229  & 998105 & 27.41\\
& $l_1$ & 1.66 $\times 10^{-8}$  & 20 & 39546 & 40000  & 998105 &
\textbf{5.56 = 2.24 + 3.32}\\\hline\hline \multirow{4}{*}{5000} &
\multicolumn{7}{|c|}{$\rank(\mathbf{L}_0)=50$, \ \
$\|\mathbf{L}_0\|_*=249432$, \ \ $\|\mathbf{S}_0\|_{l_0}=250000$,
\ \ $\|\mathbf{S}_0\|_{l_1}=6246093$}\\\cline{2-8}
& S-ADM          & 7.13 $\times 10^{-9}$ & 50 & 249432 & 249995  & 6246093 & 1093.96 \\
& L-ADM       & 4.28 $\times 10^{-7}$ & 50 & 249432 & 250636  & 6246158 & 195.79\\
& $l_1$ & 5.07 $\times 10^{-9}$ & 50 & 249432 &
250000  & 6246093 & \textbf{42.34=19.66 +
22.68}\\\hline\hline
 \multirow{4}{*}{10000} &
\multicolumn{7}{|c|}{$\rank(\mathbf{L}_0)=100$, \ \
$\|\mathbf{L}_0\|_*=997153$, \ \ $\|\mathbf{S}_0\|_{l_0}=1000000$,
\ \ $\|\mathbf{S}_0\|_{l_1}=25004070$}\\\cline{2-8}
& S-ADM           & 1.23 $\times 10^{-8}$  & 100 & 997153 & 1000146  & 25004071 & 11258.51 \\
& L-ADM        & 4.26 $\times 10^{-7}$  & 100 & 997153 & 1000744  & 25005109 & 1301.83\\
& $l_1$ & 2.90 $\times 10^{-10}$ & 100 & 997153 &
1000023  & 25004071 & \textbf{276.54 = 144.38 +
132.16}\\\hline
\end{tabular}
\end{center}
\end{table*}

\begin{figure}[ht]
\centering
\begin{tabular}{c@{\extracolsep{0.2em}}c}
\includegraphics[width=0.22\textwidth,
keepaspectratio]{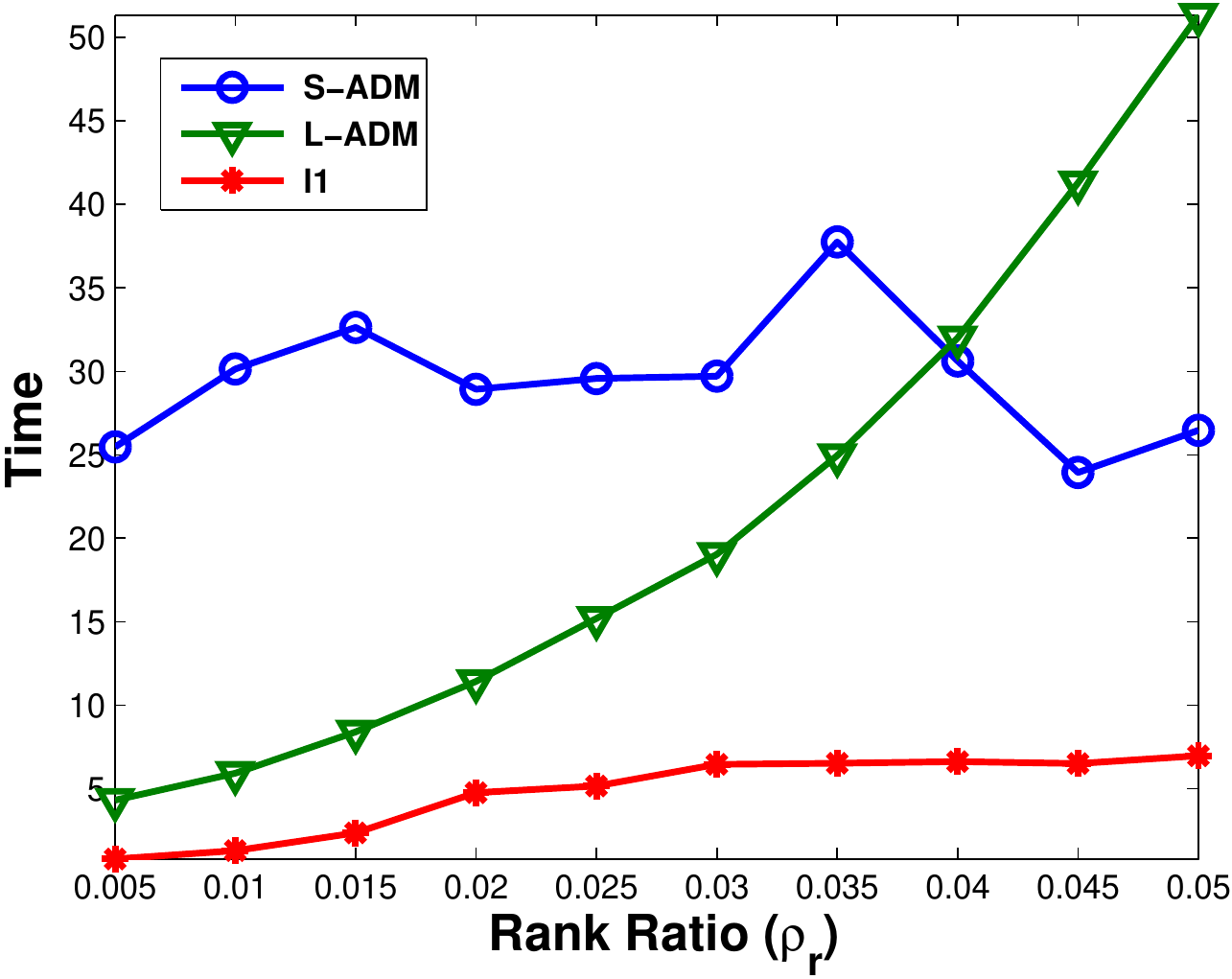}
&\includegraphics[width=0.22\textwidth,
keepaspectratio]{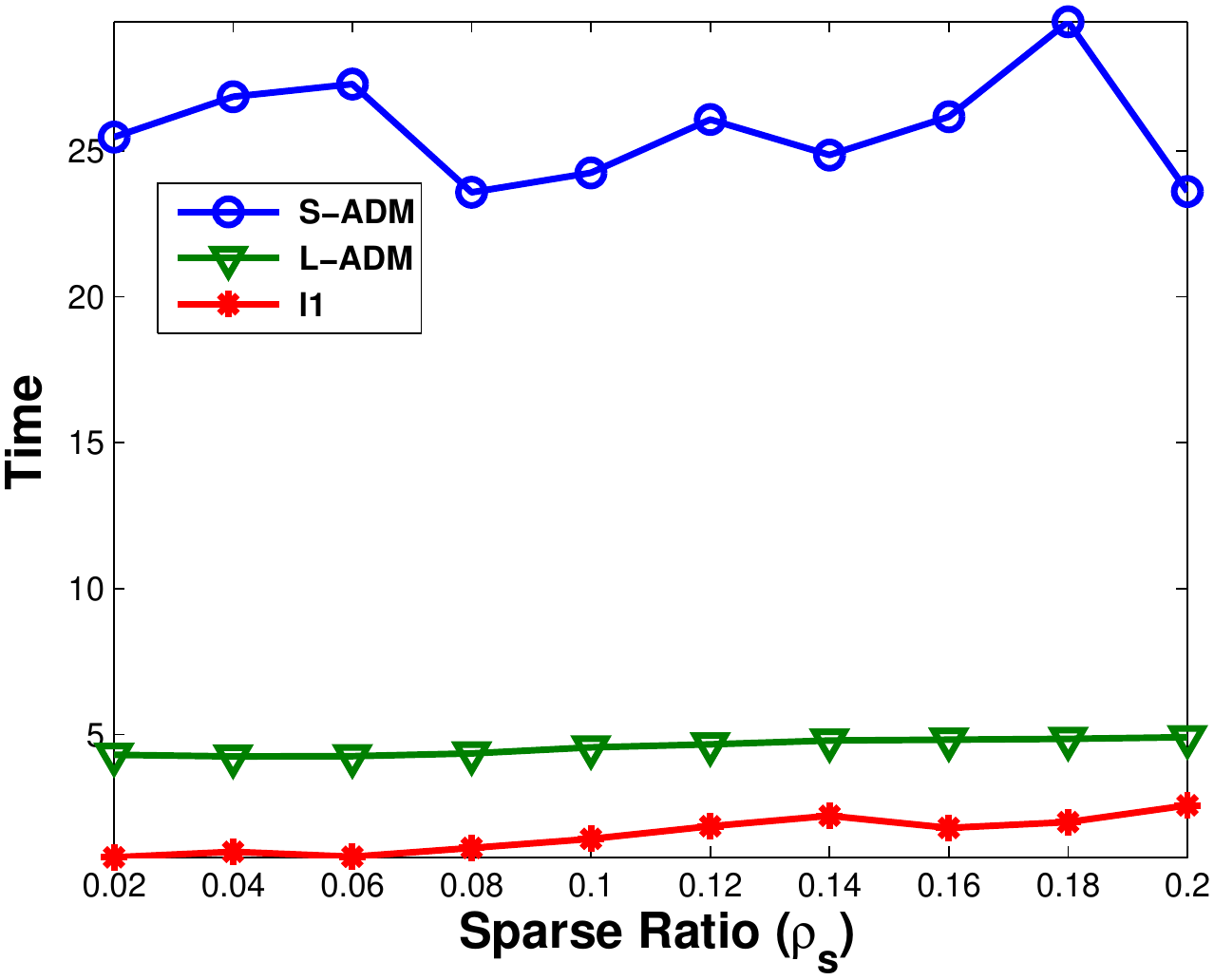}\\
(a) & (b) \\
\end{tabular}
\caption{Performance of the standard ADM (S-ADM for short), LTSVD
ADM (L-ADM for short) and $l_1$ filtering under different rank
ratios $\rho_r$ and sparsity ratios $\rho_s$, where the matrix
size is $1000\times 1000$. The $x$-axis represents the rank ratio
(a) or sparsity ratio (b). The $y$-axis represents the CPU time
(in seconds).}\label{fig:rank_ratio}
\end{figure}

We also present in Figure~\ref{fig:rank_ratio} the CPU times of
the three methods when the rank ratio $\rho_r$ and sparsity ratio
$\rho_s$ increases, respectively. The observed matrices are
generated using the following parameter settings: $m = 1000$, vary
$\rho_r$ from 0.005 to 0.05 with fixed $\rho_s = 0.02$ and vary
$\rho_s$ from 0.02 to 0.2 with fixed $\rho_r = 0.005$. It can be
seen from Figure~\ref{fig:rank_ratio} (a) that LTSVD ADM is faster
than the standard ADM when $\rho_r < 0.04$. However, the computing
time of LTSVD ADM grows quickly when $\rho_r$ increases. It even
becomes slower than the standard ADM when $\rho_r \geq 0.04$. This
is because LTSVD cannot guarantee the accuracy of partial SVD in
each iteration. So its number of iterations is larger than that of
the standard ADM. In comparison, the time cost of our $l_1$
filtering method is much less than the other two methods for all
the rank ratios. However, when $\rho_r$ further grows the
advantage of $l_1$ filtering will be lost quickly, because $l_1$
filtering has to compute the PCP on the $(s_r r)\times (s_c
r)=(10r)\times(10r)$ submatrix $\mathbf{M}^s$. In contrast,
Figure~\ref{fig:rank_ratio} (b) indicates that the CPU time of
these methods grows very slowly with respect to the sparsity
ratio.


\subsubsection{$l_1$ Filtering vs. Factorization Method}

We then compare the proposed $l_1$ filtering with a factorization
method (i.e., LMaFit\footnote{The Matlab code of LMaFit is
provided by the authors of (\cite{shen-MatrixFac}) and all the
parameters in this code are set to their default values.}) on
solving (\ref{eq:PCP}). To test the ability of these algorithms in
coping with corruptions with large magnitude, we multiply a scale
$\mathbf{\sigma}$ to the sparse matrix, i.e., $\mathbf{M} =
\mathbf{L}_0 + \sigma\mathbf{S}_0$. We fix other parameters of the
data ($m=1000$, $r = 0.01m$ and $\rho_s = 0.01$) and vary the
scale parameter $\sigma$ from 1 to 10 to increase the magnitude of
the sparse errors.

\begin{figure*}[ht]
\centering
\begin{tabular}{c@{\extracolsep{0.2em}}c@{\extracolsep{0.2em}}c@{\extracolsep{0.2em}}c}
\includegraphics[width=0.22\textwidth,
keepaspectratio]{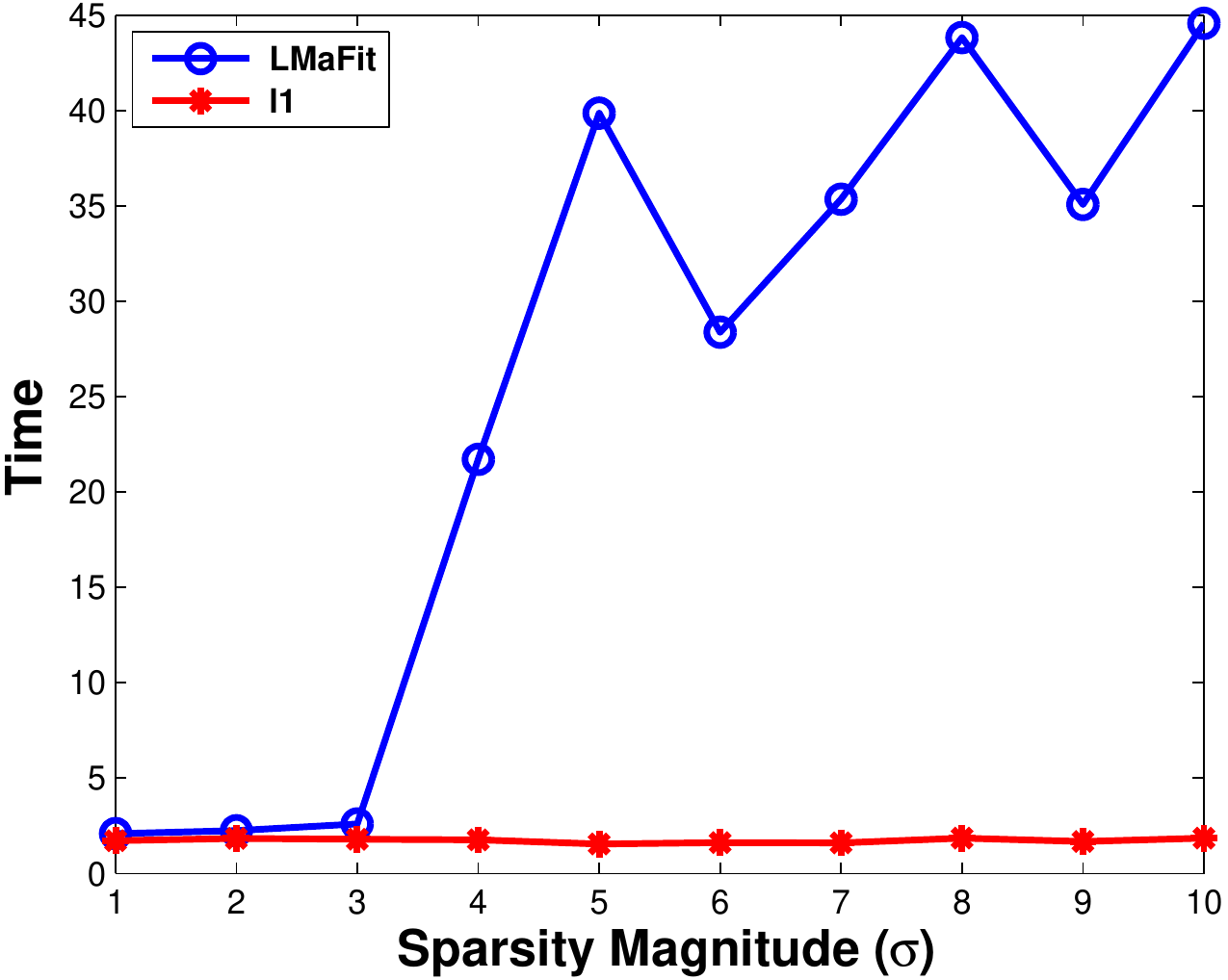}
&\includegraphics[width=0.22\textwidth,
keepaspectratio]{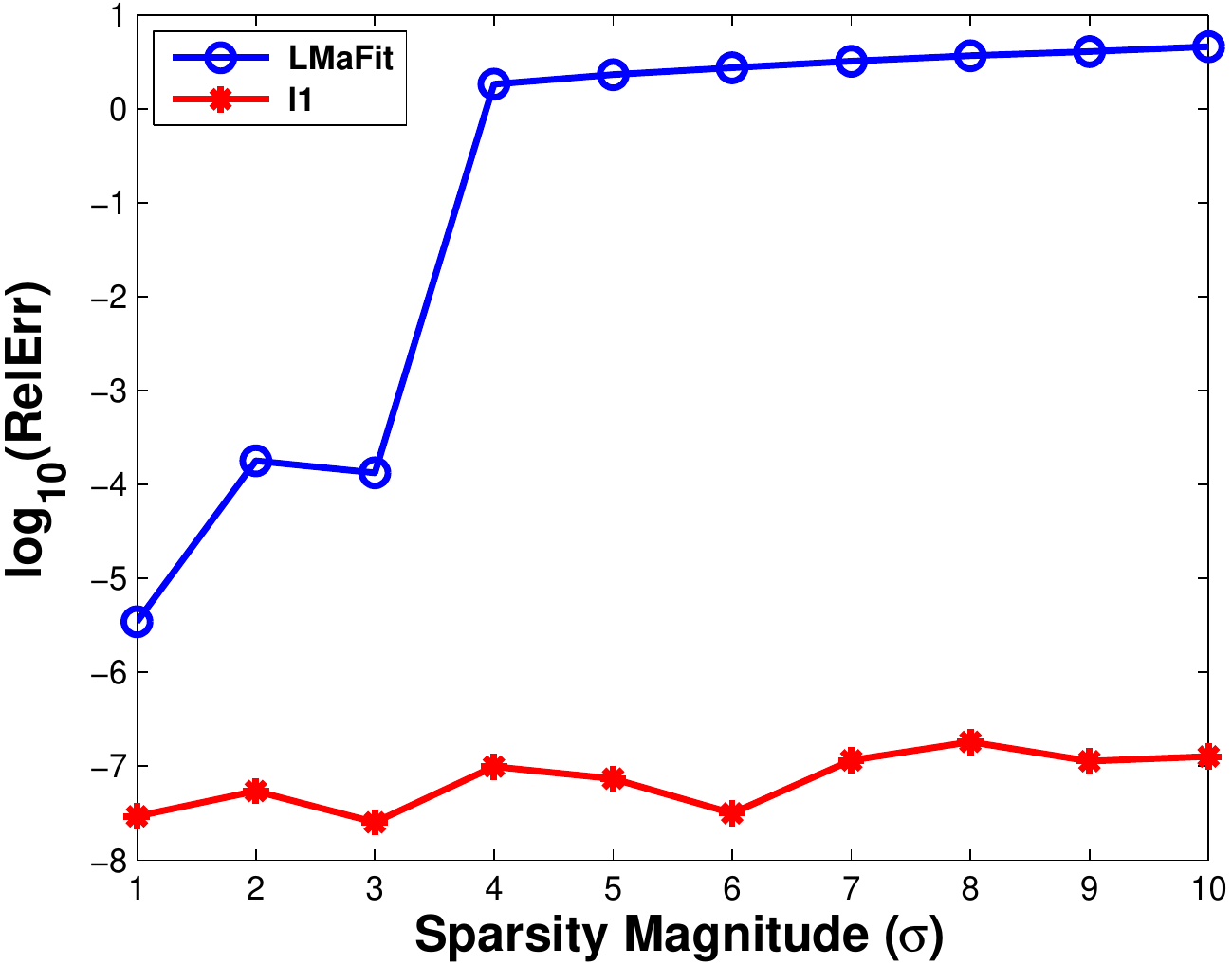}
&\includegraphics[width=0.22\textwidth,
keepaspectratio]{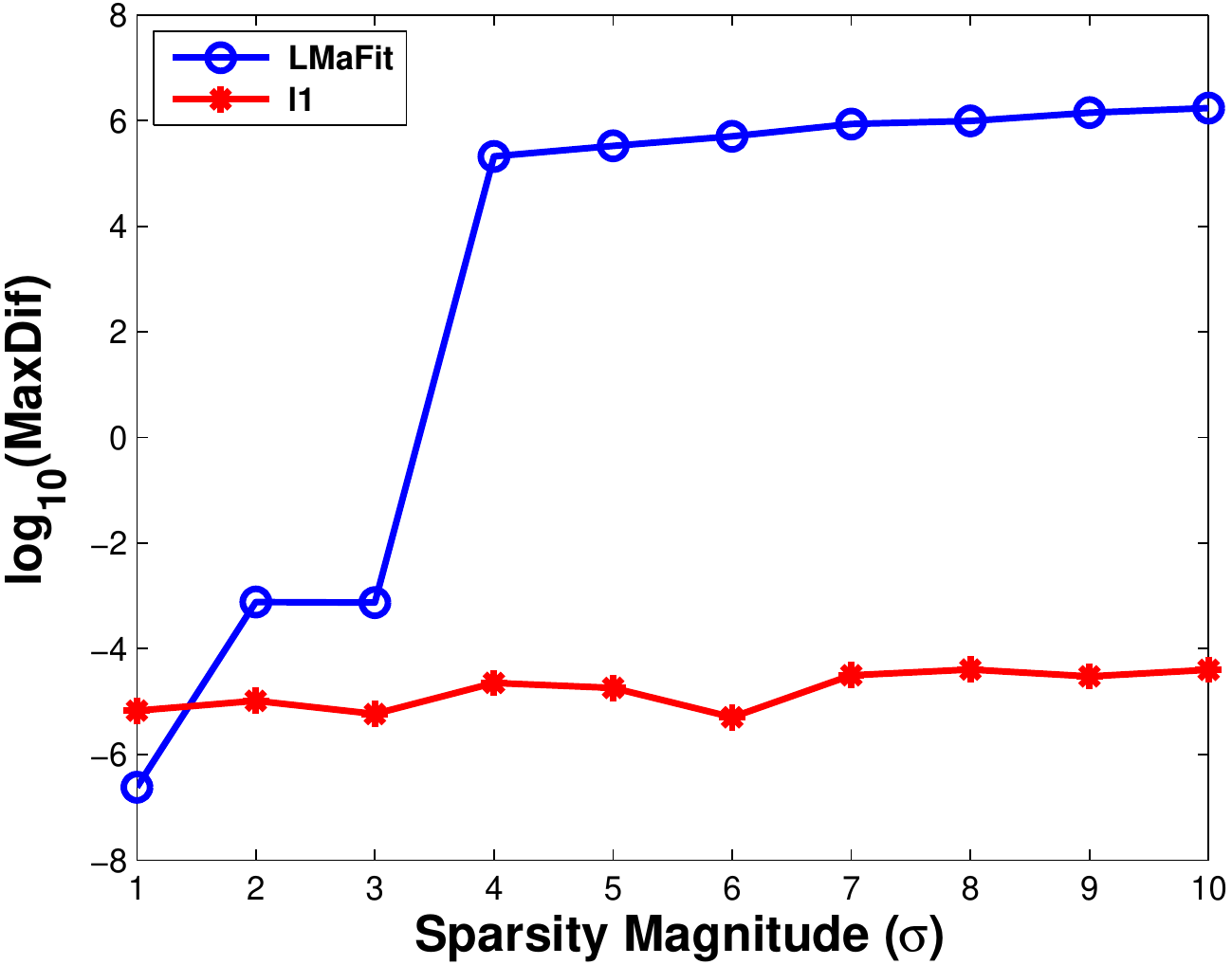}
&\includegraphics[width=0.22\textwidth,
keepaspectratio]{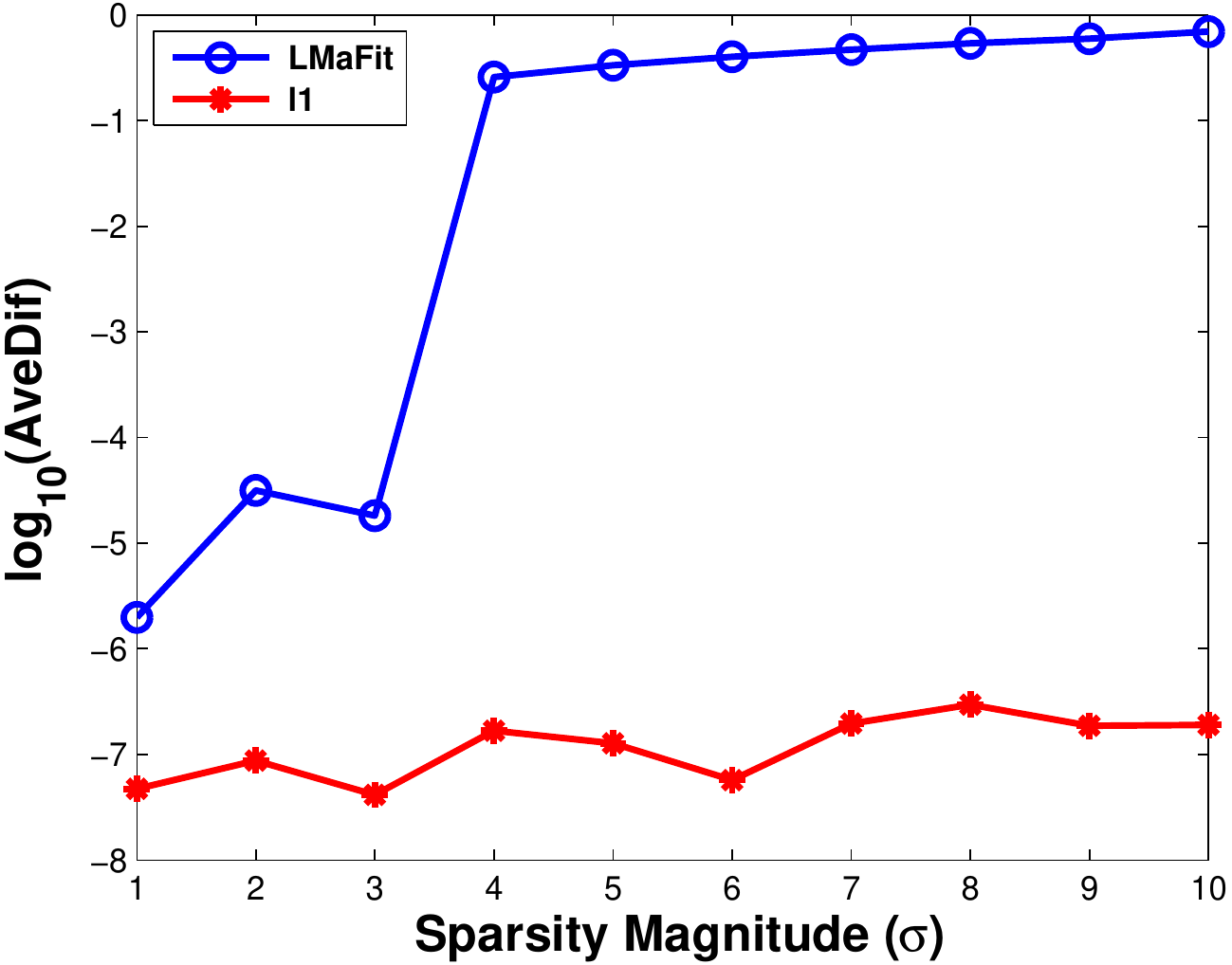}\\
(a) & (b) & (c) & (d)
\end{tabular}
\caption{Performance of LMaFit and $l_1$ filtering under different
sparsity magnitudes ($\sigma\in [1,10]$). The $x$-axes represent
the sparsity magnitudes and the $y$-axes represent the CPU time
(in seconds) (a), ``RelErr'' (b), ``MaxDif'' (c) and ``AveDif''
(d) in log scale, respectively.}\label{fig:factorization}
\end{figure*}

The computational comparisons are presented in
Figure~\ref{fig:factorization}. Besides the CPU time and relative
error, we also measure the quality of the recovered $\mathbf{L}^*$
by its maximum difference (MaxDif) and average difference (AveDif)
to the true low-rank matrix $\mathbf{L}_0$, which are respectively
defined as $\mbox{MaxDif} = \max(|\mathbf{L}^*-\mathbf{L}_0|)$ and
$\mbox{AveDif} =(\sum_{ij}|\mathbf{L}^*-\mathbf{L}_0|)/m^2$. One
can see that the performance of LMaFit dramatically decreases when
$\sigma \geq 3$. This experiment suggests that the factorization
method fails when the sparse matrix dominates the low-rank one in
magnitude. This is because a sparse matrix with large magnitudes
makes rank estimation difficult or impossible for LMaFit. Without
a correct rank, the low-rank matrix cannot be recovered exactly.
In comparison, our $l_1$ filtering always performs well on the
test data.

\begin{figure*}[ht]
\centering
\begin{tabular}{c@{\extracolsep{0.5em}}c@{\extracolsep{0.5em}}c@{\extracolsep{3em}}c}
\includegraphics[width=0.18\textwidth,
keepaspectratio]{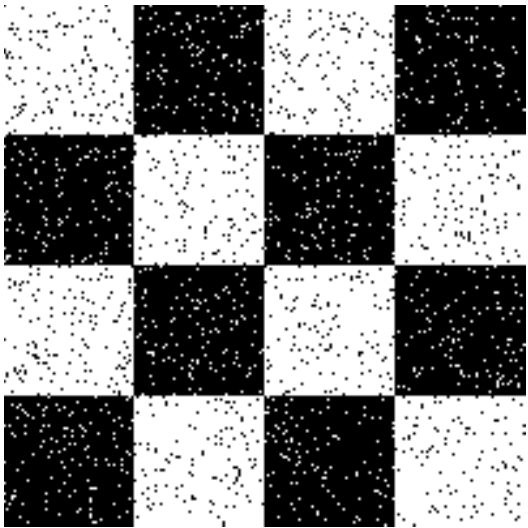}
&\includegraphics[width=0.18\textwidth,
keepaspectratio]{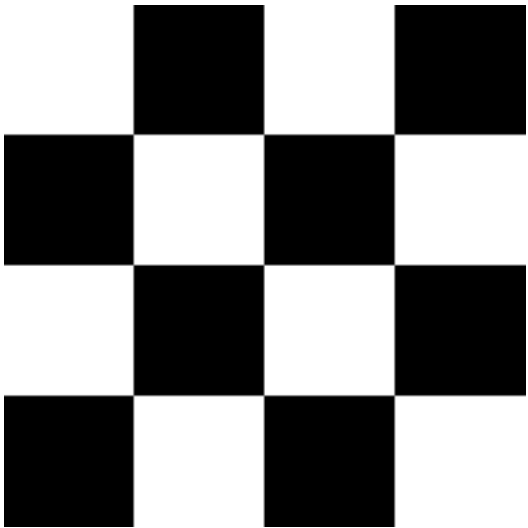}
&\includegraphics[width=0.18\textwidth,
keepaspectratio]{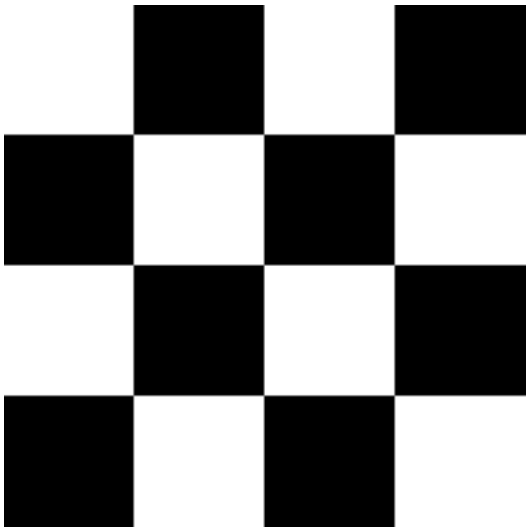}
&\includegraphics[width=0.24\textwidth,
keepaspectratio]{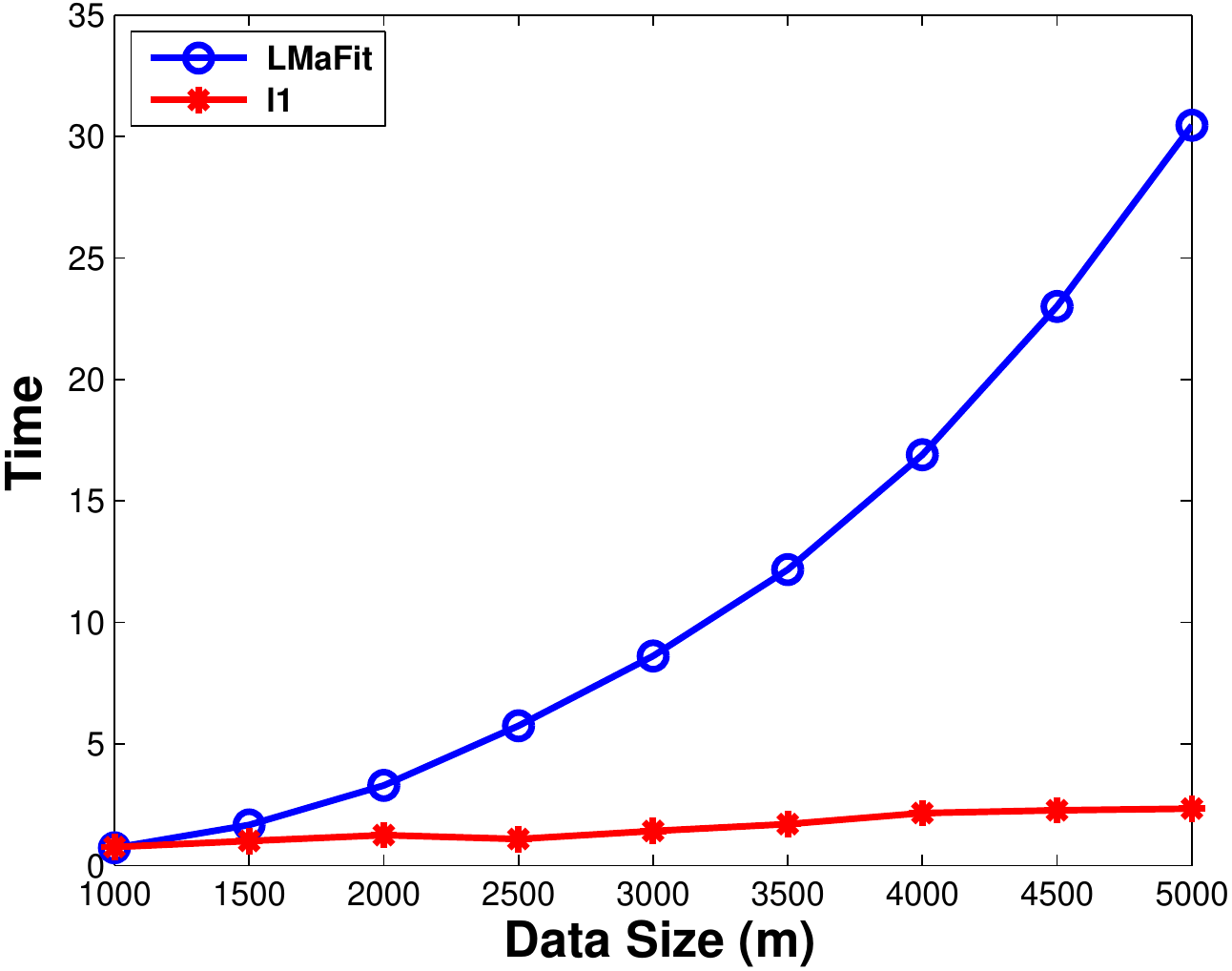}\\
(a) Corrupted & (b) LMaFit & (c) $l_1$ & (d) CPU Time
\end{tabular}
\caption{Recovery results for ``checkerboard''. (a) is the image
corrupted by 10$\%$ impulsive noise. (b) is the image recovered by
LMaFit. (c) is the image recovered by $l_1$ filtering ($l_1$). (d)
CPU time (in seconds) vs. data size ($m \in
[1000,5000]$).}\label{fig:check}
\end{figure*}

In the following, we consider the problem of recovering
\emph{deterministic} low-rank matrix from corruptions. We generate
an $m \times m$ ``checkerboard'' image (see
Figure~\ref{fig:check}), whose rank is 2, and corrupt it by adding
10$\%$ impulsive noise to it. The corruptions (nonzero entries of
the sparse matrix) are sampled uniformly at random. The image size
$m$ ranges from 1000 to 5000 with an increment 500.

The results for this test are shown in Figure~\ref{fig:check},
where the first image is the corrupted checkerboard image, the
second image is recovered by LMaFit and the third by $l_1$
filtering. A more complete illustration for this test can be seen
from Figure~\ref{fig:check}(d), where the CPU time corresponding
to all tested data matrix sizes are plotted. It can be seen that
the images recovered by LMaFit and $l_1$ filtering are visually
comparable in quality. The speeds of these two methods are very
similar when the data size is small, while $l_1$ filtering runs
much faster than LMaFit when the matrix size increases. This
concludes that our approach has significant speed advantage over
the factorization method on large scale data sets.

\subsubsection{$l_1$ Filtering vs. Compressed Optimization}

Now we compare $l_1$ filtering with a compressed optimization
method (i.e., random projection\footnote{The Matlab code of random
projection (RP) is provided by the author of
(\cite{mu-2011-RandomProj}) and all the parameters in this code
are set to their default values.}). This experiment is to study
the performance of these two methods with respect to the rank of
the matrix and the data size. The parameters of the test matrices
are set as follows: $\rho_s = 0.01$, $\rho_r$  varying from 0.05
to 0.15 with fixed $m = 1000$, and $m$ varying from 1000 to 5000
with fixed $\rho_r = 0.05$. For the dimension of the projection
matrix (i.e., $p$), we set it as $p = 2r$ for all the experiments.

\begin{figure}[ht]
\centering
\begin{tabular}{c@{\extracolsep{0.2em}}c}
\includegraphics[width=0.22\textwidth,
keepaspectratio]{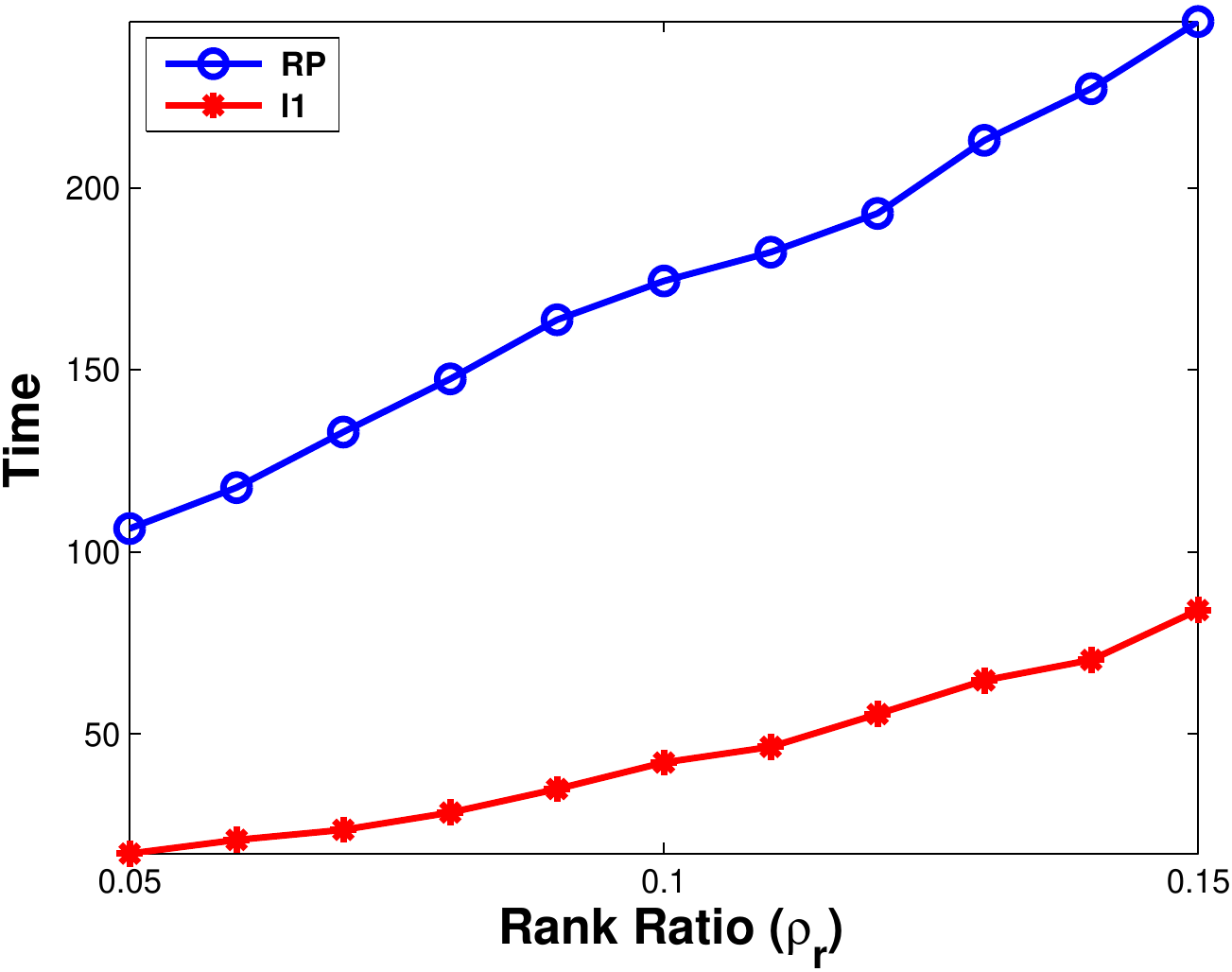}
&\includegraphics[width=0.22\textwidth,
keepaspectratio]{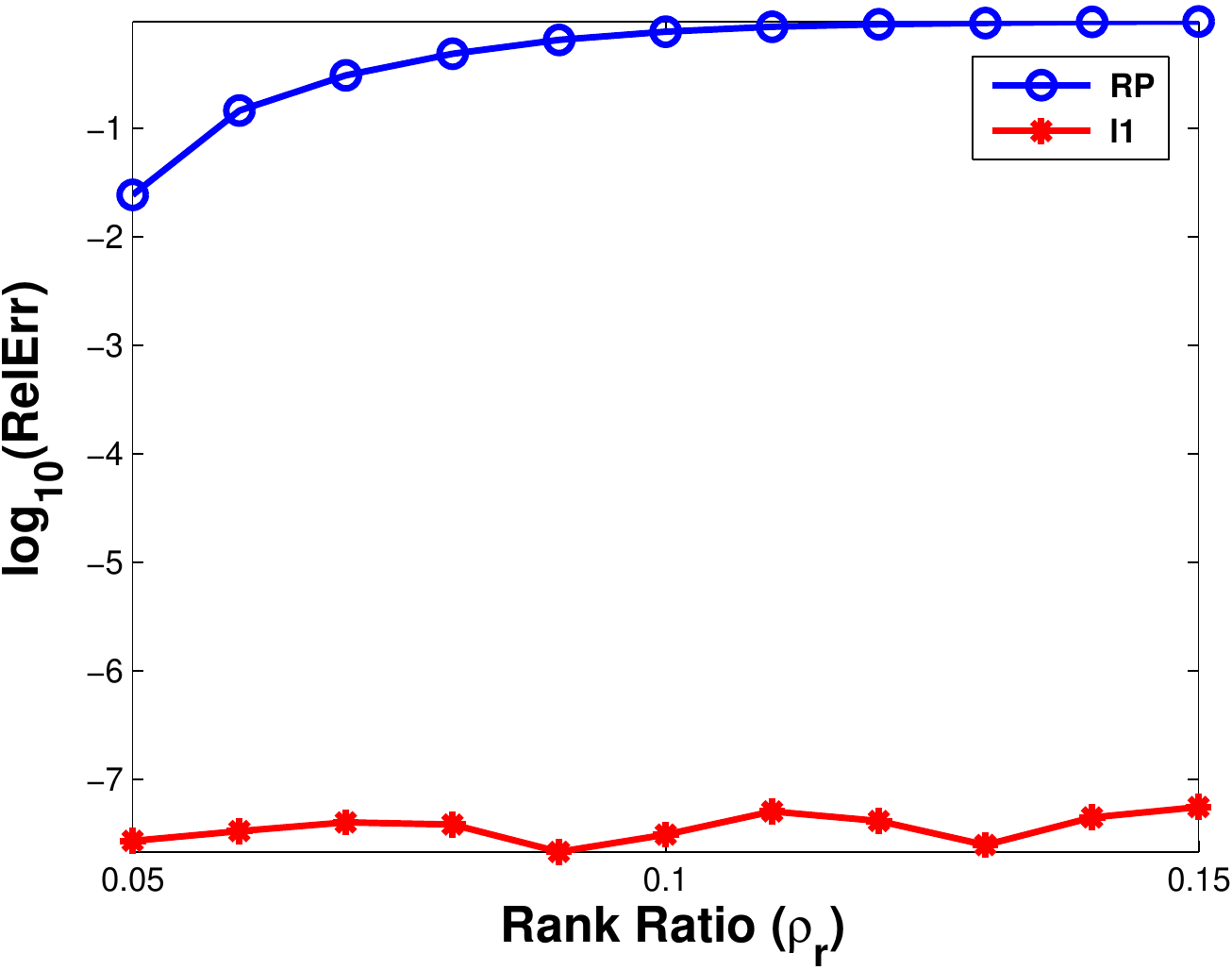}\\
(a) & (b)\\
\includegraphics[width=0.22\textwidth,
keepaspectratio]{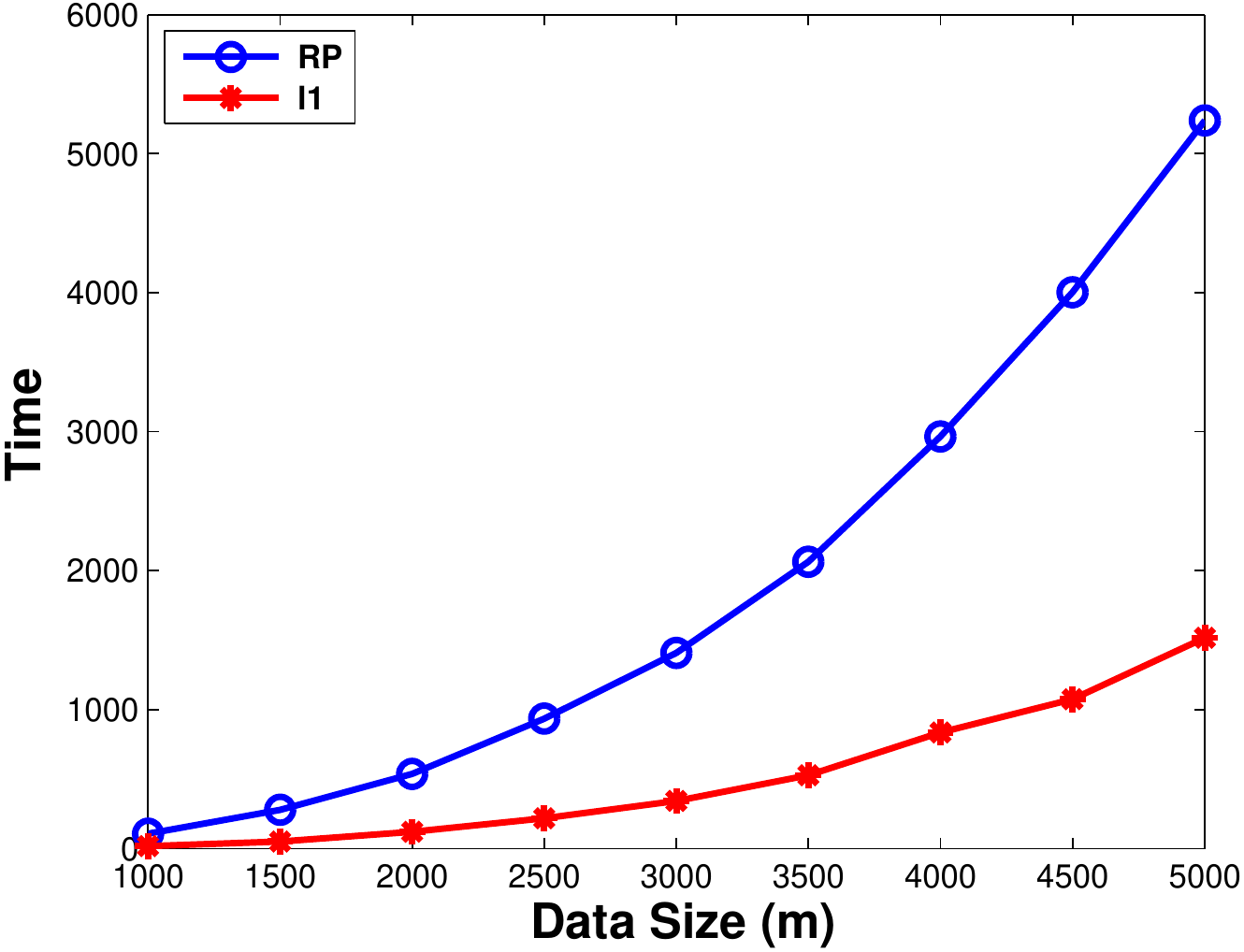}
&\includegraphics[width=0.22\textwidth,
keepaspectratio]{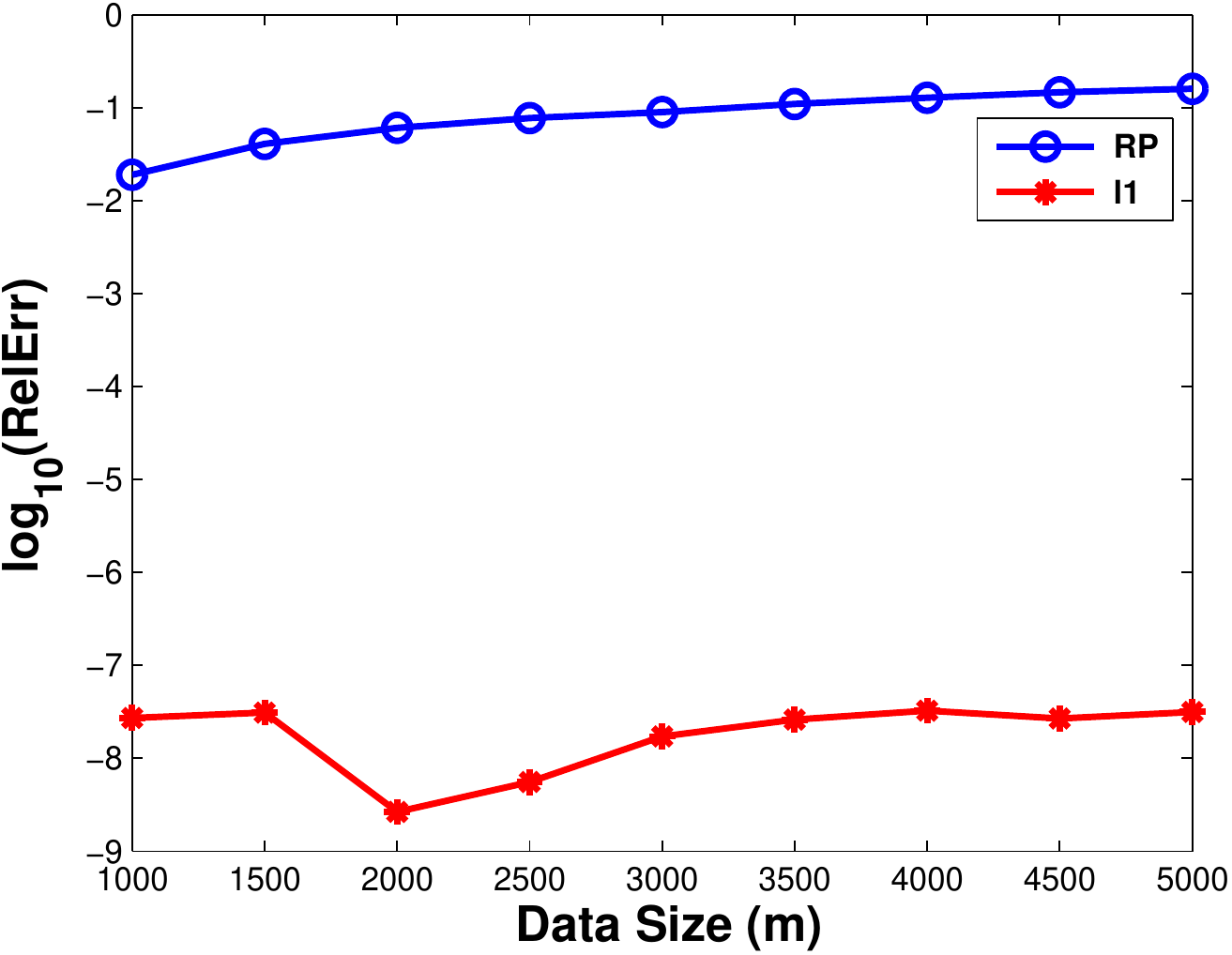}\\
(c) & (d) \\
\end{tabular}
\caption{Performance of random projection (RP for short) and $l_1$
filtering. (a)-(b) are the comparison under different rank ratios
($\rho_r\in [0.05,0.15]$). (c)-(d) are the comparison under
different data sizes ($m \in [1000,5000]$). In (a) and (c), the
$y$-axes are the CPU times (in seconds). In (b) and (d), the
$y$-axes are the relative errors in log
scale.}\label{fig:projection}
\end{figure}

As shown in Figure~\ref{fig:projection}, in all cases the speed
and the numerical accuracy of $l_1$ filtering are always much
higher than those of random projection. \comment{This is because
when fixing the data size, the complexity of $l_1$ filtering is
with respect to the rank of the data matrix $r$ while that of
random projection is with respect to the dimension of the
projection $p$ (note that $p \geq r$). In contrast, when fixing
the rank, the complexity of $l_1$ filtering and random projection
are respectively linear and quadratic with respect to the data
size.}

\subsection{Structure from Motion}
In this subsection, we apply $l_1$ filtering to a real world
vision application, namely structure from motion (SfM). The
problem of SfM is to automatically recover the 3D structure of an
object from a sequence of images of the object. Suppose that the
object is rigid, there are $F$ frames and $P$ tracked feature
points (i.e., $\mathbf{L}_0 = \begin{bmatrix}\mathbf{X} \\
\mathbf{Y}
\end{bmatrix}_{2F\times P}$), and the camera intrinsic parameters do not change. As shown
in~(\cite{Rao-2010-motion}), the trajectories of feature points
from a single rigid motion of the camera all lie in a liner
subspace of $\mathbb{R}^{2F}$, whose dimension is at most four
(i.e., $\rank(\mathbf{L}_0) \leq 4$). It has been shown that
$\mathbf{L}_0$ can be factorized as
$\mathbf{L}_0=\mathbf{A}\mathbf{B}$, where
$\mathbf{A}\in\mathbb{R}^{2F\times 4}$ recovers the rotations and
translations while the first three rows of $\mathbf{B}\in
\mathbb{R}^{4\times P}$ encode the relative 3D positions for each
feature point in the reconstructed object. However, when there
exist errors (e.g., occlusion, missing data or outliers) the
feature matrix is no longer of rank 4. Then recovering the full 3D
structure of the object can be posed as a low-rank matrix recovery
problem.

\begin{figure*}[ht]
\centering
\begin{tabular}{c@{\extracolsep{0.2em}}c@{\extracolsep{0.2em}}c@{\extracolsep{0.2em}}c@{\extracolsep{0.2em}}c@{\extracolsep{0.2em}}c@{\extracolsep{0.2em}}c}
\includegraphics[width=0.12\textwidth,
keepaspectratio]{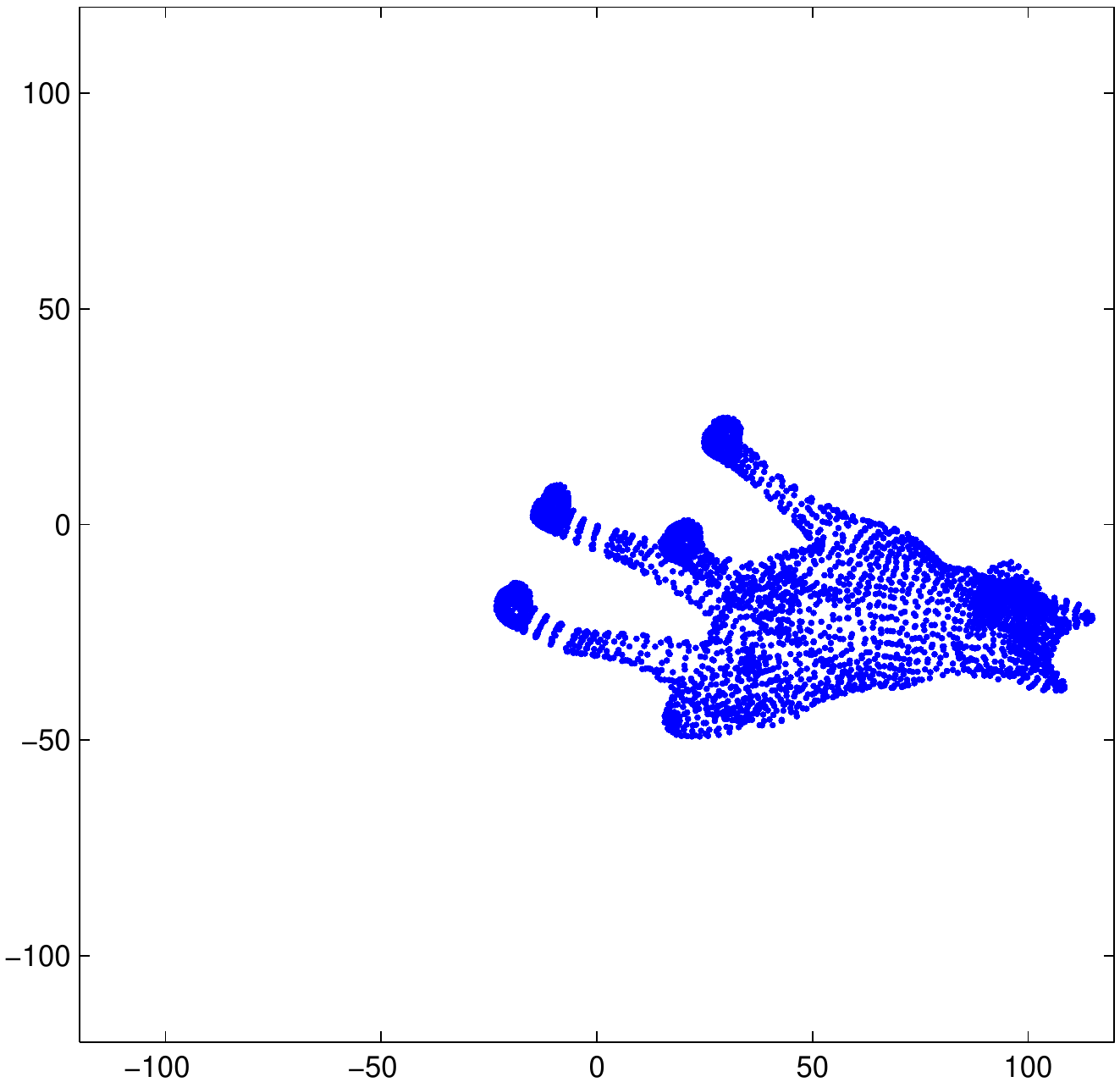}
&\includegraphics[width=0.12\textwidth,
keepaspectratio]{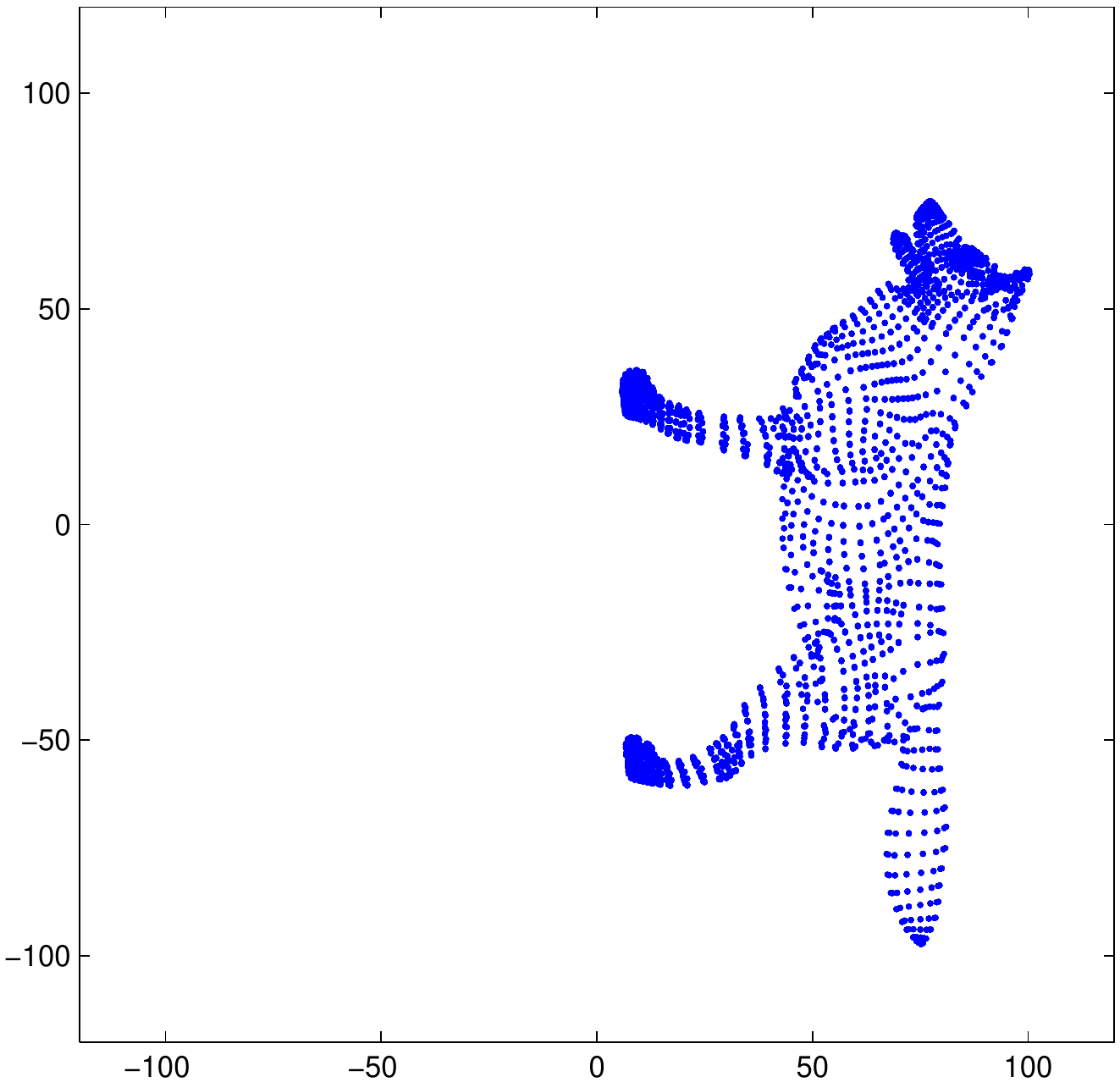}
&\includegraphics[width=0.12\textwidth,
keepaspectratio]{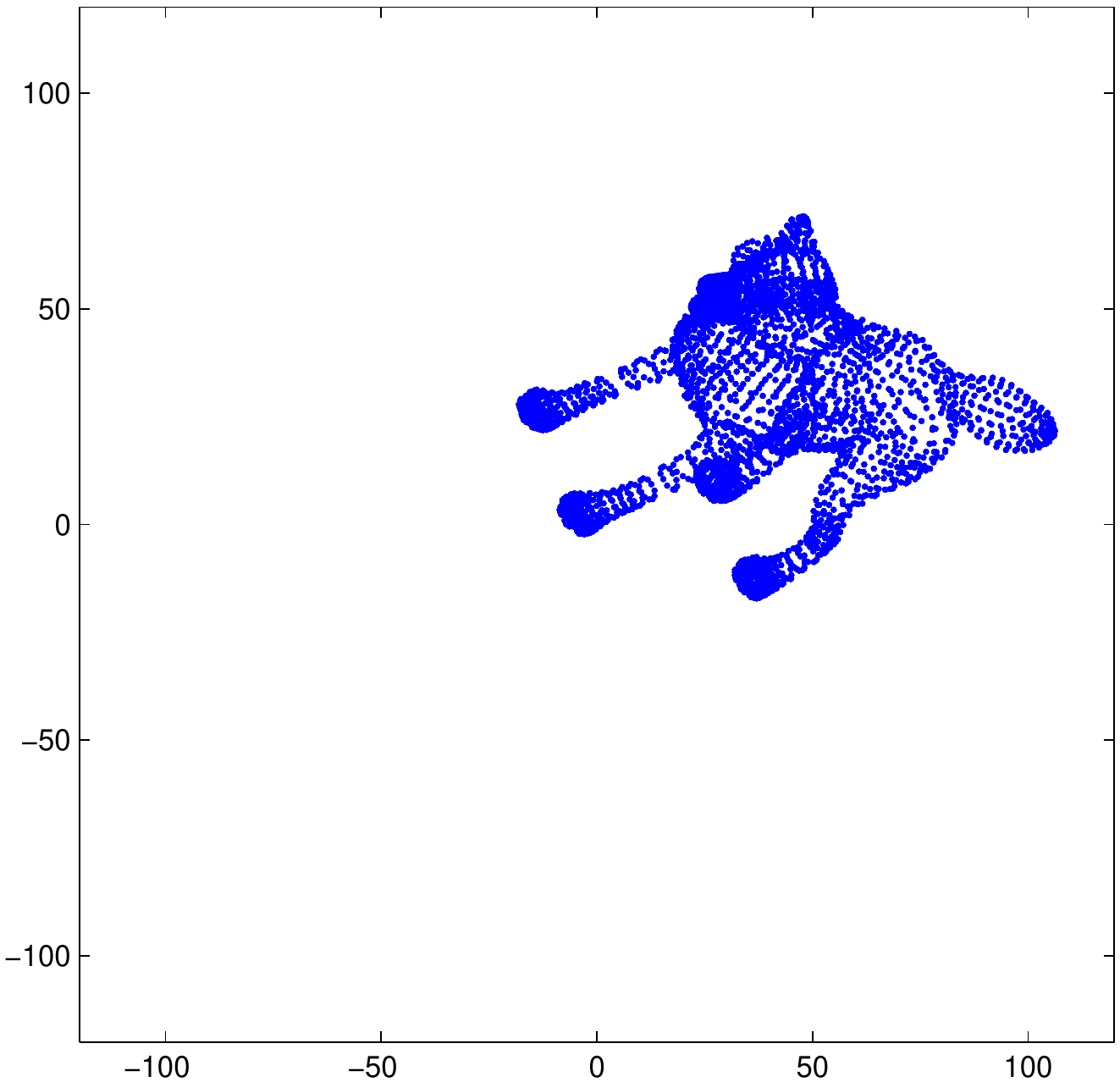}
&\includegraphics[width=0.12\textwidth,
keepaspectratio]{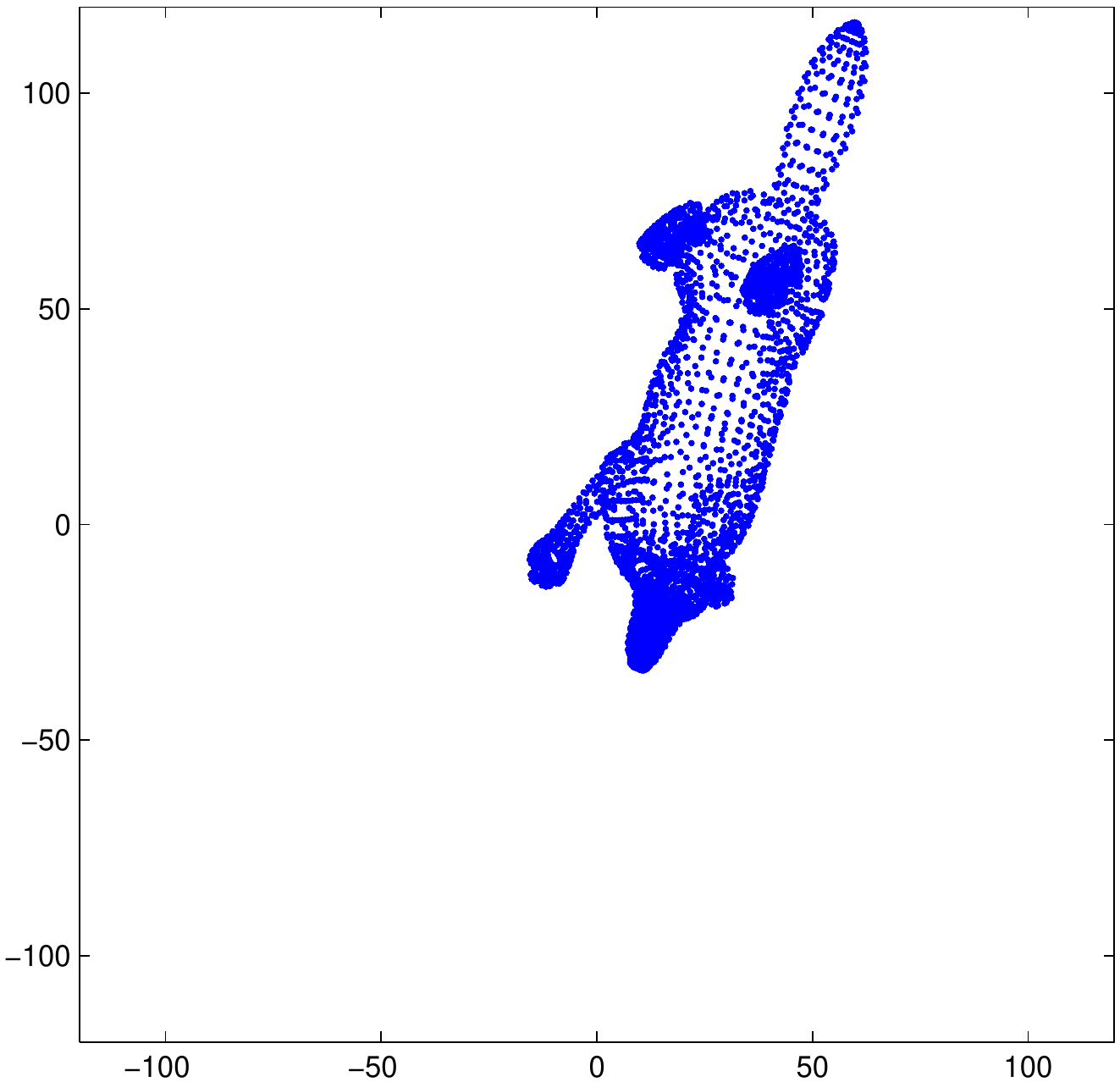}
&\includegraphics[width=0.12\textwidth,
keepaspectratio]{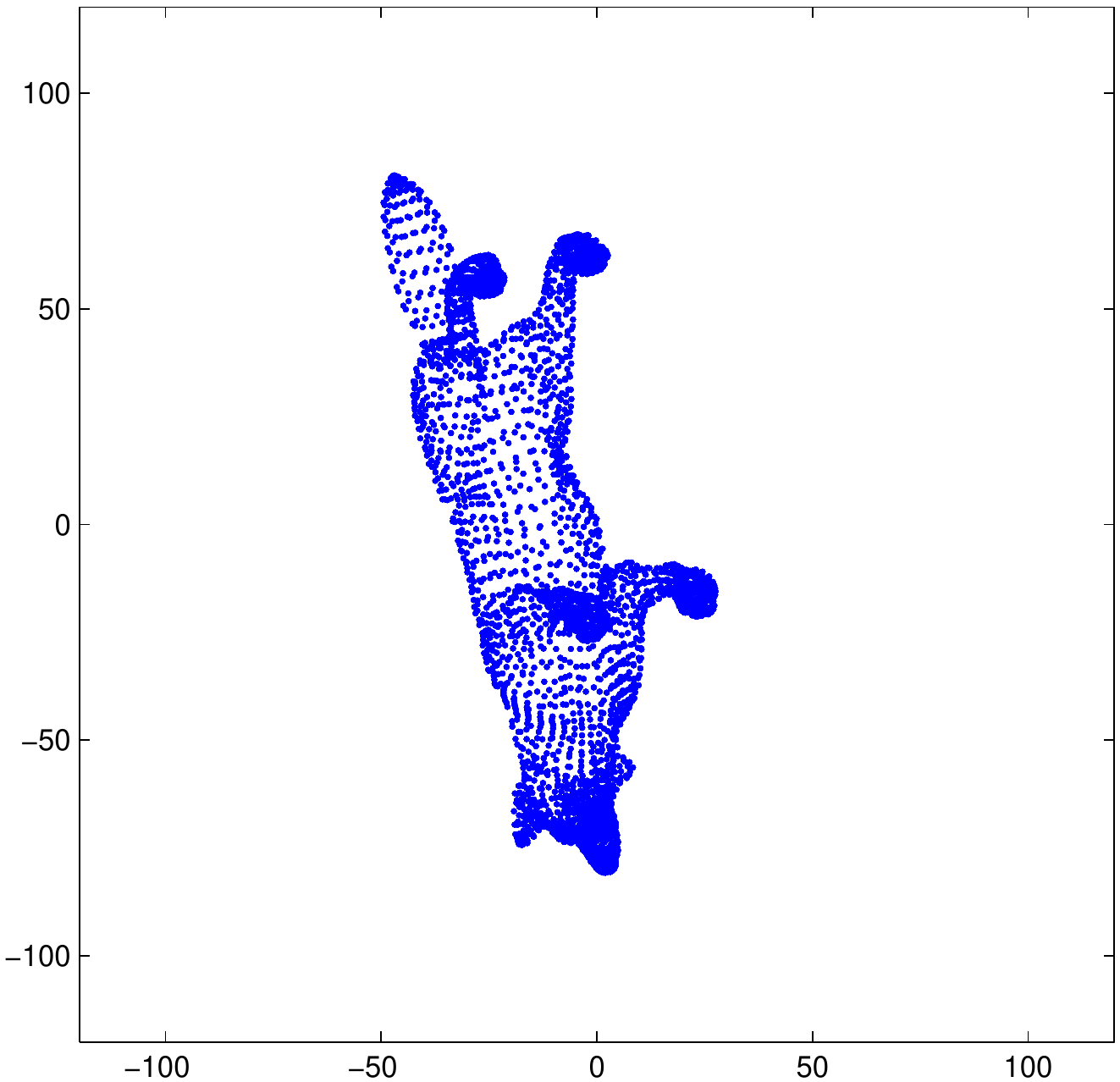}
&\includegraphics[width=0.12\textwidth,
keepaspectratio]{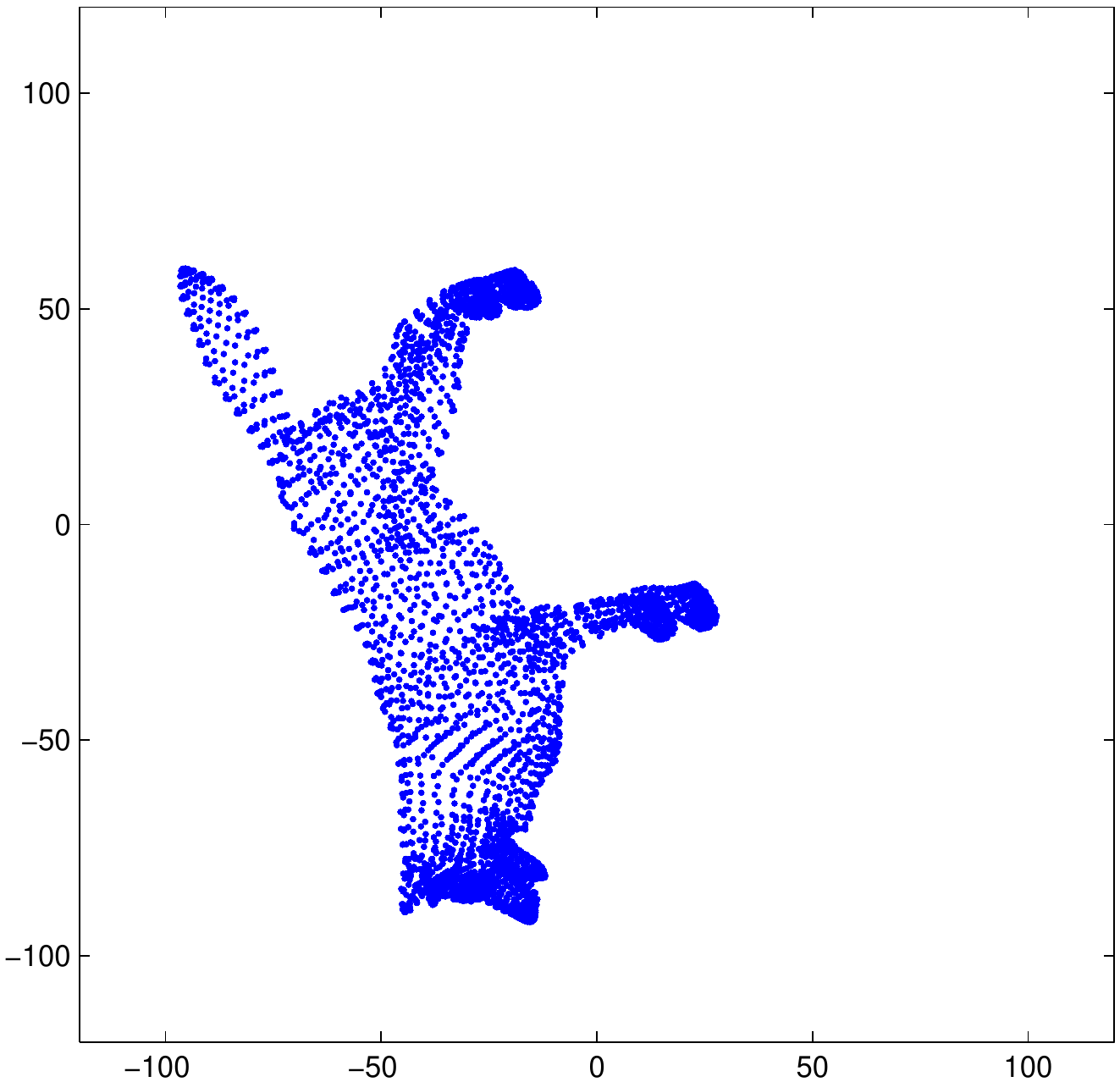}
&\includegraphics[width=0.12\textwidth,
keepaspectratio]{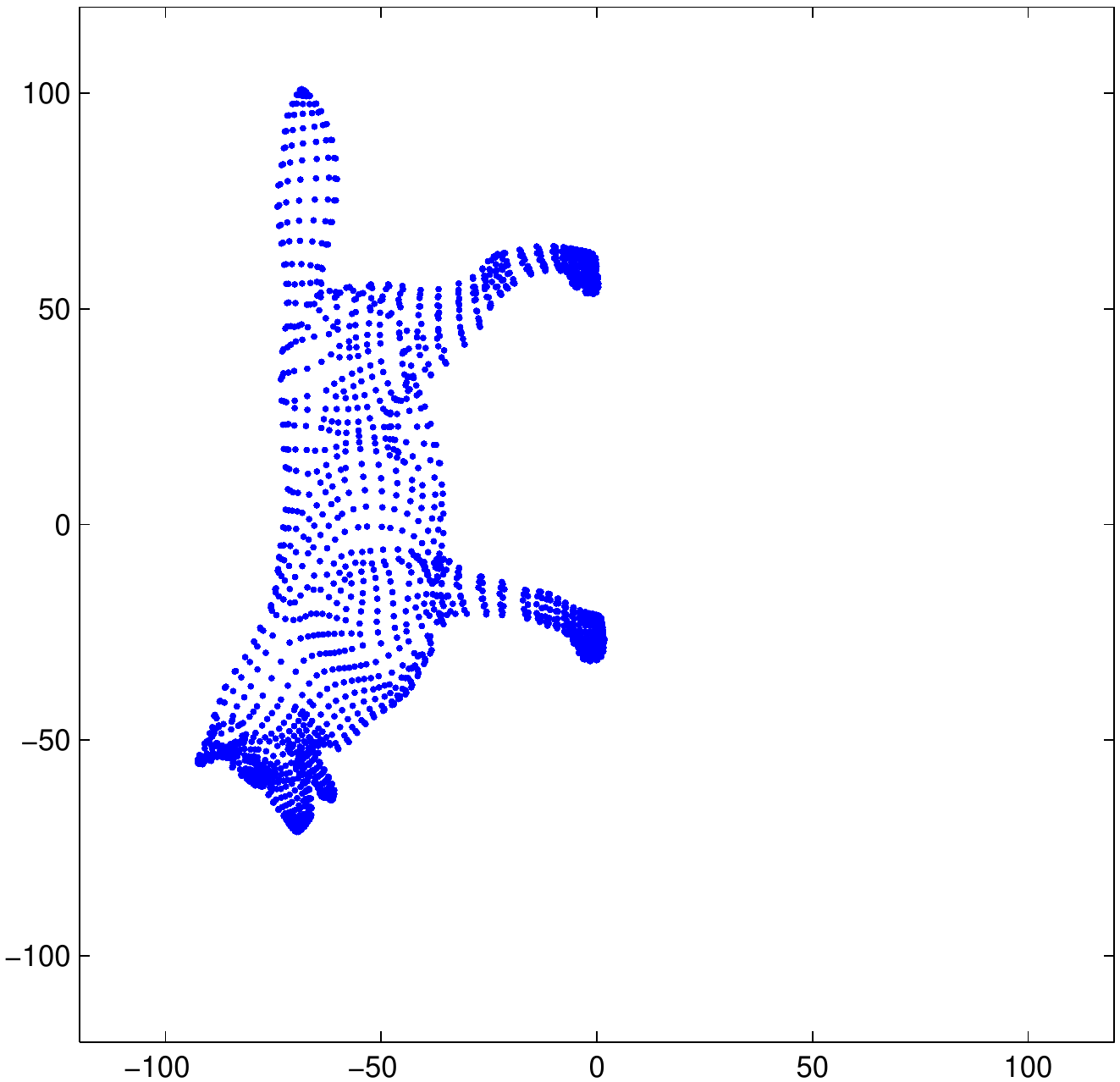}\\
\includegraphics[width=0.12\textwidth,
keepaspectratio]{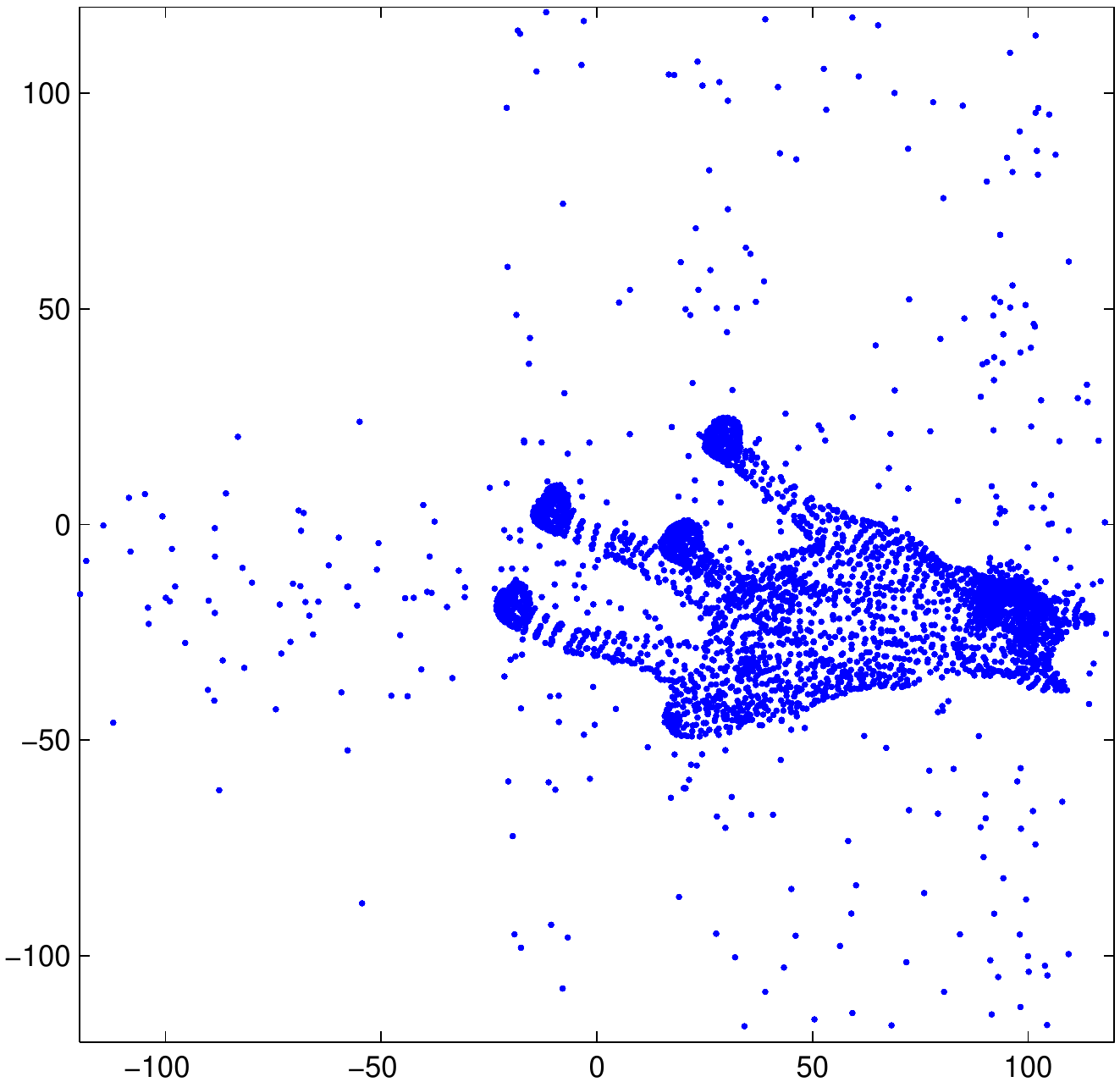}
&\includegraphics[width=0.12\textwidth,
keepaspectratio]{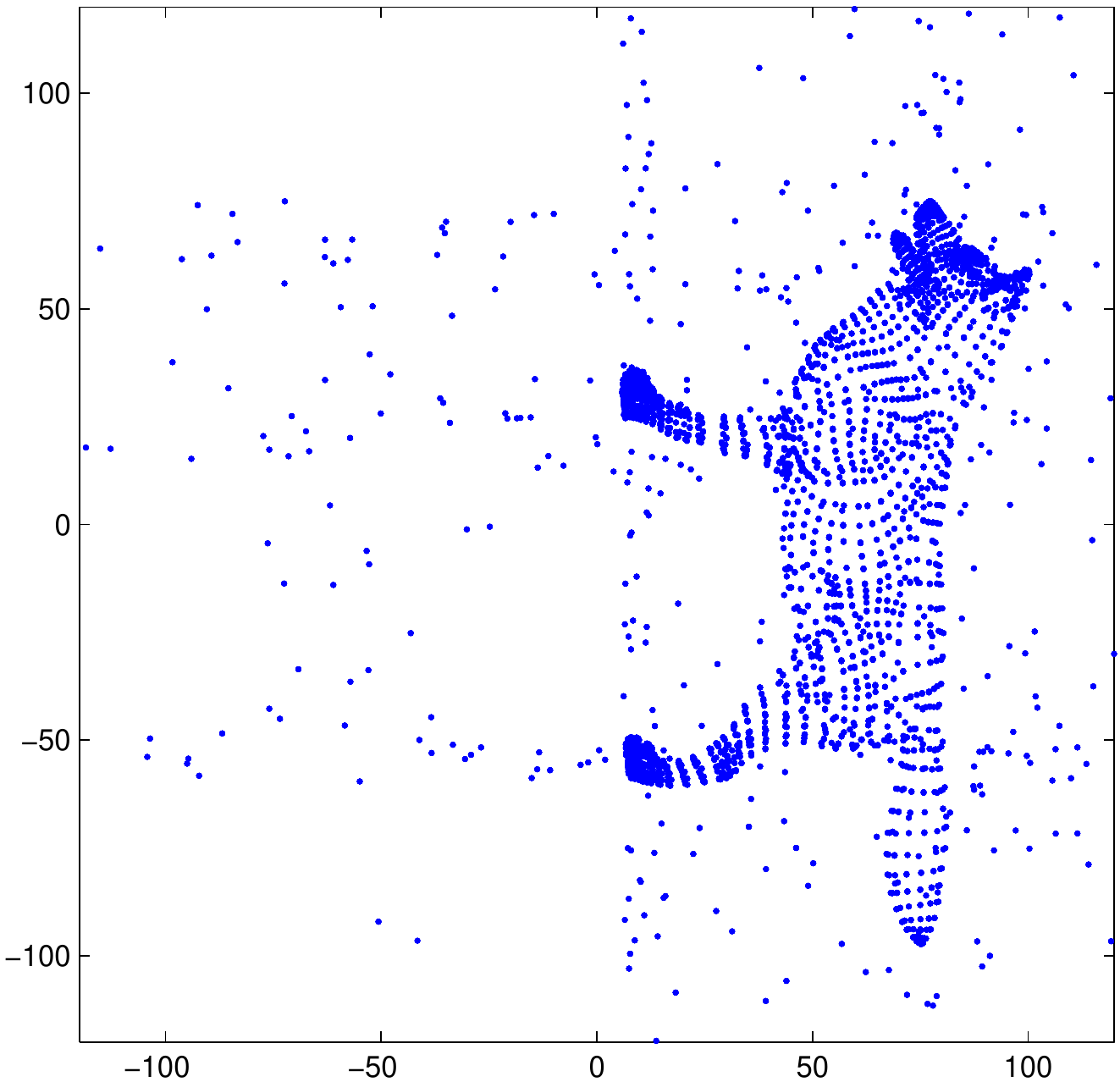}
&\includegraphics[width=0.12\textwidth,
keepaspectratio]{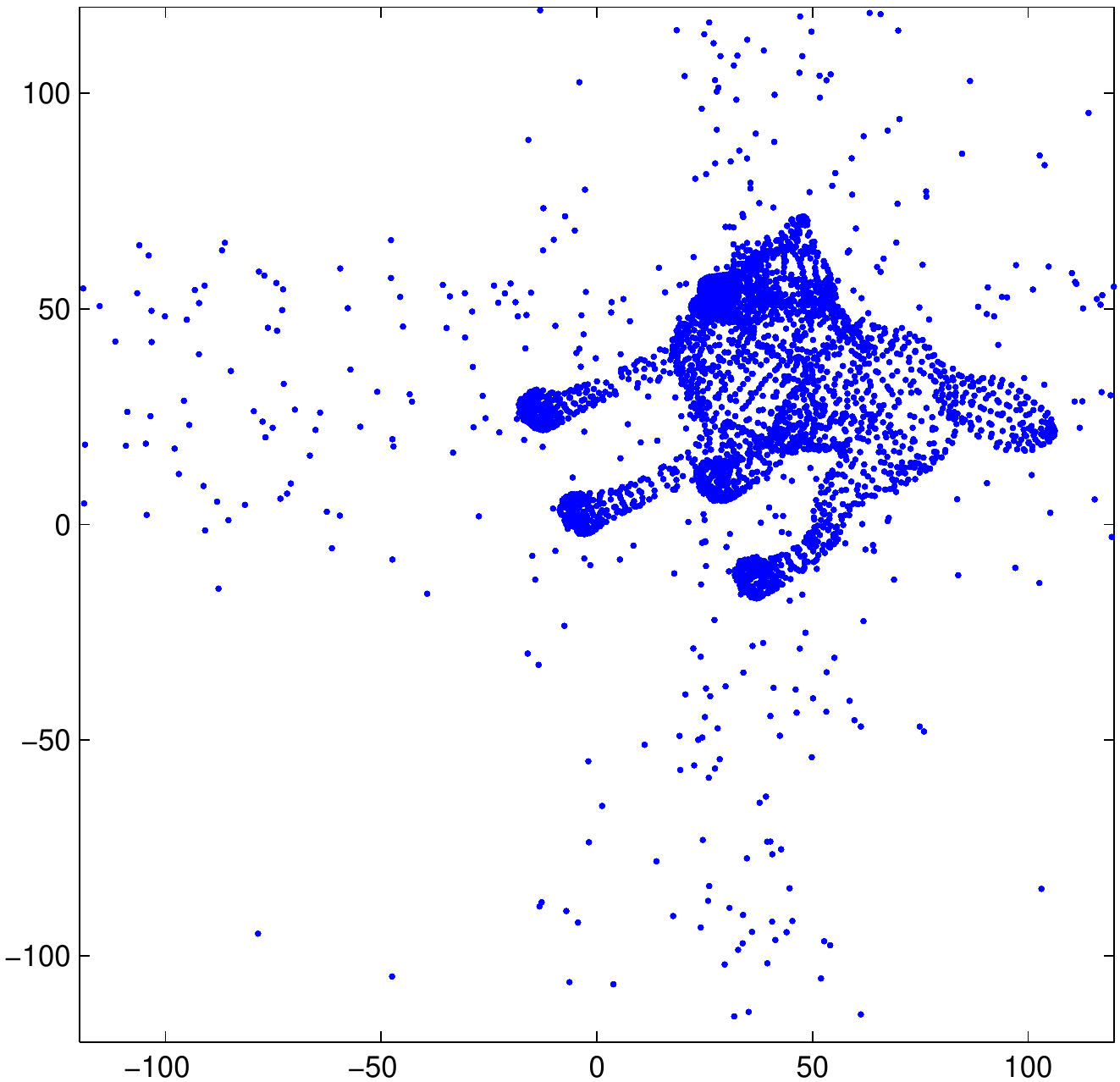}
&\includegraphics[width=0.12\textwidth,
keepaspectratio]{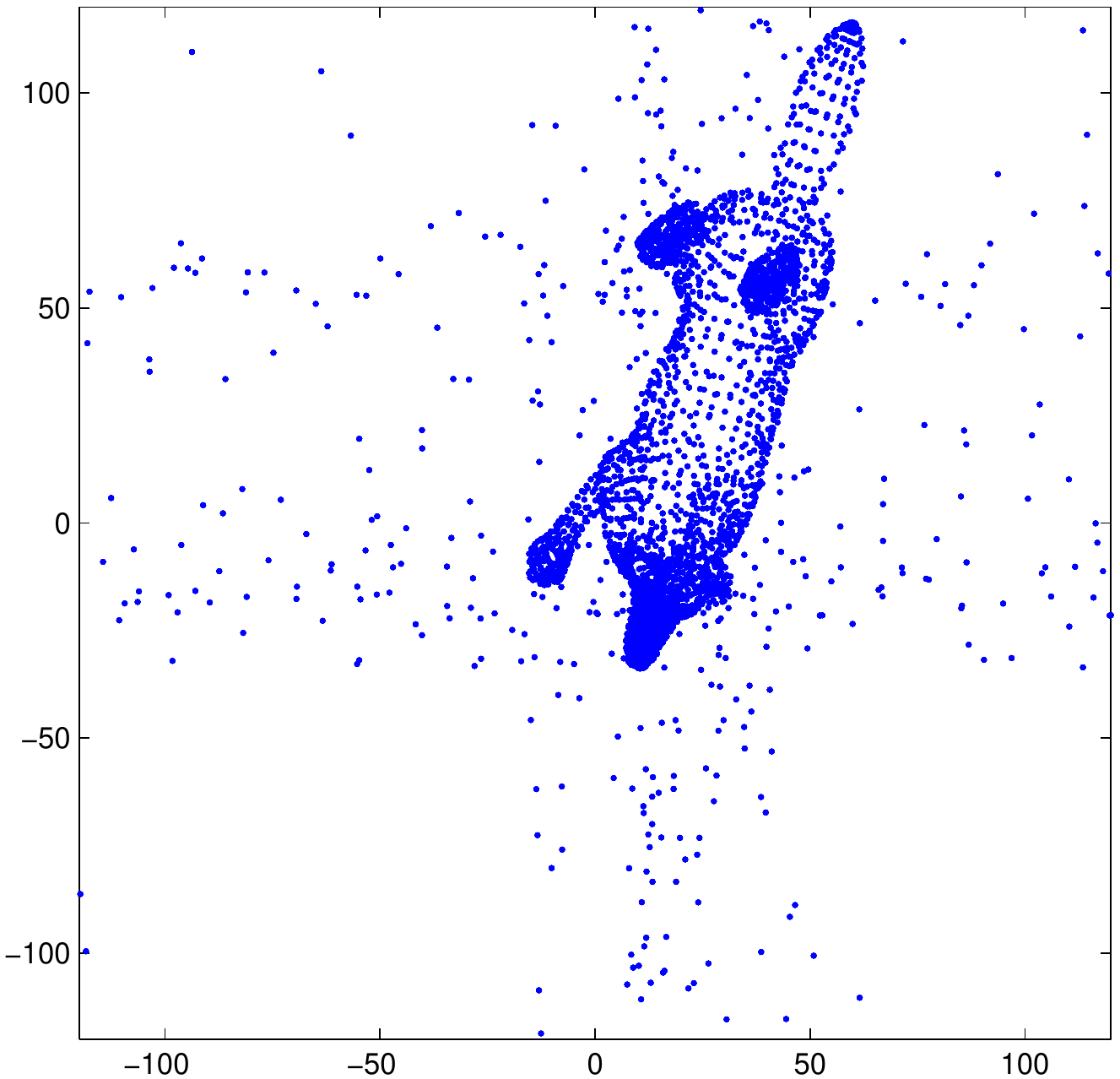}
&\includegraphics[width=0.12\textwidth,
keepaspectratio]{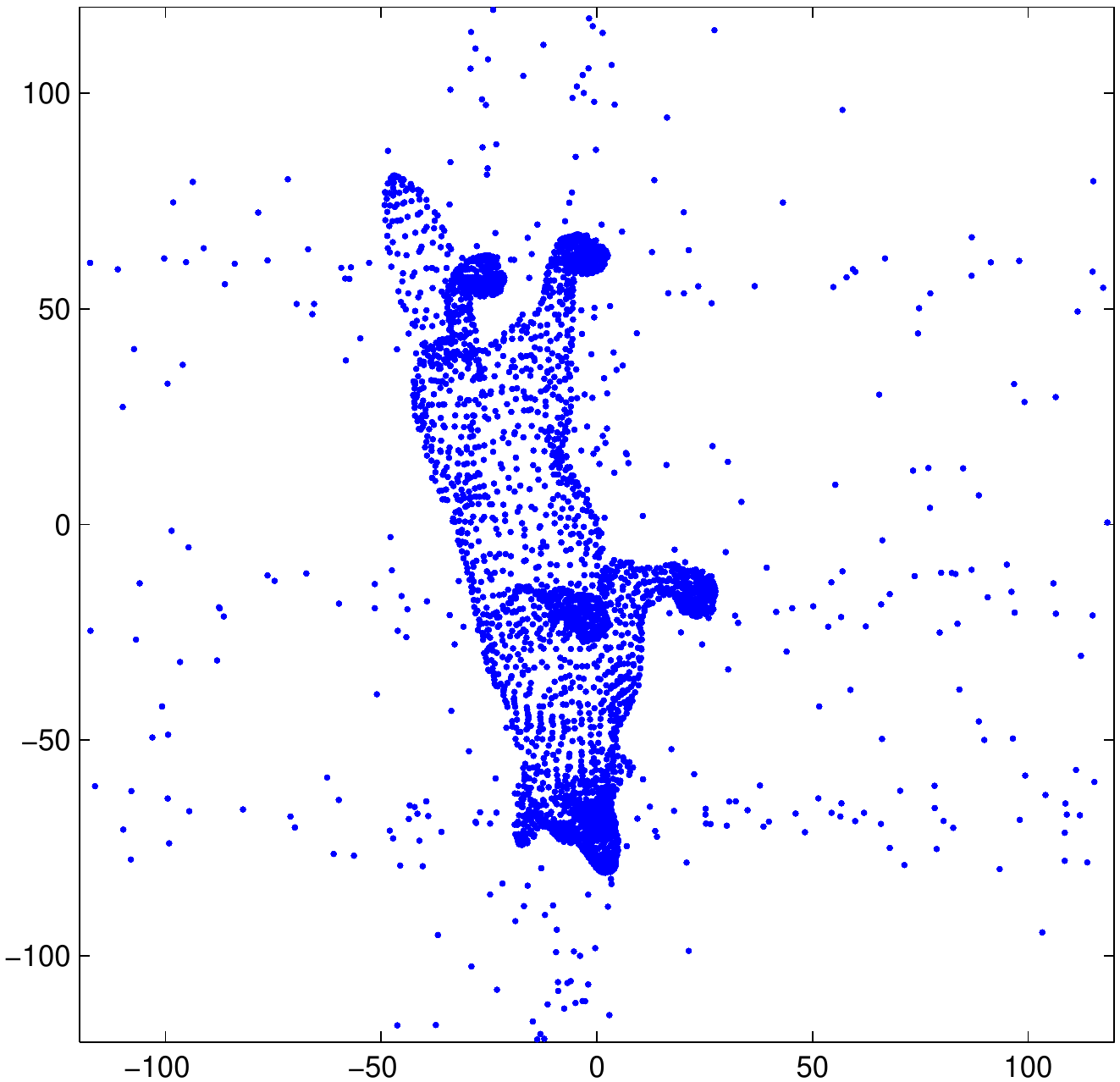}
&\includegraphics[width=0.12\textwidth,
keepaspectratio]{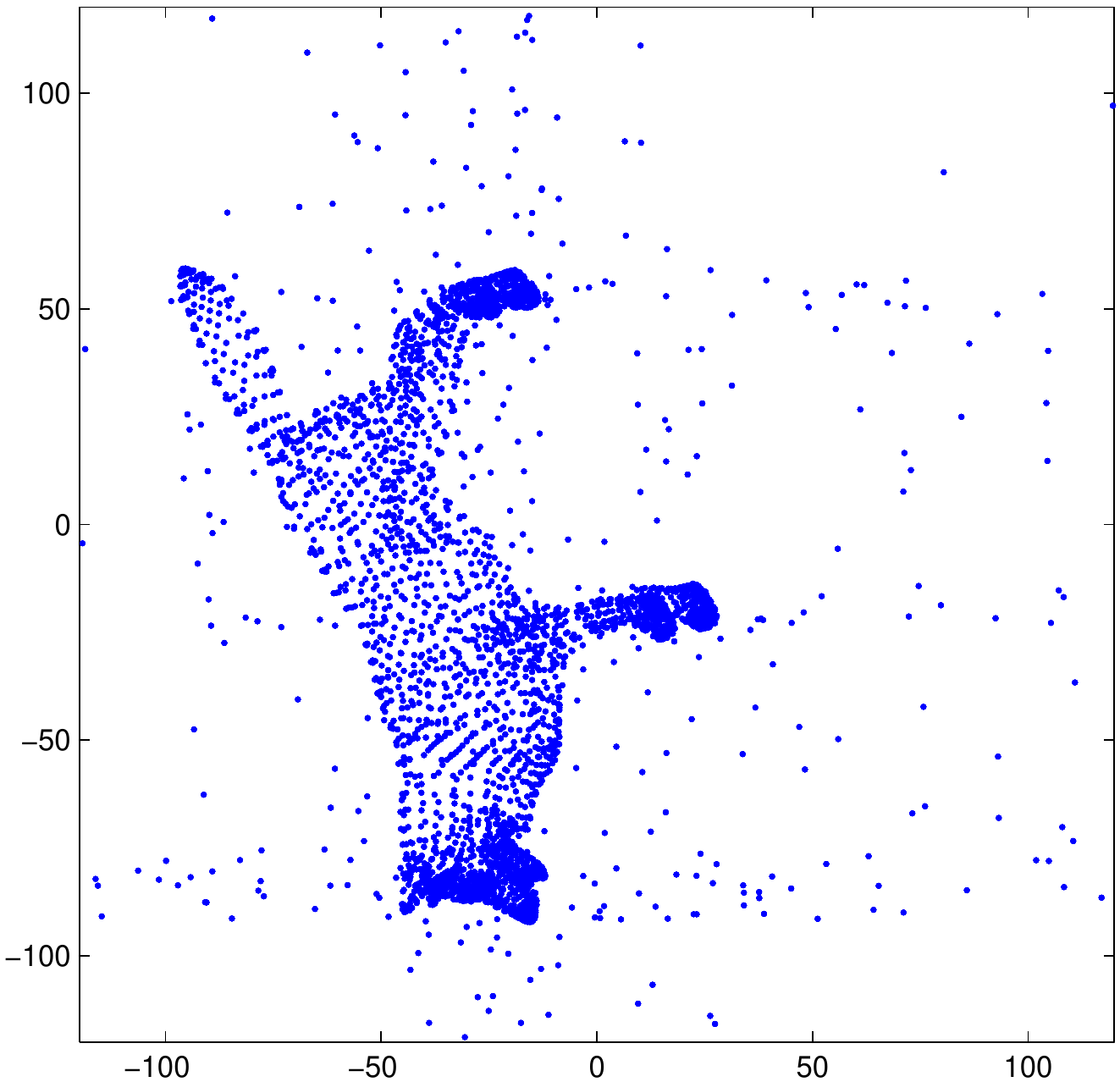}
&\includegraphics[width=0.12\textwidth,
keepaspectratio]{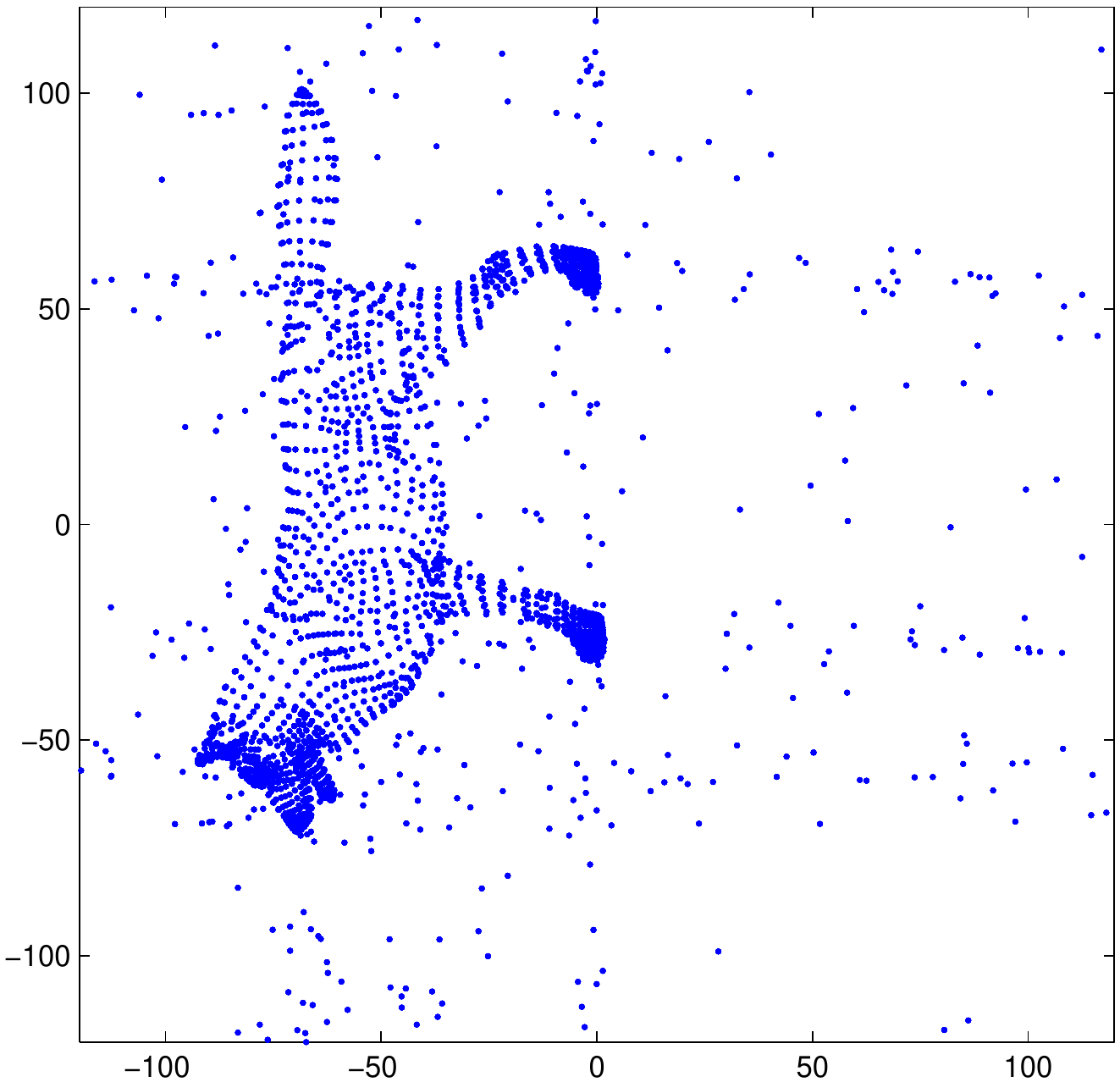}\\
\end{tabular}
\caption{The illustrations of some trajectories (2D image frames)
generated by the 3D ``wolf'' object (300th, 500th, $\cdots$,
1300th, 1500th frames). Top row: the ground truth trajectories.
Bottom row: $10\%$ corrupted trajectories.}\label{fig:sfm_im}
\end{figure*}

For this experiment, we first generate the 2D feature points
$\mathbf{L}_0$ by applying an affine camera model (with rotation
angles between 0 and $2\pi$, with a step size $\pi/1000$, and
uniformly randomly generated translations) to the 3D ``wolf''
object~\footnote{The 3D ``wolf'' data is available at:
http://tosca.cs.technion.ac.il/.}, which contains $4344$ 3D
points. Then we add impulsive noises $\mathbf{S}_0$ (the locations
of the nonzero entries are uniformly sampled at random) to part
(e.g., $5\%$ or $10\%$) of the feature points (see
Figure~\ref{fig:sfm_im}). In this way, we obtain corrupted
observations $\mathbf{M}=\mathbf{L}_0 + \mathbf{S}_0$ with a size
$4002 \times 4344$.

\begin{figure*}[ht]
\centering
\begin{tabular}{c@{\extracolsep{1em}}c@{\extracolsep{1em}}c@{\extracolsep{1em}}c@{\extracolsep{1em}}c}
\includegraphics[width=0.18\textwidth,
keepaspectratio]{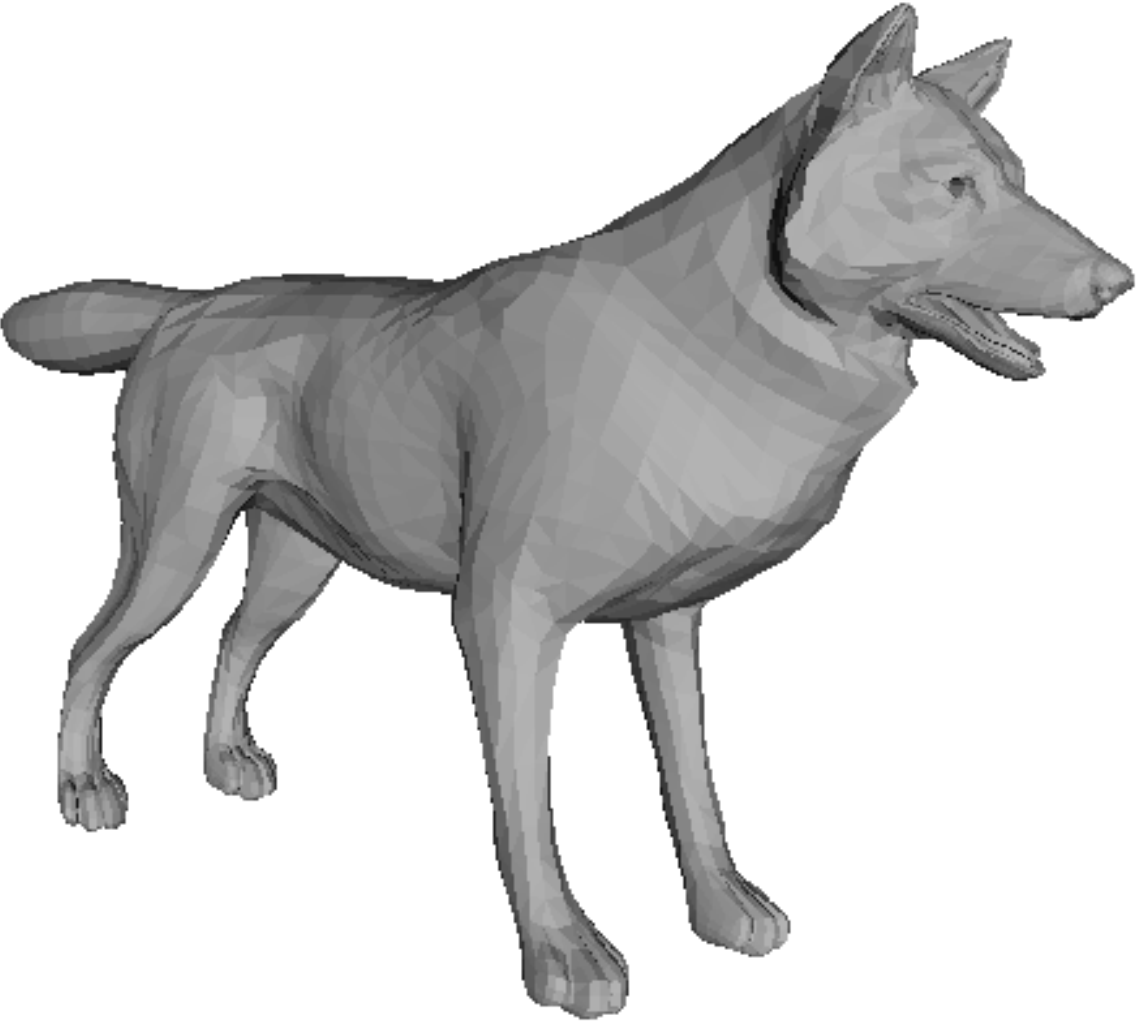}
&\includegraphics[width=0.18\textwidth,
keepaspectratio]{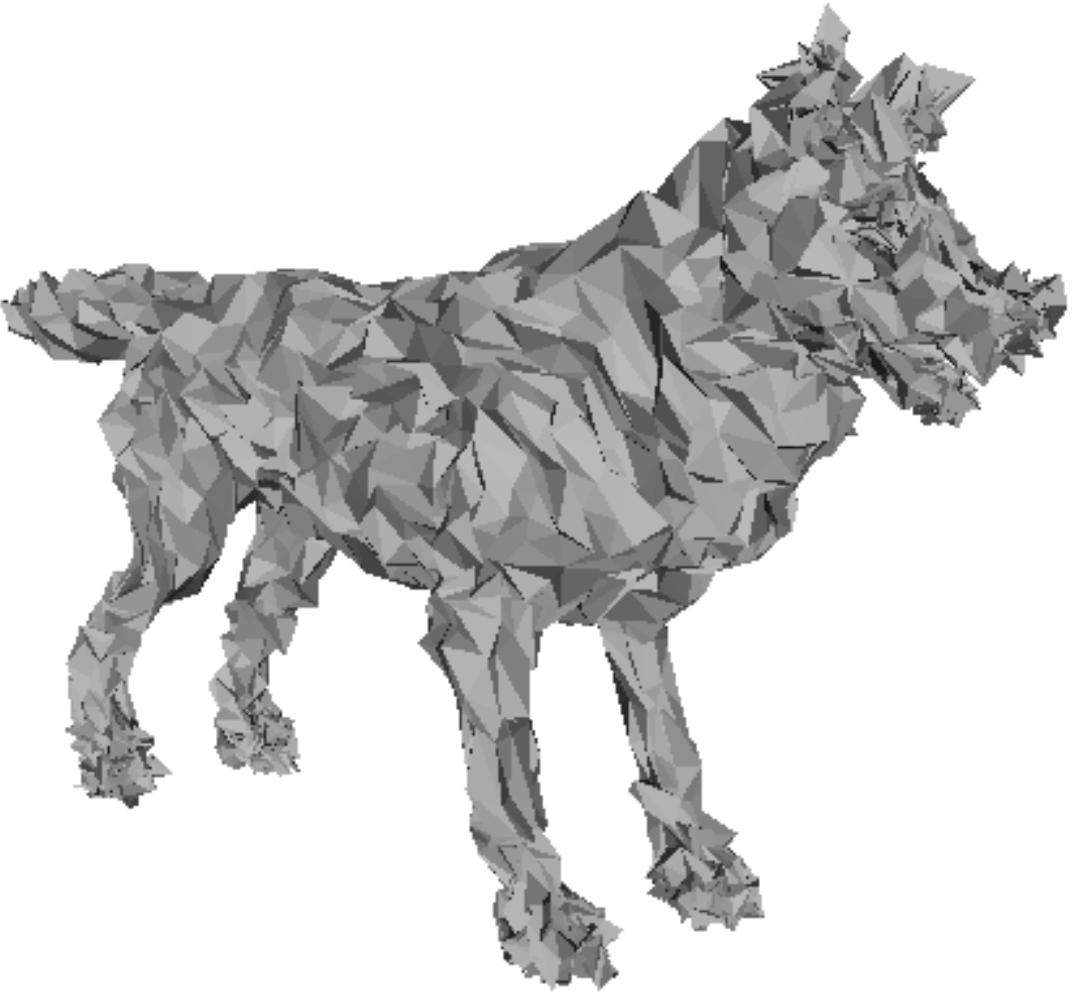}
&\includegraphics[width=0.18\textwidth,
keepaspectratio]{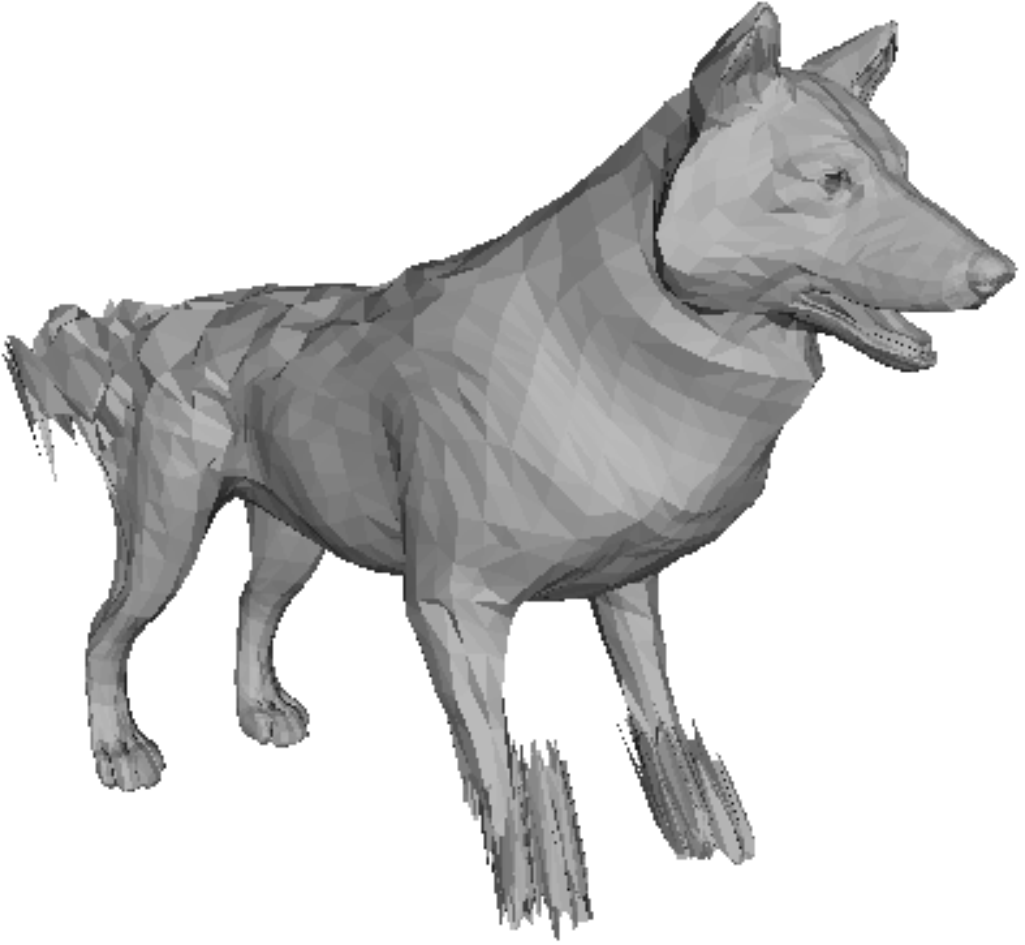}
&\includegraphics[width=0.18\textwidth,
keepaspectratio]{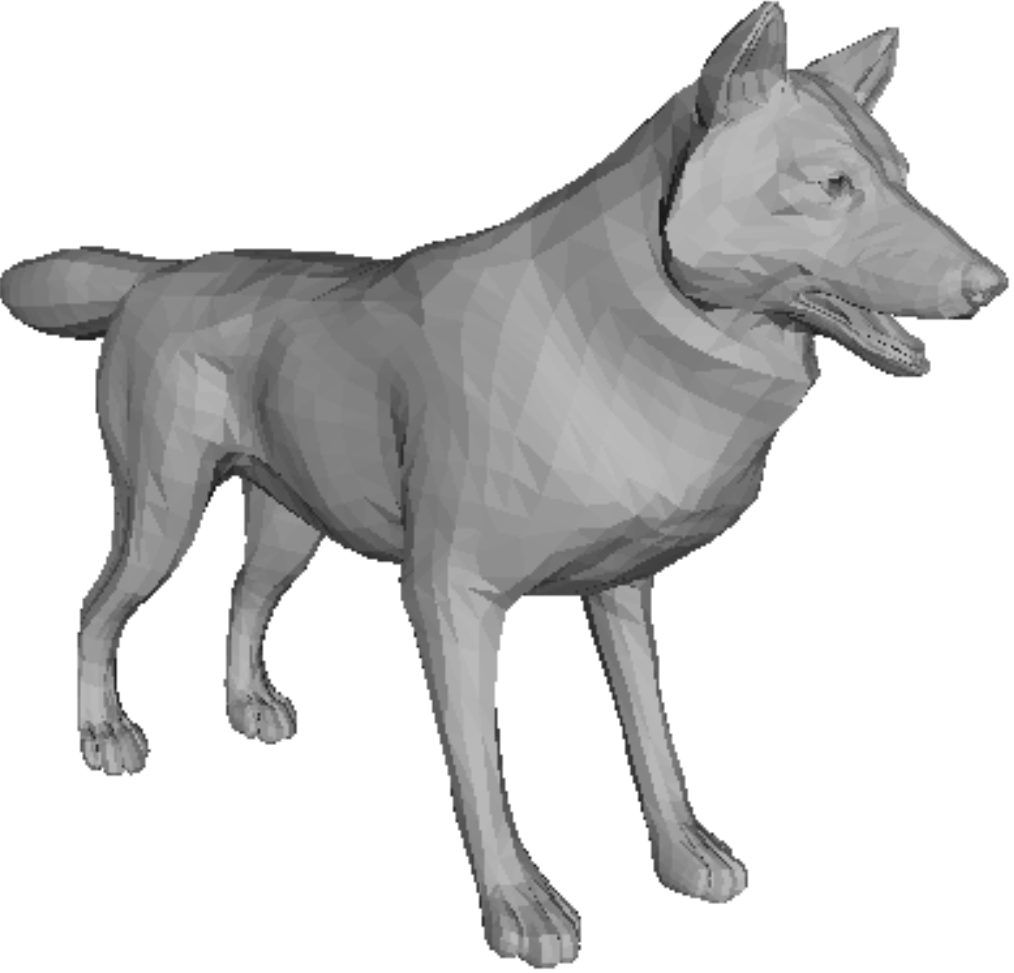}
&\includegraphics[width=0.18\textwidth,
keepaspectratio]{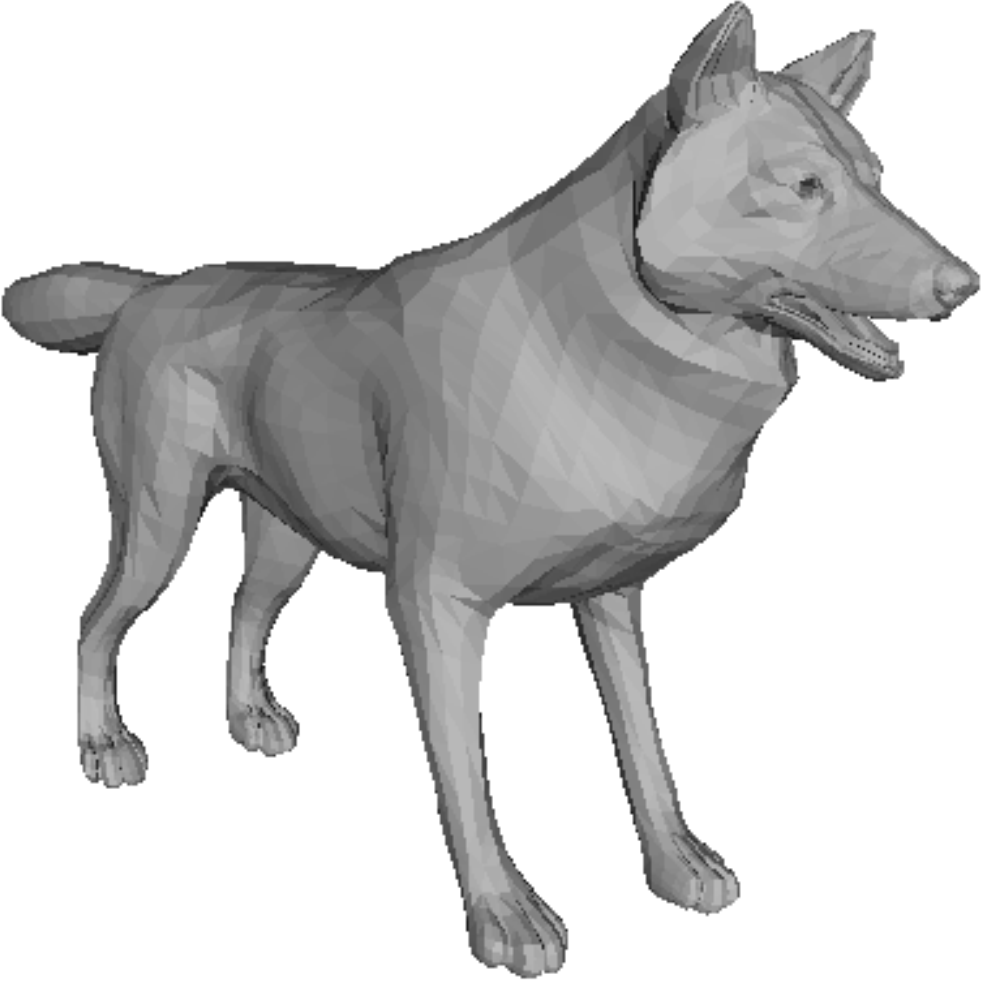}\\
(a) Original & (b) Corrupted & (c) RSL & (d) S-PCP & (e)
$l_1$-filtering
\end{tabular}
\caption{The SfM reconstruction. (a) is the original 3D object.
(b)-(e) are SfM results using corrupted trajectory and the
trajectories recovered by RSL, standard PCP (S-PCP for short) and
$l_1$ filtering, respectively.}\label{fig:sfm}
\end{figure*}

We apply our $l_1$ filtering to remove outliers (i.e.,
$\mathbf{S}_0$) and compute the affine motion matrix $\mathbf{A}$
and the 3D coordinates $\mathbf{B}$ from the recovered features
(i.e., $\mathbf{L}_0$). For comparison, we also include the
results from the robust subspace learning (RSL)
(\cite{fernando2003rsl})\footnote{The Matlab code of RSL is
available at
http://www.salleurl.edu/$\sim$ftorre/papers/rpca/rpca.zip and the
parameters in this code are set to their default values.} and
standard PCP (i.e., S-ADM based PCP). In Figure~\ref{fig:sfm}, we
show the original 3D object, SfM results based on noisy
trajectories and trajectories recovered by RSL, standard PCP and
$l_1$ filtering, respectively. It is easy to see that the 3D
reconstruction of RSL fails near the front legs and tail. In
contrast, the standard PCP and $l_1$ filtering provide results
with almost the same quality. Table~\ref{tab:sfm} further compares
the numerical behaviors of these methods. We measure the
quantitative performance for SfM by the well-known mean 2D
reprojection error, which is denoted as ``ReprojErr'' and defined
by the mean distance of the ground truth 2D feature points and
their reprojections. We can see that the standard PCP provides the
highest numerical accuracy while its time cost is extremely high
(9 times slower than RSL and more than 100 times slower than $l_1$
filtering). Although the speed of RSL is faster than standard PCP,
its numerical accuracy is the worst among these methods. In
comparison, our $l_1$ filtering achieves almost the same numerical
accuracy as standard PCP and is the fastest.

\begin{table*}[th]
\begin{center}
\caption{Comparison among RSL,
S-PCP and $l_1$ filtering on the structure from motion problem.
We present CPU time (in seconds) and the numerical
accuracy of tested algorithms. $\mathbf{L}_0$ and $\mathbf{S}_0$
are the ground truth and $\mathbf{L}^*$ and $\mathbf{S}^*$ are the
solution computed by different methods.}\label{tab:sfm}
\begin{tabular}{|c||c|c|c|c|c|c|c|c|}\hline
Noisy Level &Method&  RelErr     & $\rank(\mathbf{L}^*)$    &
$\|\mathbf{S}^*\|_{l_0}$   &   Time & MaxDif($\mathbf{L}^*$) &
AveDif($\mathbf{L}^*$) & ReprojErr\\\hline\hline
\multirow{4}{*}{$5\%$} &
\multicolumn{8}{c|}{$\rank(\mathbf{L}_0)=4$, \qquad
$\|\mathbf{S}_0\|_{l_0}=869234$} \\\cline{2-9} &RSL   & 0.0323 & 4
(fixed)   & 15384229& 93.05  & 32.1731                & 0.4777 &
0.9851\\\cline{2-9} &S-PCP & $5.18 \times 10^{-9}$ & 4 & 869200 &
848.09 & $1.70 \times 10^{-5}$ & $2.47\times 10^{-8}$ & $4.18
\times 10^{-8}$\\\cline{2-9} &$l_1$ & $1.16 \times 10^{-8}$ & 4
& 869644 &  \textbf{6.46}   & $1.80 \times 10^{-5}$ & $3.61 \times
10^{-7}$ & $4.73 \times 10^{-7}$
\\\hline\hline \multirow{4}{*}{$10\%$} &
\multicolumn{8}{c|}{$\rank(\mathbf{L}_0)=4$, \qquad
$\|\mathbf{S}_0\|_{l_0}=1738469$} \\\cline{2-9} &RSL   & 0.0550 &
4 (fixed)   & 16383294&  106.65 & 38.1621               & 0.9285 &
1.8979\\\cline{2-9} &S-PCP & $6.30 \times 10^{-9}$ & 4 & 1738410 &
991.40 & $1.57 \times 10^{-5}$ & $4.09 \times 10^{-8}$ & $6.82
\times 10^{-7}$\\\cline{2-9} &$l_1$ & $3.18 \times 10^{-8}$ & 4
& 1739912 &  \textbf{6.48}   & $5.61 \times 10^{-5}$ & $9.03
\times 10^{-7}$ & $1.26 \times 10^{-6}$
\\\hline
\end{tabular}
\end{center}
\end{table*}

\subsection{Background Modeling}

In this subsection, we consider the problem of background modeling
from video surveillance. The background of a group of video
surveillance frames are supposed to be exactly the same and the
foreground on each frame is recognized as sparse errors. Thus this
vision problem can be naturally formulated as recovering the
low-rank matrix from its sum with sparse errors
(\cite{candes2009robust}). We compare our $l_1$ filtering with the
baseline median filter\footnote{Please refer to
http://en.wikipedia.org/wiki/Median\_filter.} and other
state-of-the-art robust approaches, such as RSL and S-PCP. When
median filtering using all the frames, the complexity is actually
also quadratic and the results are not good. So we \emph{only
buffer 20 frames when using median filter} to compute the
background. For $l_1$ filtering, we set the size of the seed
matrix as $20\times 20$.

For quantitative evaluation, we perform all the compared methods
on the ``laboratory'' sequence from a public surveillance
database~(\cite{Benedek2008VS}) which has ground truth foreground.
Both the false negative rate (FNR) and the false positive rate
(FPR) are calculated in the sense of foreground detection. FNR
indicates the ability of the method to correctly recover the
foreground while the FPR represents the power of a method on
distinguishing the background. These two scores correspond to the
Type I and Type II errors in the statistical test
theory\footnote{Please refer to
http://en.wikipedia.org/wiki/Type\_I\_and\_type\_II\_errors.} and
are judged by the criterion that the smaller the better. One can
see from Table~\ref{tab:vs} that RSL has the lowest FNR but the
highest FPR among the compared methods. This reveals that RSL
could not exactly distinguish the background. Although the speed
of our $l_1$ filtering is slightly slower than median filtering on
20 frames, its performance is as good as S-PCP, which achieves the
best results but with the highest time cost.

\begin{table}[ht]
\begin{center}
\caption{Comparison among median filter (Median for short), RSL, S-PCP, and $l_1$ filtering on background modeling problem.
``Resolution'' and ``No. Frames'' denote the size of each frame and the number of
frames in a video sequence, respectively. We present FNR, FPR and the CPU time (in seconds) for the ``laboratory'' data set.
For our collected ``meeting'' data set, we only report the CPU time because there is no ground truth foreground for this video sequence.}\label{tab:vs}
\begin{tabular}{|c|c|c|c|c|c|}
\hline
Video &  -    & Median &  RSL     & S-PCP    & $l_1$  \\\hline\hline
\multirow{4}{*}{``laboratory''} & \multicolumn{5}{|c|}{Resolution: $240\times 320$, \quad No. Frames: 887}\\\cline{2-6}
                          &FNR   & 9.85   &  7.31    & 8.61     & 8.62   \\\cline{2-6}
                          &FPR   & 9.18   &  10.83   & 8.72     & 8.76   \\\cline{2-6}
                          &Time  & 42.90  &  3159.92 & 10897.96 & \textbf{48.99}  \\\hline\hline
\multirow{2}{*}{``meeting''} & \multicolumn{5}{|c|}{Resolution: $576 \times 720$, \quad No. Frames: 700}\\\cline{2-6}
                          &Time & 179.19   & N.A.     & N.A.     & \textbf{178.74}\\\hline
\end{tabular}
\end{center}
\end{table}

\begin{figure*}[ht]
\centering
\begin{tabular}{c@{\extracolsep{0.2em}}c@{\extracolsep{0.2em}}c@{\extracolsep{0.2em}}c@{\extracolsep{0.2em}}c}
\includegraphics[width=0.18\textwidth,
keepaspectratio]{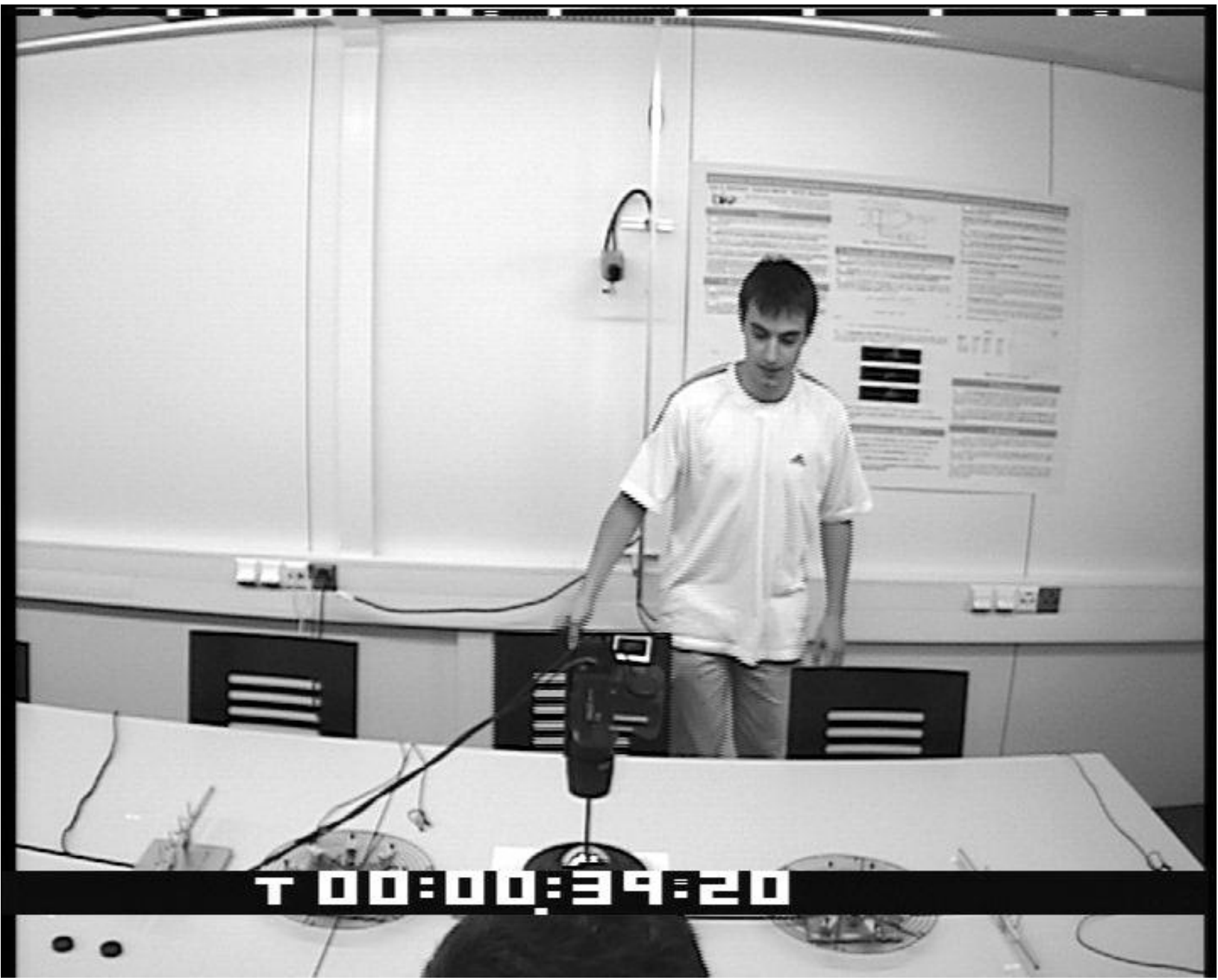}
&\includegraphics[width=0.18\textwidth,
keepaspectratio]{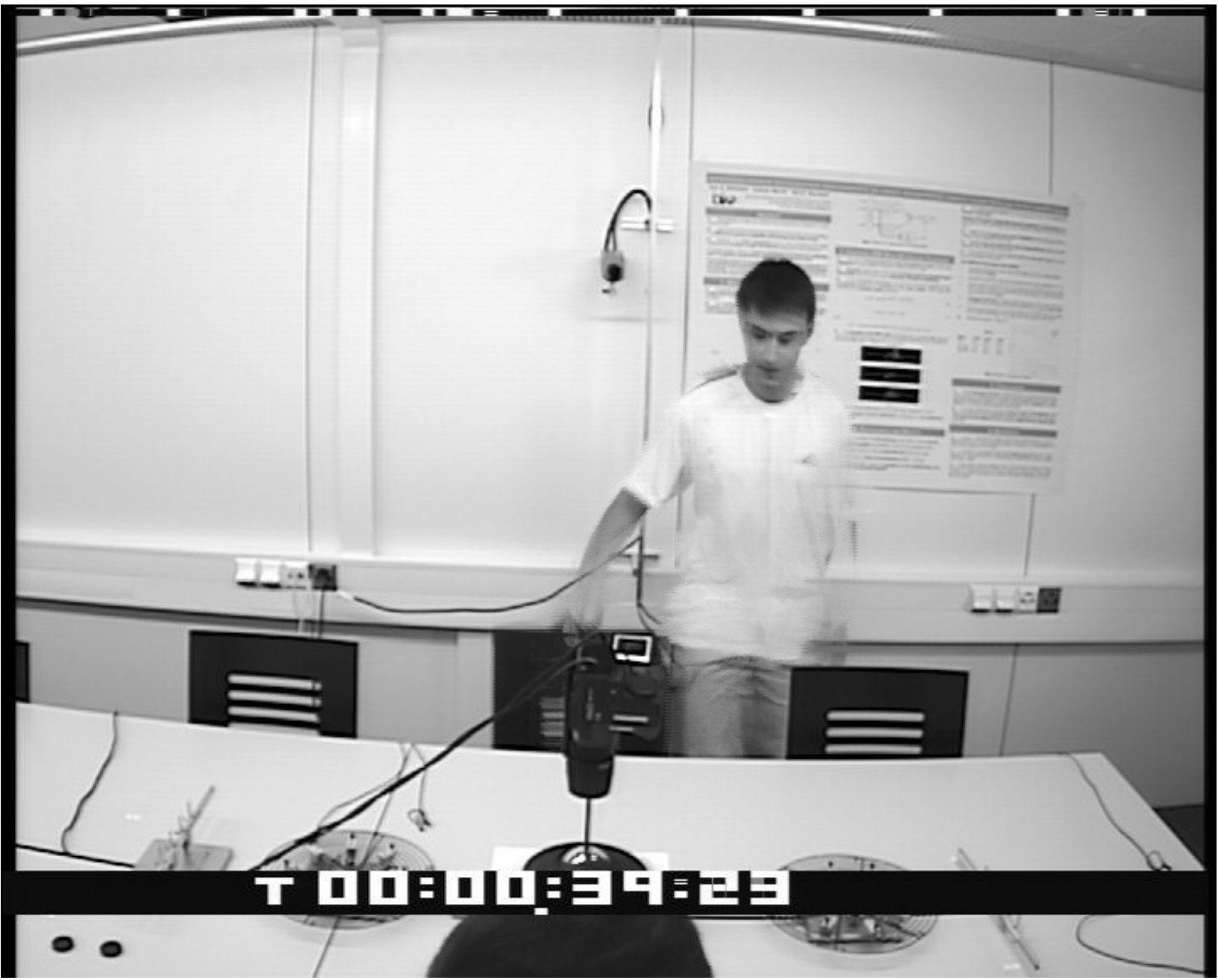}
&\includegraphics[width=0.18\textwidth,
keepaspectratio]{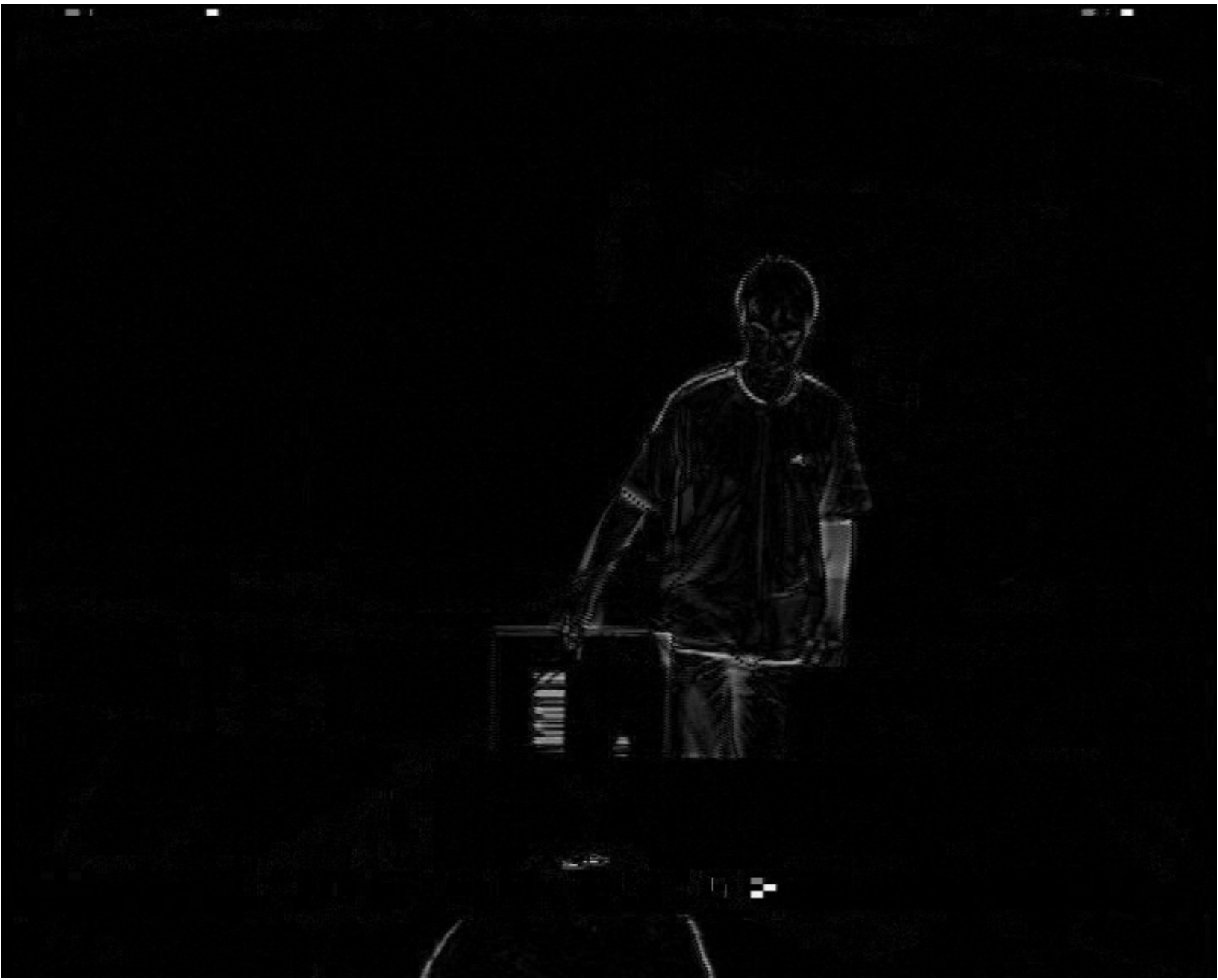}
&\includegraphics[width=0.18\textwidth,
keepaspectratio]{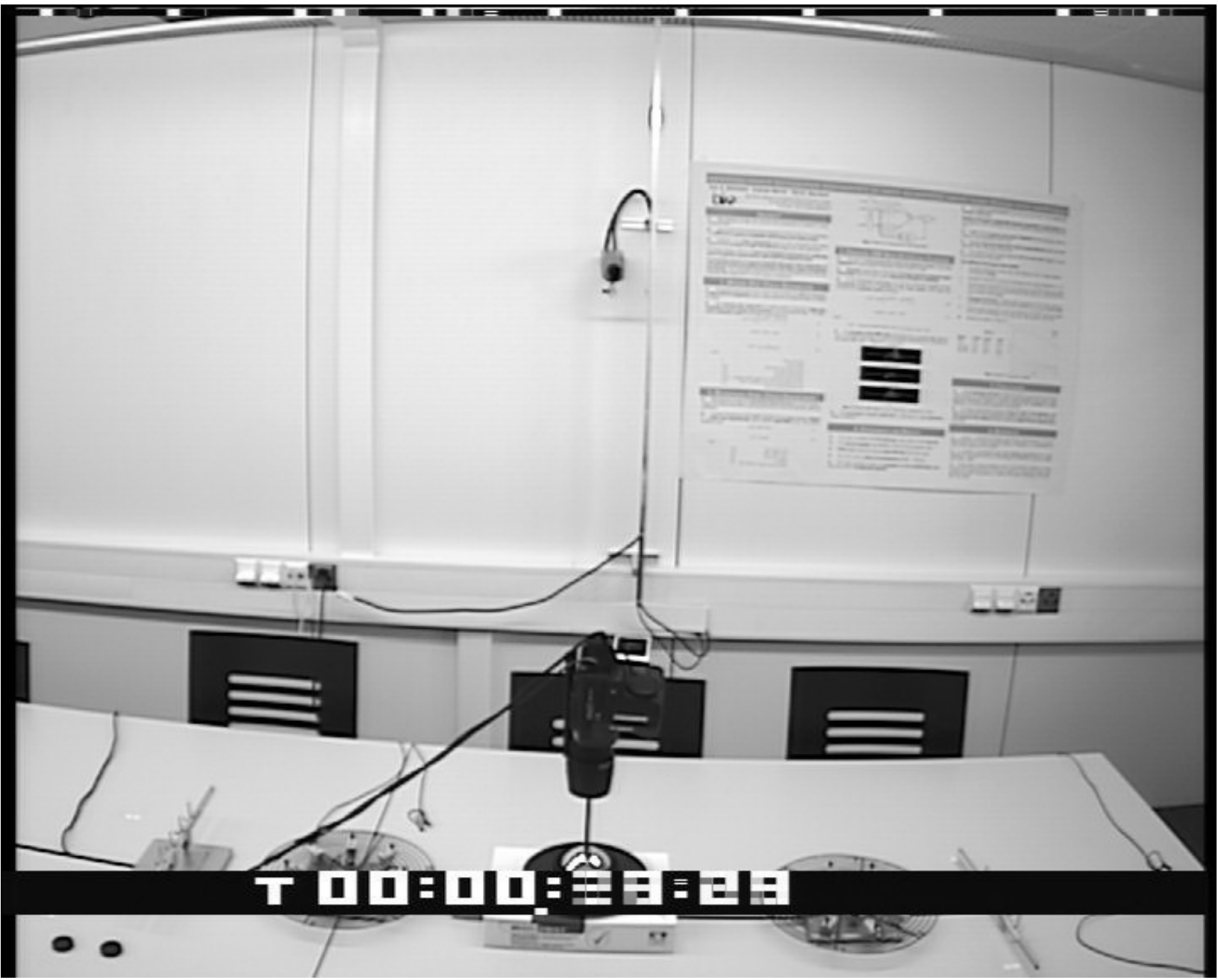}
&\includegraphics[width=0.18\textwidth,
keepaspectratio]{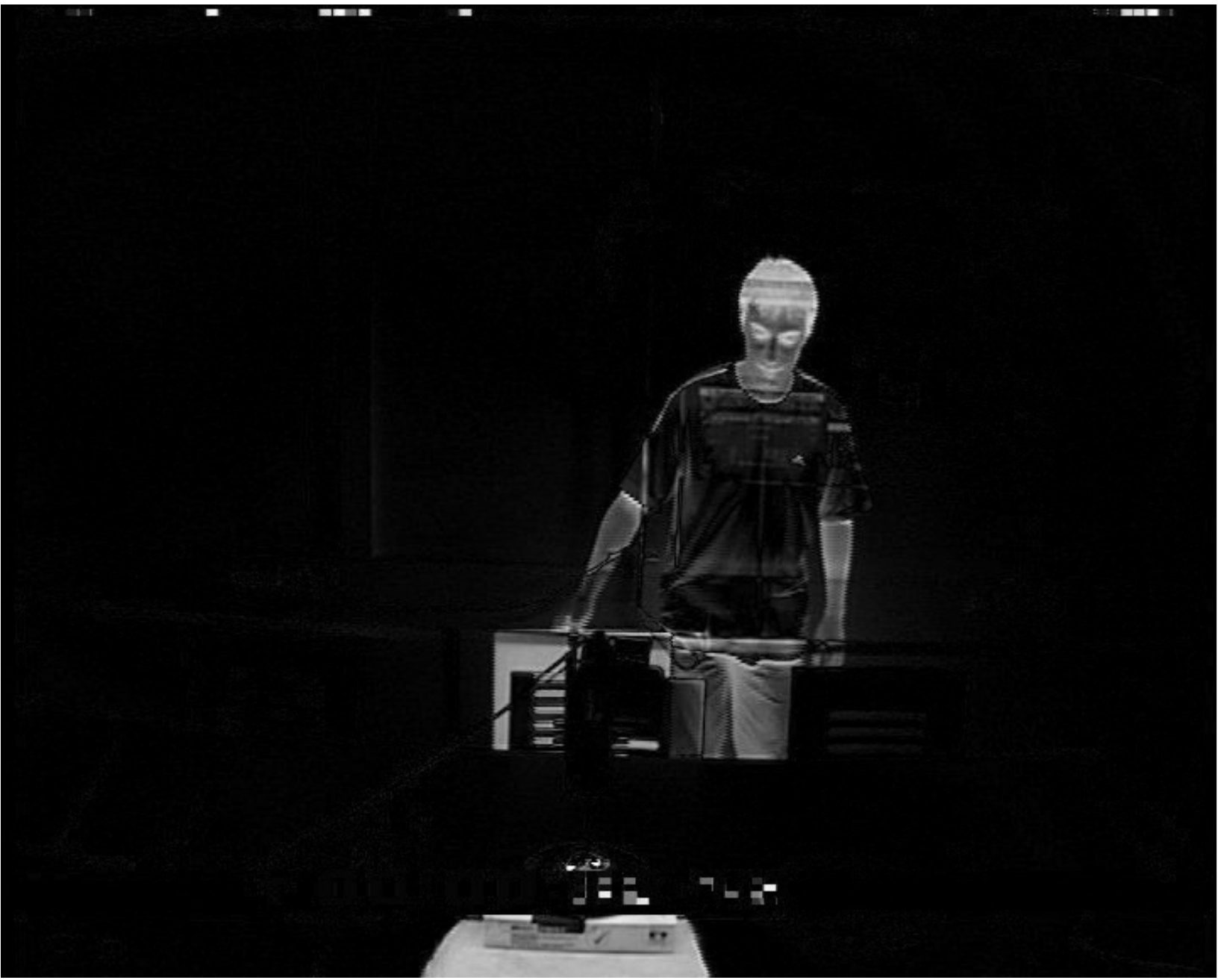}\\
\includegraphics[width=0.18\textwidth,
keepaspectratio]{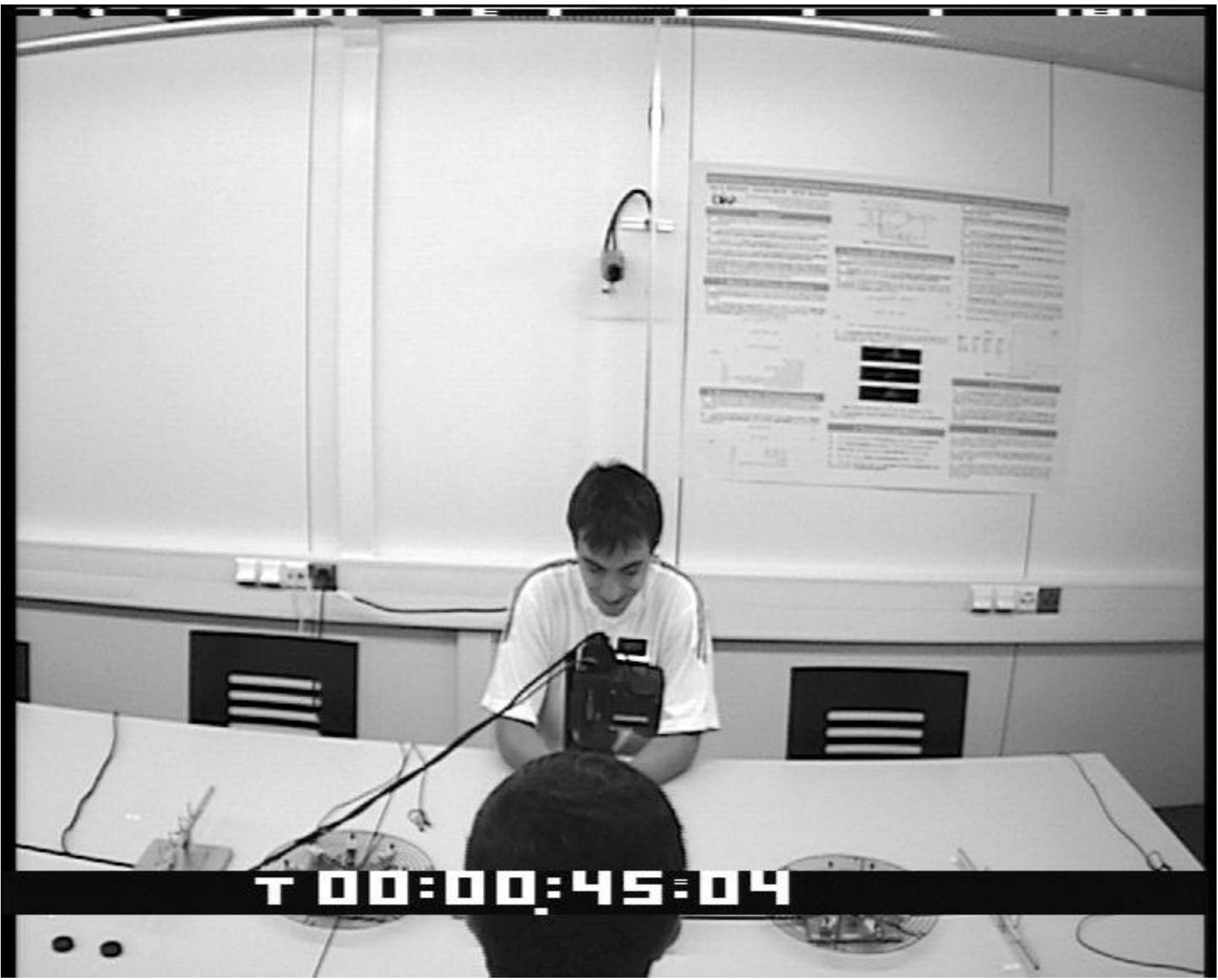}
&\includegraphics[width=0.18\textwidth,
keepaspectratio]{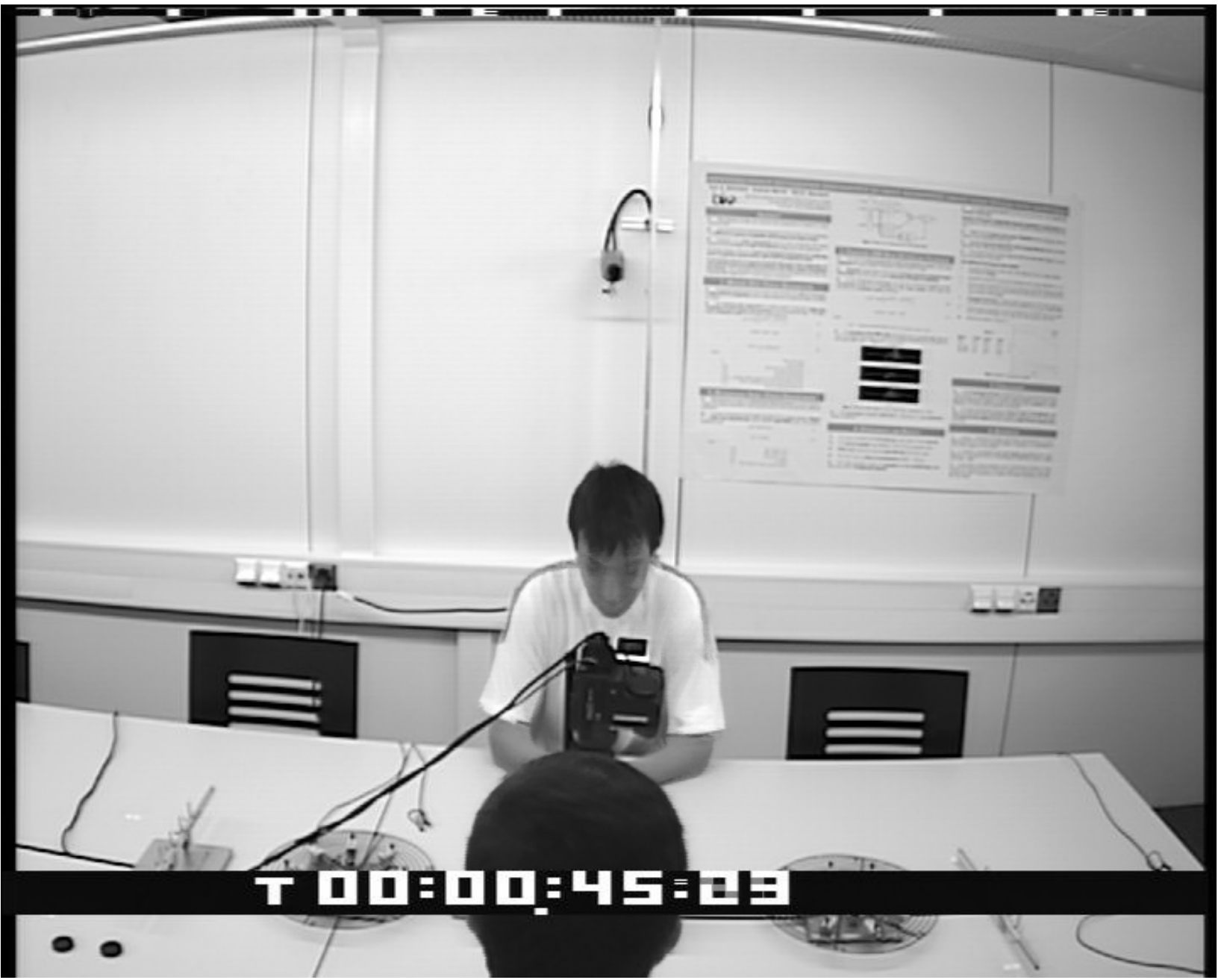}
&\includegraphics[width=0.18\textwidth,
keepaspectratio]{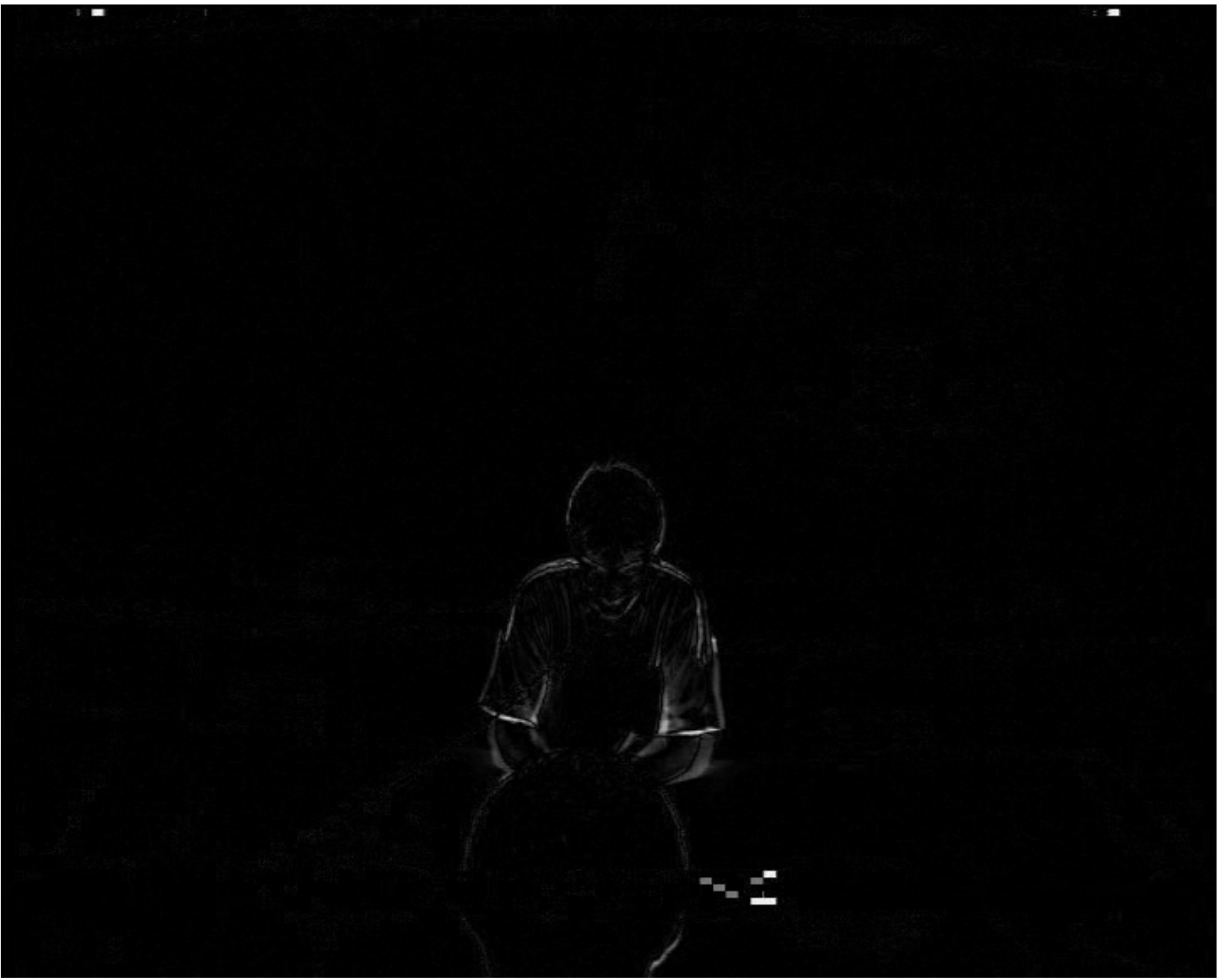}
&\includegraphics[width=0.18\textwidth,
keepaspectratio]{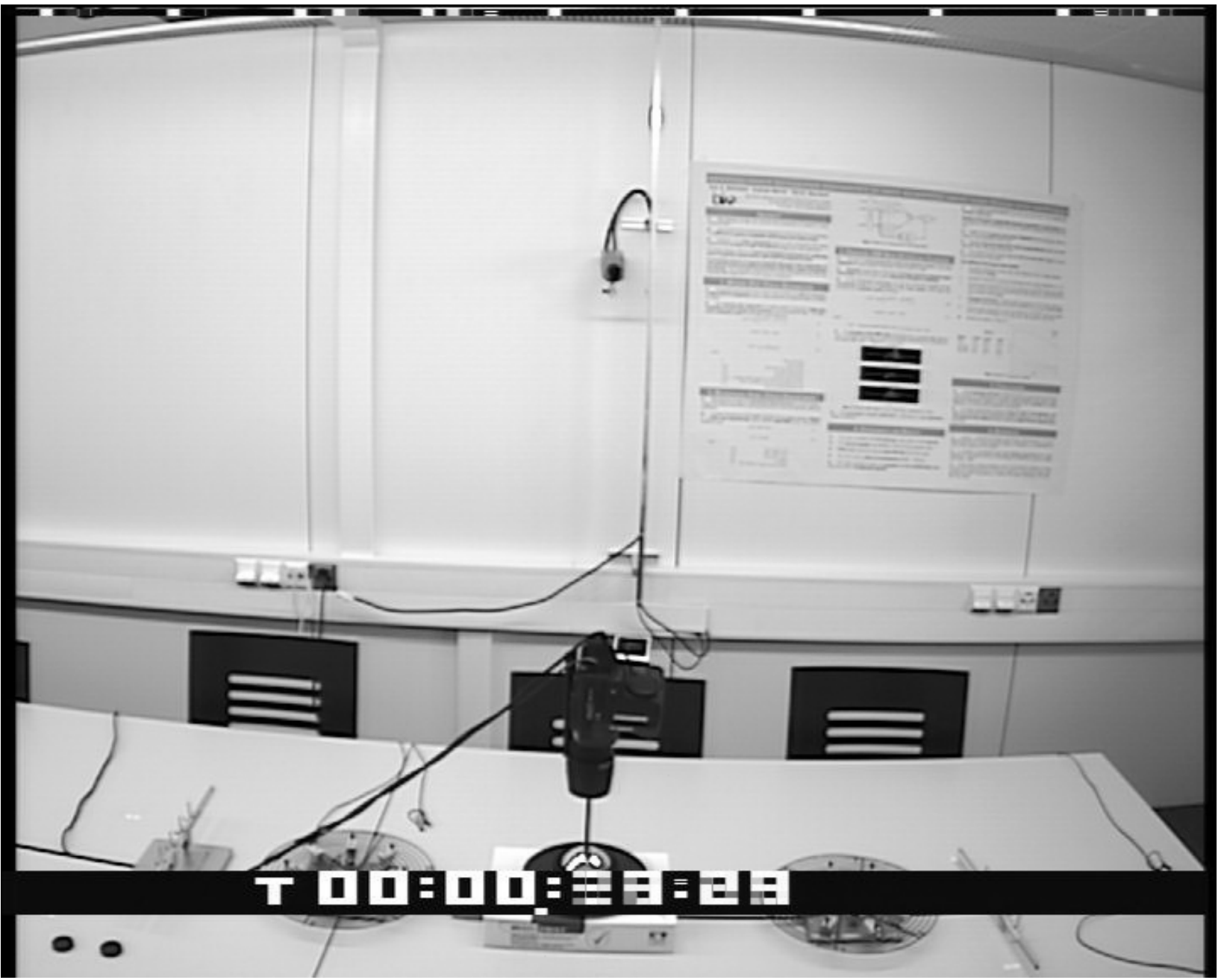}
&\includegraphics[width=0.18\textwidth,
keepaspectratio]{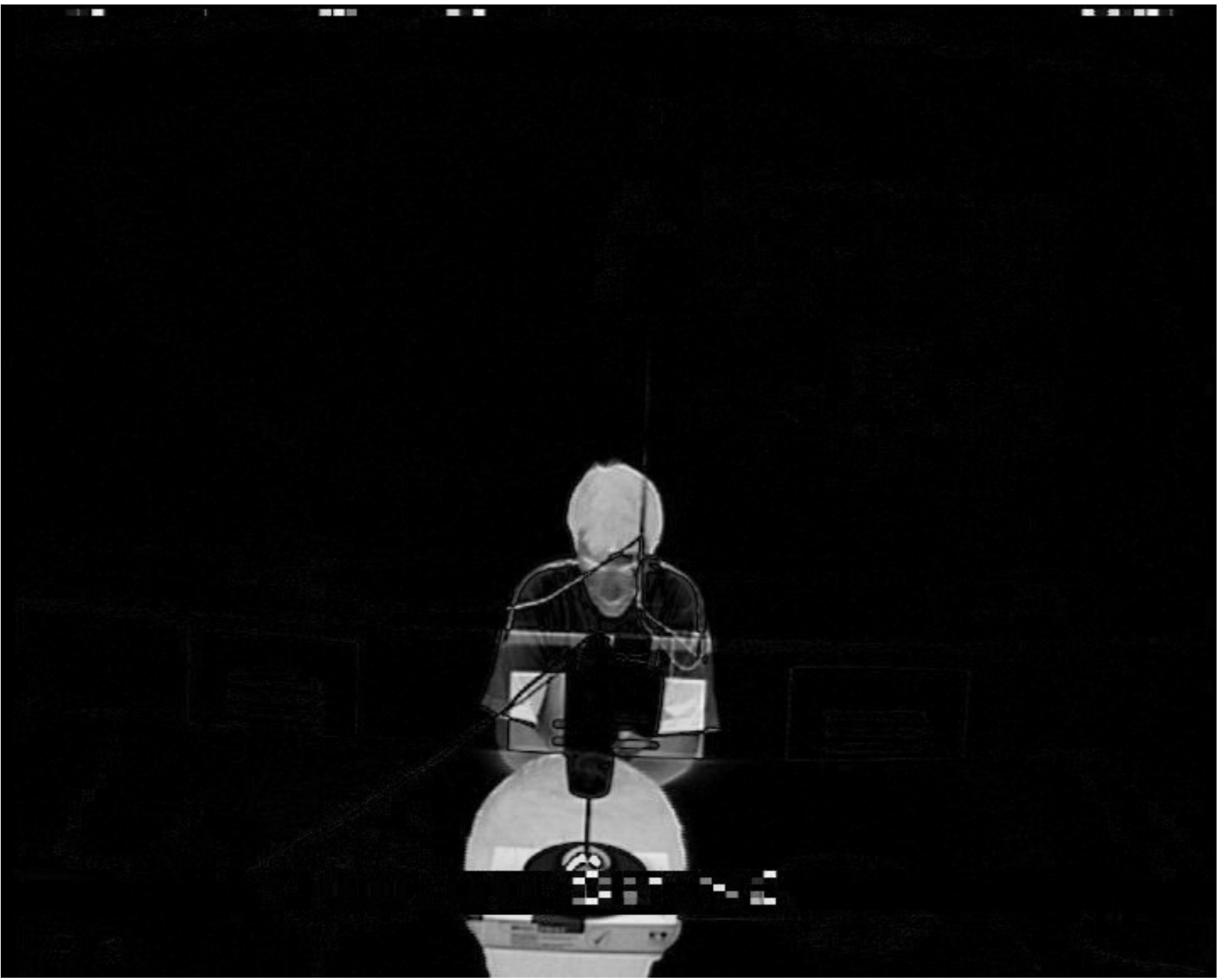}\\
\includegraphics[width=0.18\textwidth,
keepaspectratio]{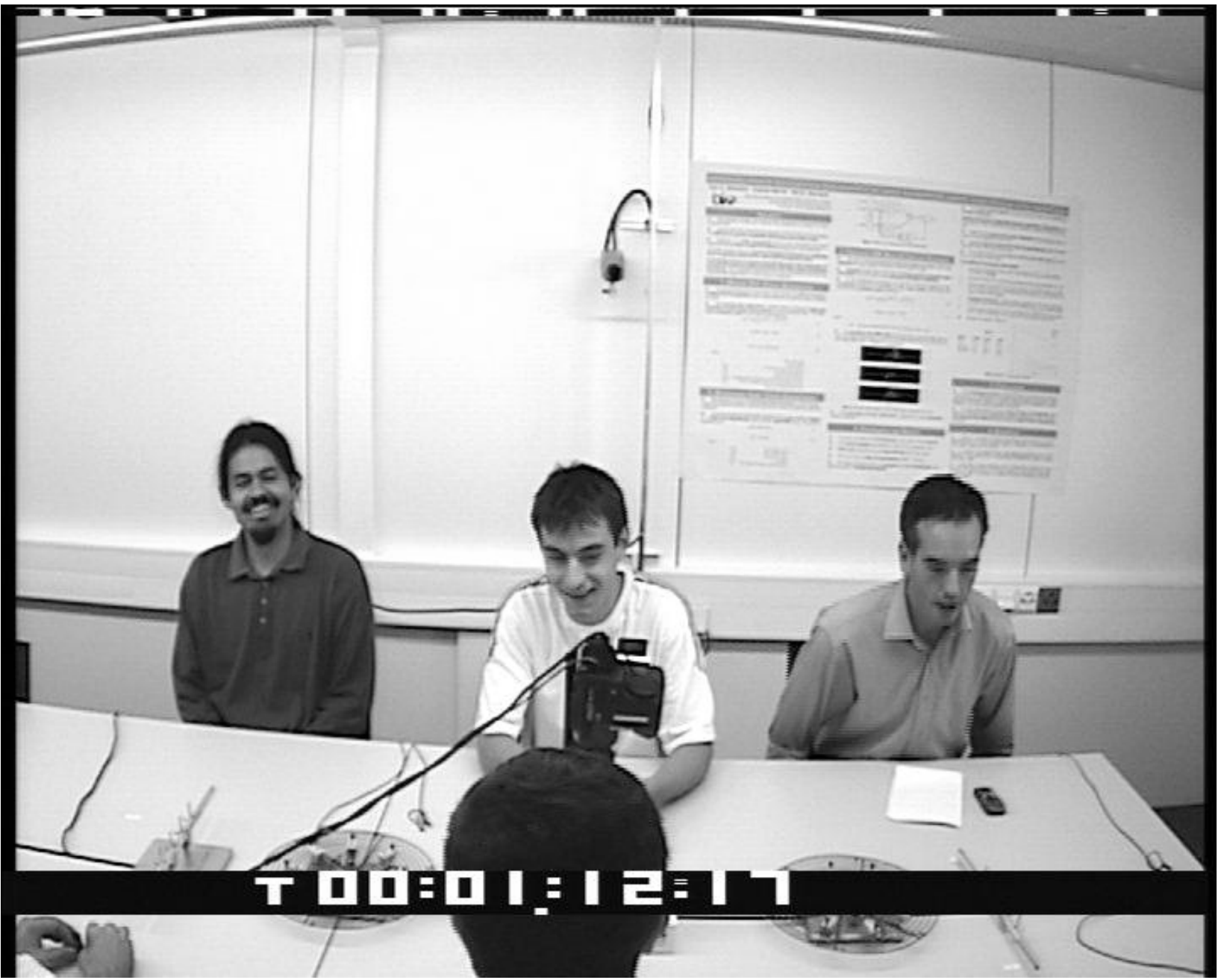}
&\includegraphics[width=0.18\textwidth,
keepaspectratio]{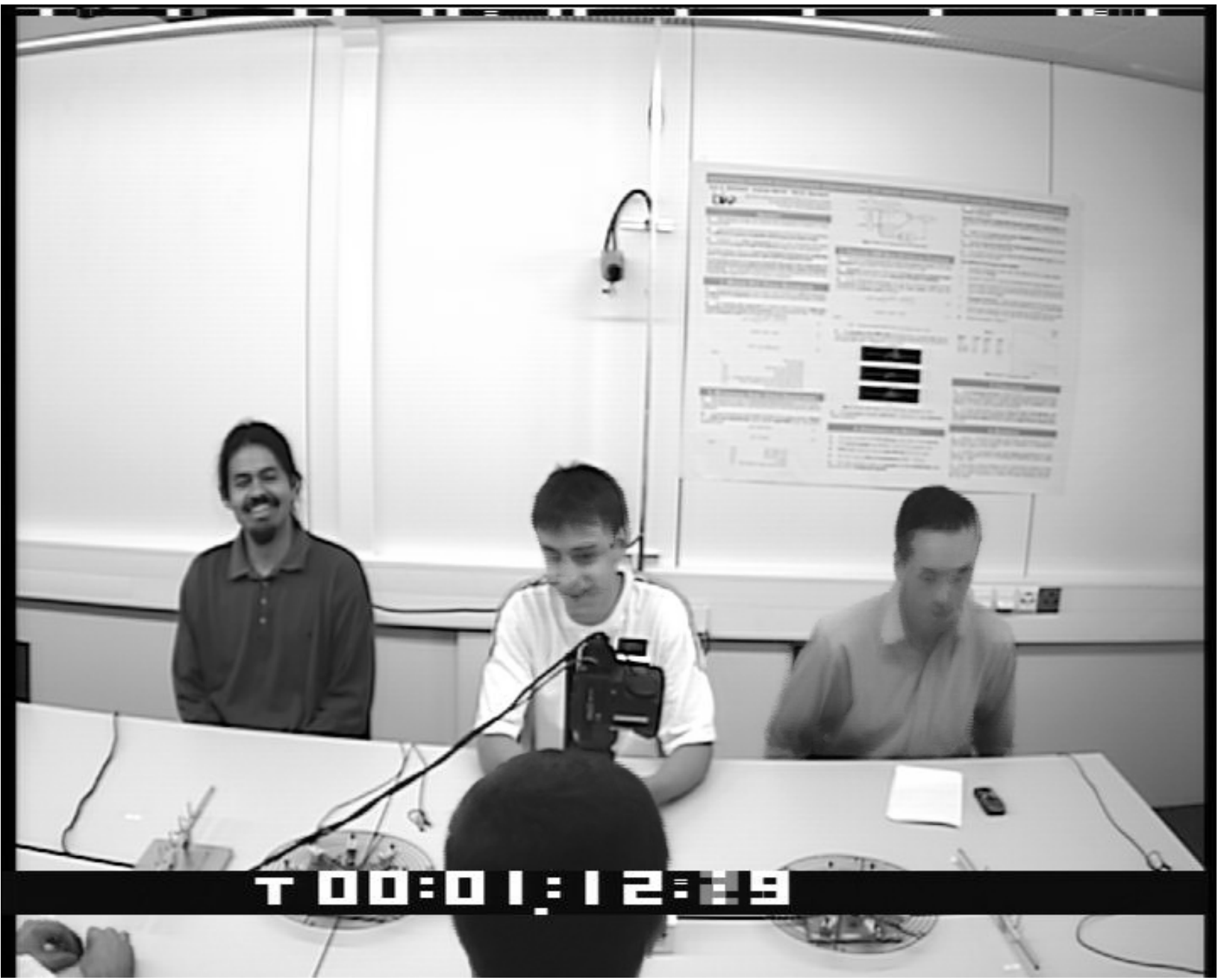}
&\includegraphics[width=0.18\textwidth,
keepaspectratio]{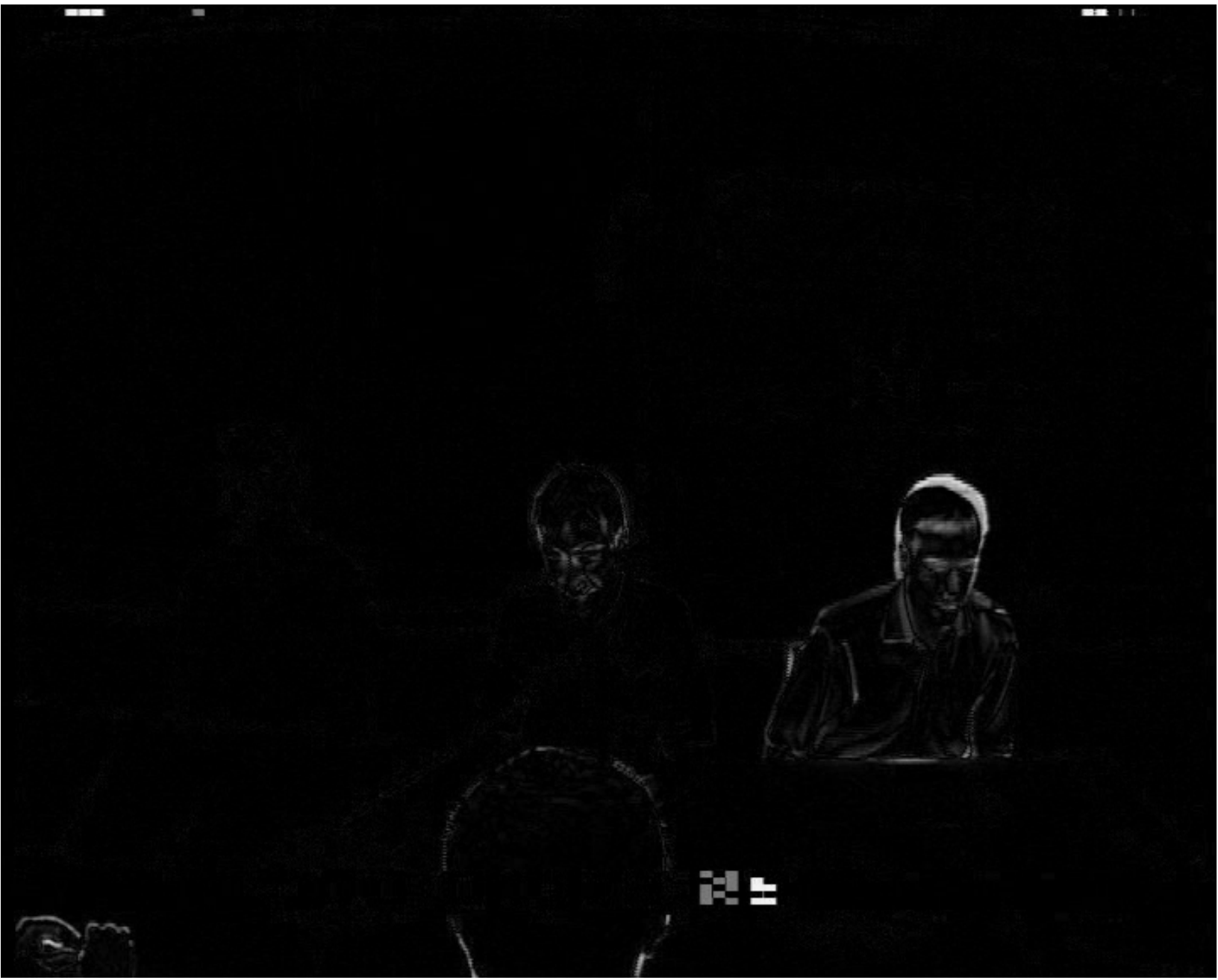}
&\includegraphics[width=0.18\textwidth,
keepaspectratio]{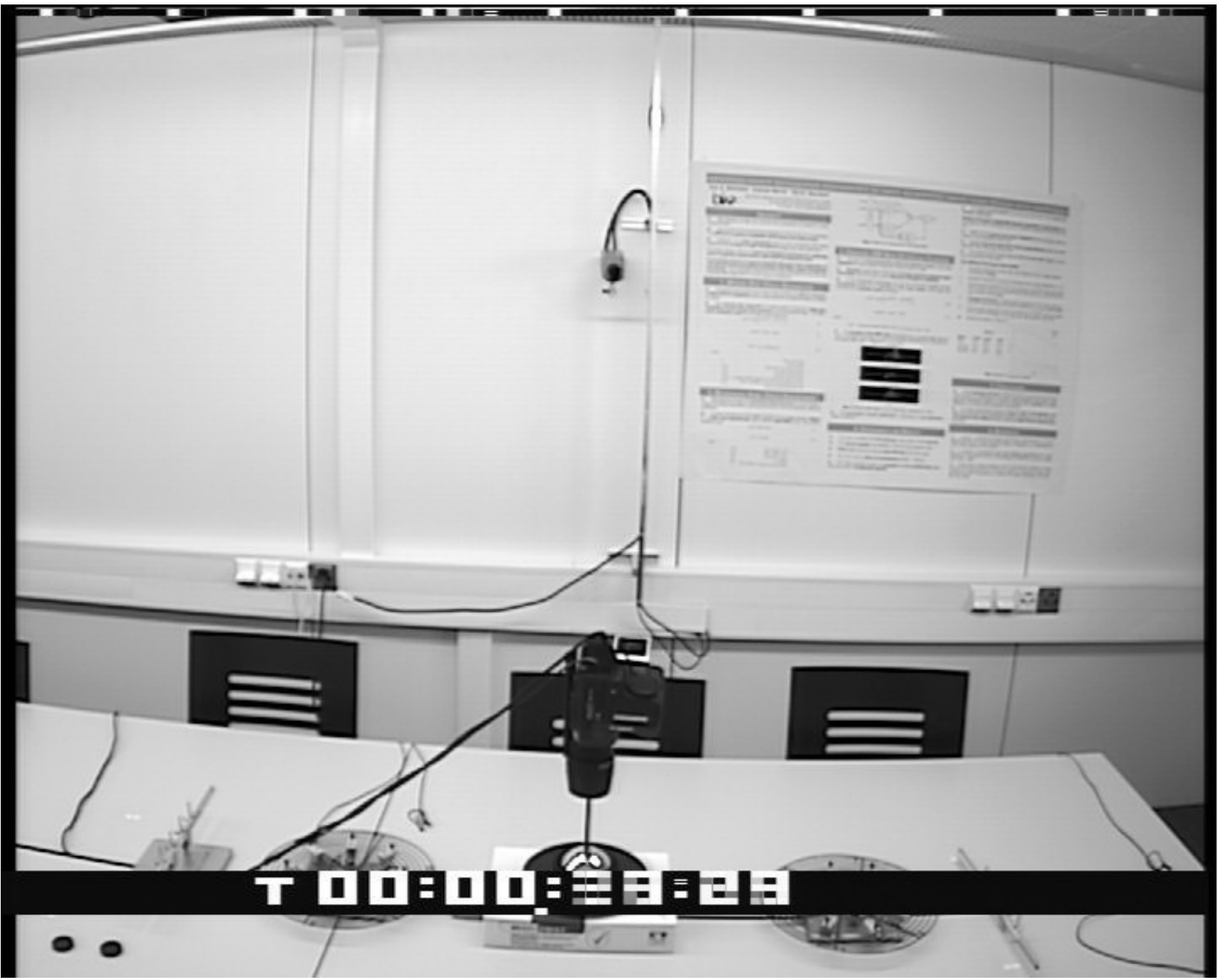}
&\includegraphics[width=0.18\textwidth,
keepaspectratio]{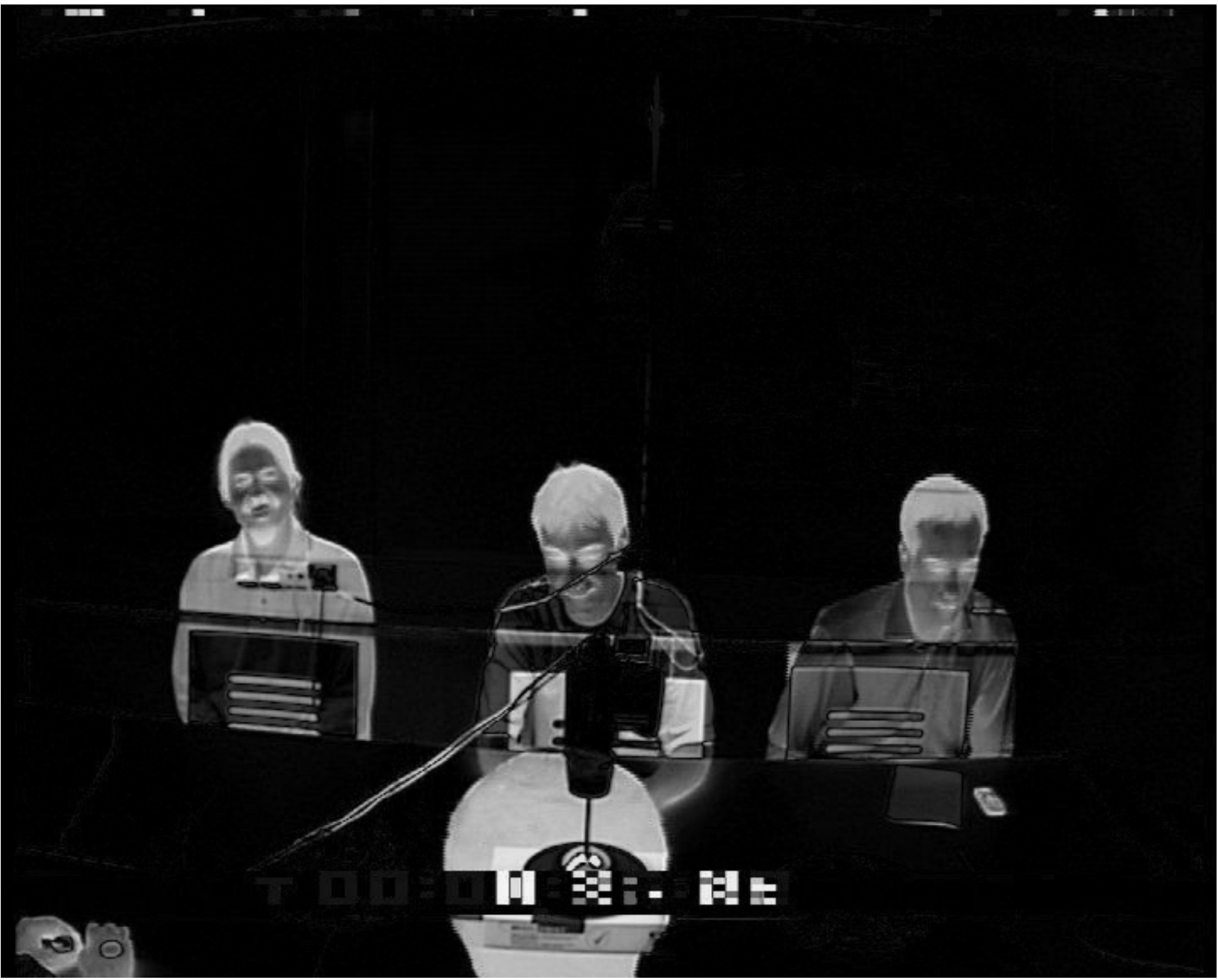}\\
(a) Original & (b) Background (Median) & (c) Foreground (Median) & (d) Background ($l_1$) & (e) Foreground ($l_1$)
\end{tabular}
\caption{The partial background modeling results of median filter
and $l_1$ filtering on the ``meeting'' video sequence. (b)-(c) and
(d)-(e) are the the background ($\mathbf{L}^*$) and the foreground
($\mathbf{S}^*$) recovered by median filter and $l_1$ filtering,
respectively.}\label{fig:vs}
\end{figure*}

To further test the performance of $l_1$ filtering on large scale
data set, we also collect a video sequence (named ``meeting'') of
700 frames, each of which has a resolution $576\times 720$. So the
data matrix is of size greater than $700\times 400000$, which
cannot be fit into the memory of our PC. As a result, we cannot
use the standard ADM to solve the corresponding PCP problem. As
for RSL, we have found that it did not converge on this data. Thus
we only compare the performance of median filter and $l_1$
filtering. The time cost is reported in Table~\ref{tab:vs} and the
qualitative comparison is shown in Figure~\ref{fig:vs}. We can see
that $l_1$ filtering is as fast as median filtering with 20
frames, and median filter fails on this data set. This is because
the mechanism of median filter is based on the (local) frame
difference. Thus when the scene contains slowly moving objects
(such as that in the ``meeting'' video sequence), median filter
will not give good results. In contrast, the background and the
foreground can be separated satisfactorily by $l_1$ filtering.
This makes sense because our $l_1$ filtering can exactly recover
the (global) low-rank structure for the background and remove the
foreground as sparse errors.

\section{Conclusion and Further Work}\label{sec:con}
In this paper, we propose the first \emph{linear time} algorithm,
named the $l_1$ filtering method, for \emph{exactly} solving very
large PCP problems, whose ranks are supposed to be very small
compared to the data size. It first recovers a seed matrix and
then uses the seed matrix to filter some rows and columns of the
data matrix. It avoids SVD on the original data matrix, and the
$l_1$ filtering step can be done in full parallelism. As a result,
the time cost of our $l_1$ filtering method is only linear with
respect to the data size, making applications of RPCA to extremely
large scale problems possible. The experiments on both synthetic
and real world data demonstrate the high accuracy and efficiency
of our method. It is possible that the proposed technique can be
applied to other large scale nuclear norm minimization problems,
e.g., matrix completion (\cite{cai2008singular}) and low-rank
representation (\cite{Liu2010LRR}). This will be our future work.

\begin{acknowledgements}
The authors would like to thank Prof. Zaiwen Wen and Dr. Yadong Mu
for sharing us their codes for LMaFit (\cite{shen-MatrixFac}) and
random projection (\cite{mu-2011-RandomProj}), respectively. This
work is partially supported by the grants of the National Nature
Science Foundation of China-Guangdong Joint Fund (No. U0935004),
the National Nature Science Foundation of China Fund (No.
60873181, 61173103) and the Fundamental Research Funds for the
Central Universities. The first author would also like to thank
the support from China Scholarship Council.
\end{acknowledgements}

\bibliographystyle{spbasic}      

\begin{thebibliography}{}

\bibitem[{Aanes et~al(2002)Aanes, Fisker, Astrom, and
  Carstensen}]{Aanes-2002-fac}
Aanes H, Fisker R, Astrom K, Carstensen J (2002) Robust factorization. IEEE
  Trans on PAMI 24(9):359--368

\bibitem[{Baccini et~al(1996)Baccini, Besse, and {de
  Falguerolles}}]{Baccini96anl1-norm}
Baccini A, Besse P, {de Falguerolles} A (1996) An $l_1$-norm {PCA} and a
  heuristic approach. In: Proceedings of the International Conference on
  Ordinal and Symbolic Data Analysis, pp 359--368

\bibitem[{Benedek and Szir\'anyi(2008)}]{Benedek2008VS}
Benedek C, Szir\'anyi T (2008) Bayesian foreground and shadow detection in
  uncertain frame rate surveillance videos. IEEE Trans on Image Processing
  17(4):608--621

\bibitem[{Cai et~al(2010)Cai, Cand\'es, and Shen}]{cai2008singular}
Cai J, Cand\'es E, Shen Z (2010) A singular value thresholding algorithm for
  matrix completion. SIAM Journal on Optimization 20(4):1956--1982

\bibitem[{Cand\'es and Wakin(2007)}]{candes2007intro}
Cand\'es E, Wakin M (2007) An introduction to compressive sampling. IEEE Signal
  Processing Magazine 25(2):21--30

\bibitem[{Cand\'es et~al(2011)Cand\'es, Li, Ma, and Wright}]{candes2009robust}
Cand\'es E, Li X, Ma Y, Wright J (2011) Robust {P}rincipal {C}omponent
  {A}nalysis? Journal of the ACM 58(3):11

\bibitem[{{De la Torre} and Black(2003)}]{fernando2003rsl}
{De la Torre} F, Black M (2003) A framework for robust subspace learning. IJCV
  54(1--3):117--142

\bibitem[{Drineas et~al(2006)Drineas, Kannan, and Mahoney}]{Drineas2006LTSVD}
Drineas P, Kannan R, Mahoney M (2006) Fast {M}onte {C}arlo algorithms for
  matrices {II}: Computing a low rank approximation to a matrix. SIAM Journal
  on Computing 36(1):158--183

\bibitem[{Ganesh et~al(2009)Ganesh, Lin, Wright, Wu, Chen, and
  Ma}]{Ganesh-2009}
Ganesh A, Lin Z, Wright J, Wu L, Chen M, Ma Y (2009) Fast algorithms for
  recovering a corrupted low-rank matrix. In: Proceedings of International
  Workshop on Computational Advances in Multi-Sensor Adaptive Processing

\bibitem[{Halko et~al(2011)Halko, Martinsson, and Tropp}]{Halko-2011-random}
Halko N, Martinsson P, Tropp J (2011) Finding structure with randomness:
  Probabilistic algorithms for constructing approximate matrix decompositions.
  SIAM Review 53(2):217--288

\bibitem[{Ji et~al(2010)Ji, Liu, Shen, and Xu}]{Hui2010Video}
Ji H, Liu C, Shen Z, Xu Y (2010) Robust video denoising using low-rank matrix
  completion. In: CVPR

\bibitem[{Ke and Kanade(2005)}]{Ke2005l1}
Ke Q, Kanade T (2005) Robust $l_1$-norm factorization in the presence of
  outliers and missing data by alternative convex programming. In: CVPR

\bibitem[{Larsen(1998)}]{Larsen-1998-PROPACK}
Larsen R (1998) Lanczos bidiagonalization with partial reorthogonalization.
  Department of Computer Science, Aarhus University, Technical report, DAIMI
  PB-357

\bibitem[{Lin et~al(2009)Lin, Chen, Wu, and Ma}]{lin2009augmented}
Lin Z, Chen M, Wu L, Ma Y (2009) The augmented {L}agrange multiplier method for
  exact recovery of corrupted low-rank matrices. UIUC Technical Report
  UILU-ENG-09-2215

\bibitem[{Liu et~al(2010)Liu, Lin, and Yu}]{Liu2010LRR}
Liu G, Lin Z, Yu Y (2010) Robust subspace segmentation by low-rank
  representation. In: ICML

\bibitem[{Mu et~al(2011)Mu, Dong, Yuan, and Yan}]{mu-2011-RandomProj}
Mu Y, Dong J, Yuan X, Yan S (2011) Accelerated low-rank visual recovery by
  random projection. In: CVPR

\bibitem[{Nie et~al(2011)Nie, Huang, Ding, Luo, and Wang}]{Nie2011L1}
Nie F, Huang H, Ding C, Luo D, Wang H (2011) Robust principal component
  analysis with non-greedy $l_1$-norm maximization. In: IJCAI

\bibitem[{Peng et~al(2010)Peng, Ganesh, Wright, Xu, and Ma}]{peng2010rasl}
Peng Y, Ganesh A, Wright J, Xu W, Ma Y (2010) {RASL}: {R}obust alignment by
  sparse and low-rank decomposition for linearly correlated images. In: CVPR

\bibitem[{Rao et~al(2010)Rao, Tron, Vidal, and Ma}]{Rao-2010-motion}
Rao S, Tron R, Vidal R, Ma Y (2010) Motion segmentation in the presence of
  outlying, incomplete, and corrupted trajectories. IEEE Trans on PAMI
  32(10):1832--1845

\bibitem[{Shen et~al(2011)Shen, Wen, and Zhang}]{shen-MatrixFac}
Shen Y, Wen Z, Zhang Y (2011) Augmented {L}agrangian alternating direction
  method for matrix separation based on low-rank factorization. preprint

\bibitem[{Skocaj et~al(2007)Skocaj, Leonardis, and
  Bischof}]{Skocaj-2007-weight}
Skocaj D, Leonardis A, Bischof H (2007) Weighted and robust learning of
  subspace representations. Pattern Recognition 40(5):1556--1569

\bibitem[{Storer et~al(2009)Storer, Roth, Urschler, and
  Bischof}]{Markus2009FRPCA}
Storer M, Roth P, Urschler M, Bischof H (2009) Fast-robust {PCA}. In: Proc.
  16th Scandinavian Conference on Image Analysis (SCIA)

\bibitem[{Wang et~al(2009)Wang, Dong, Tong, Lin, and Guo}]{Wang2009Nystrom}
Wang J, Dong Y, Tong X, Lin Z, Guo B (2009) Kernel {N}ystr{\"o}m method for
  light transport. ACM Transactions on Graphics 28(3)

\bibitem[{Wright et~al(2009)Wright, Ganesh, Rao, Peng, and
  Ma}]{Wright-NIPS2009}
Wright J, Ganesh A, Rao S, Peng Y, Ma Y (2009) Robust principal component
  analysis: Exact recovery of corrupted low-rank matrices via convex
  optimization. In: NIPS

\bibitem[{Wu et~al(2010)Wu, Ganesh, Shi, Matsushita, Wang, and Ma}]{wu-robust}
Wu L, Ganesh A, Shi B, Matsushita Y, Wang Y, Ma Y (2010) Robust photometric
  stereo via low-rank matrix completion and recovery. preprint

\bibitem[{Yuan and Yang(2009)}]{yuan2009sparse}
Yuan X, Yang J (2009) Sparse and low-rank matrix decomposition via alternating
  direction methods. preprint

\bibitem[{Zhang et~al(2012)Zhang, Ganesh, Liang, and Ma}]{zhang2011texture}
Zhang Z, Ganesh A, Liang X, Ma Y (2012) {TILT}: Transform-invariant low-rank
  textures. accepted by IJCV

\end{thebibliography}


\end{document}